 \documentclass[preprint,3p]{elsarticle}
\usepackage{amsmath}
\usepackage{amssymb}
\usepackage{tabularx}
\usepackage{bm,color,soul}
\input psfig.sty

\usepackage{epsfig,epstopdf,color}

\usepackage{hyperref}

\newcommand{\vep}{\varepsilon}
\newcommand{\gm}{\gamma}

\newcommand{\tx}{{\widetilde x} }
\newcommand{\ty}{{\widetilde y} }
\newcommand{\tz}{{\widetilde z} }

\newcommand\bx{\mathbf{x}}
\newcommand\bz{\mathbf{z}}
\newcommand\bn{\mathbf{n}}
\newcommand\bbm{\mathbf{m}}

\newcommand\bk{\mathbf{k}}

\newcommand\by{\mathbf{y}}
\newcommand\bau{{\textbf B}_1}

\newcommand{\p}{\partial}

\newcommand{\Og}{\Omega}
\newcommand{\fl}[2]{\frac{#1}{#2}}

\newcommand{\nn}{\nonumber}

\newcommand{\bt}{\beta}

\newcommand{\Dt}{\Delta}
\newcommand{\tbx}{\widetilde{\bf x} }
\newcommand{\tby}{\widetilde{\bf y} }

\renewcommand{\theequation}{\arabic{section}.\arabic{equation}}
\newcommand{\be}{\begin{equation}}
\newcommand{\ee}{\end{equation}}
\newcommand{\ba}{\begin{array}}
\newcommand{\ea}{\end{array}}
\def\bea{\begin{eqnarray}}
\def\eea{\end{eqnarray}}
\def \beas{\begin{eqnarray*}}
\def \eeas{\end{eqnarray*}}
\newtheorem{remark}{Remark}[section]
\newtheorem{exmp}{Example}[section]
\newtheorem{thm}{Theorem}[section]
\newtheorem{lemma}{Lemma}[section]

\begin{document}

\begin{frontmatter}

\title{On  the ground states and dynamics of space fractional  nonlinear \\ Schr\"{o}dinger/Gross-Pitaevskii
  equations with rotation \\ term  and  nonlocal nonlinear interactions}
\author[1]{Xavier ANTOINE}
\ead{xavier.antoine@univ-lorraine.fr}
\ead[url]{http://iecl.univ-lorraine.fr/\~{}xantoine/}

\address[1]{Institut Elie Cartan de Lorraine, Universit\'e de Lorraine, Inria Nancy-Grand Est,
F-54506 Vandoeuvre-l\`es-Nancy Cedex, France.}

\author[1,2]{Qinglin TANG\corref{5}}
\ead{qinglin.tang@inria.fr}

\address[2]{Laboratoire de Math\'ematiques Rapha\"el Salem, Universit\'e de Rouen,
Technop\^{o}le du Madrillet, 76801 Saint-Etienne-du-Rouvray, France.}
\cortext[5]{Corresponding author.}

\author[3,4]{Yong ZHANG}
 \ead{yong.zhang@univ-rennes1.fr}

\address[3]{Universit\'e de Rennes 1, IRMAR, Campus de Beaulieu, 35042 Rennes Cedex, France}
\address[4]{Wolfgang Pauli Institute c/o Fak. Mathematik, University Wien, Oskar-Morgenstern-Platz 1, 1090 Vienna, Austria}


\begin{abstract}
In this paper, we propose  some efficient and robust numerical methods
to compute the ground states and dynamics of  Fractional   
Schr\"{o}dinger Equation  (FSE) with a rotation
term and nonlocal nonlinear interactions. In particular,  a newly developed 
Gaussian-sum (GauSum) solver is used for the nonlocal interaction evaluation \cite{EMZ2015}.  
To compute the ground states, we integrate the preconditioned Krylov subspace pseudo-spectral method \cite{AD1} and the 
GauSum solver. For the dynamics simulation,
using the rotating Lagrangian coordinates transform \cite{BMTZ2013}, 
we first reformulate the FSE into a new equation without rotation. 
Then, a time-splitting pseudo-spectral scheme incorporated with the 
GauSum solver is proposed to simulate the new FSE.  
In parallel to the numerical schemes, we also prove some  
 existence and nonexistence results for the ground states. Dynamical laws 
 of some standard quantities, including the mass, energy, angular momentum and the center of mass, 
are stated. The ground states properties with respect to the fractional order
and/or rotating frequencies, dynamics involving decoherence and turbulence together with some interesting phenomena are reported.

\end{abstract}

\begin{keyword}
fractional Schr\"{o}dinger equation, rotation, nonlocal nonlinear interaction, 
rotating Lagrangian coordinates, Gaussian-sum solver, ground state, dynamics

\end{keyword}

\end{frontmatter}

\tableofcontents


\section{Introduction}
\noindent

Recently, a great deal of attention has been directed towards the derivation of a powerful
generalization of PDEs through
the inclusion of fractional order operators. These developments now
impact strongly most areas of physics and engineering
 \cite{GuoPuHuangBook,HLSBookChapter,SpecialIssueJCP,Pod1999,WBG2003}. Additionally, some
new applications are also emerging  in biology, molecular dynamics, finance, etc. Due to the fact
that extremely important applications are related to these models,
a  significative  effort has been made in the last few years to obtain some mathematical properties  
and numerical tools \cite{GuoPuHuangBook} Êfor the generalized systems of PDEs. An example of such
a keen interest is the recent
Journal of Computational Physics \cite{SpecialIssueJCP}  special issue  in 2015 that is
dedicated to
``Fractional PDEs: Theory, Numerics, and Applications''.
The aim of this paper is to contribute to this new hot area for fractional quantum physics, with possible applications, e.g.
in Bose-Einstein condensation (BEC).

During the last decades, the classical Schr\"{o}dinger Equation (SE) has been widely investigated and
applied to many areas in physics (optics, electromagnetic, superfluidity, etc.). It is known as the fundamental equation of classical 
quantum mechanics which can be interpreted by the Feynman path integral approach over Brownian-like quantum 
paths \cite{FH1965}.  Brownian motion (Wiener process) represents a simple diffusion random walk process.
More general and complex stochastic  processes (L\'{e}vy motion) exist and  can be still presumably modeled 
 by modifying the standard diffusion equation using a fractional Laplacian operator $(-\Delta)^{s}:=(-\nabla^{2})^{s}$ 
 (with $s> 0$ being the fractional order) \cite{MK2000, SKZ1999}. 
 L\'{e}vy processes provide a general framework to study \textit{anomalous diffusion}.   
    Fractional diffusion has been widely studied by many authors 
   \cite{Bouchard,DGLZ2012, HLSBookChapter,MK2000,MetzlerKlafter,Pod1999, Vlahos, WBG2003}
   and is now  considered as a suitable way to describe spatially disordered systems (such as porous media and fractal media),
   turbulent fluids and plasmas, biological media with traps, etc.
Fractional diffusion for classical mechanics is referred to as \textit{subdiffusion}
if $s<1$  and as \textit{superdiffusion} if $s>1$.
     For more details, we refer  to \cite{HLSBookChapter}
where a concise table shows the scaling laws for fractional diffusion.

Analogously, in the context of quantum mechanics, fractional quantum models, based in particular on
Schr\"{o}dinger-type equations,  are now emerging while being however more limited in terms of publications and studies
compared with classical fractional mechanics.
 Laskin extended the Feynman path integral approach
over L\'{e}vy-like quantum  paths and derived a Fractional Schr\"{o}dinger Equation (FSE) that modifies the SE
by involving the fractional Laplacian $(-\Delta)^{s}$ \cite{Laskin0, Laskin3,Laskin4,Laskin5}. The FSE
was applied  to represent the Bohr atom, fractional oscillator \cite{Laskin5}, and it is a new fractional approach to study the quantum chromodynamics (QCD) 
problem of quarkonium \cite{Laskin0}. The FSE also arises in the continuum limit of the discrete SE 
with long-range dispersive interaction \cite{KLS2012}, in the mathematical
description of  boson stars \cite{ES2007} and in some models of water wave dynamics \cite{IP2014}. 
 It has also  been  proposed to study  BEC
of which the particles obey  a non-Gaussian distribution law \cite{ESDB2012,UB2013,UHTA2012}, where 
FSE was named as Fractional Gross-Pitaevskii Equation (FGPE) and BEC as Fractional BEC (FBEC).
Compared with the SE, the literature on FSE is quite limited but growing quickly to
understand its mathematical and physical properties.


More precisely, we consider here  the following generalized dimensionless  
 (space-)Fractional NonLinear Schr\"{o}dinger equation (FNLSE) with a rotation term and 
a  nonlocal nonlinear interaction
\bea
\label{SFSE}
&&i\p_t\psi(\bx,t)=\left[\fl{1}{2}\left(-\nabla^2+m^2\right)^{s}
+V(\bx)+\beta|\psi(\bx,t)|^2+\lambda\Phi(\bx,t)
-\Omega L_z
\right]\psi(\bx,t),\\
\label{NonLocalPot}
&& \Phi(\bx,t)=\mathcal{U}\ast|\psi(\bx,t)|^2,
\qquad \bx\in{\Bbb R}^d,\ t>0, \  d\ge2.
\eea
In the context of BEC, this equation is also called as FGPE.
Here,  $\psi(\bx,t)$ is the complex-valued wave-function, 
$t>0$ is the time variable and $\bx\in{\Bbb R}^d$ is the spatial coordinate.
The constant $m\ge0$  denotes 
the scaled particle mass, with $m=0$ representing the massless particle. The parameter
$s>0$ is the space fractional order characterizing the nonlocal dispersive interaction.
The fractional kinetic operator is defined \textit{via} a Fourier integral operator
%
\bea
\left(-\nabla^2+m^2\right)^{s} \psi =\frac{1}{(2\pi)^d} \int_{{\Bbb R}^d} \widehat{\psi}(\bk) \, (|\bk|^2+m^2)^s e^{i \bk\cdot\bx}d\bk,
\eea
where the Fourier transform is given by $ \widehat{\psi}(\bk) = \int_{{\Bbb R}^d}  \psi(\bx) e^{-i\bk \cdot \bx}  d \bx$.
The  potential $V(\bx)$ is supposed to be trapping, a standard example is the harmonic potential  given by 
\be
\label{harm_poten}
V(\bx)=
\left\{
\begin{array}{cr}
\fl{\gm_x^2 x^2+ \gm_y^2 y^2}{2}, &  \qquad d=2, \\ [0.5em]
\fl{\gm_x^2 x^2+ \gm_y^2 y^2 +\gm_z^2 z^2}{2},   & \qquad  d=3,
\end{array}
\right.
\ee
where $\gm_v$ ($v=x, y, z$) is the trapping frequency in the $v$-direction. The real-valued constants
$\beta$  and $\lambda$   characterize  the local and nonlocal
interaction strengths (positive/negative for repulsive/attractive interaction), respectively. 
The local interaction is supposed to be cubic, but other choices may also be considered.
Concerning the nonlocal interaction (\ref{NonLocalPot}), 
the convolution kernel  $\mathcal{U}(\bx)$ can be chosen as either the kernel of a Coulomb-type interaction or
a Dipole-Dipole Interaction (DDI) \cite{BJTZ2015,BTZ2015,CHHO2013}
\be\label{Kernel}
\mathcal{U}(\bx)=
\left\{
\begin{array}{ll}
\fl{1}{2^{d-1}\pi |\bx|^{\mu}}, \quad 0<\mu \le d-1, &	 \\[0.5em]
-\delta(\bx)-3\,\partial_{\bn\bn}   
\left( \fl{1}{4\pi|\bx|} \right), &     \\[0.5em]
-\fl{3}{2} \left(\partial_{\bn_{\perp}\bn_{\perp}}-n_3^2 \nabla_{\perp}^2\right)\left( \fl{1}{2\pi|\bx|} \right),
 &    
\end{array}
\right.
\Longleftrightarrow
\quad
\widehat{\mathcal{U}}(\bk)=
\left\{
\begin{array}{ll}
 \fl{C}{|\bk|^{d-\mu}}, \quad 0<\mu \le d-1, &	 {\rm Coulomb},\\[0.5em]
-1+\fl{3 (\bn\cdot\bk)^2}{|\bk|^2}, &   {\rm 3D \ DDI}, \\[0.5em]
\fl{3[(\bn_\perp\cdot\bk)^2-n_3^2 |\bk|^2]}{2 |\bk|},  &    {\rm 2D \ DDI},
\end{array}
\right.
\ee  
where $C = \pi^{d/2-1} 2^{1-\mu} \Gamma(\frac{d-\mu}{2})/\Gamma(\frac{\mu}{2})$ 
($\Gamma (t) := \int_0^\infty x^{t-1}e^{-x} d x$ is the Gamma function),
$\bn=(n_1, n_2, n_3)^T\in \mathbb R^3$ is a unit vector representing the dipole orientation
and $\bn_\perp=(n_1, n_2)^T$.
In addition, $L_z = -i (x\partial_y -y\partial_x)= -i \partial_\theta$ is the $z$-component of the 
angular momentum, $\Omega$  represents the rotating frequency. 
%

The FNLSE conserves two important physical quantities (see Section \ref{sec: dyn_prop}):  the {\sl mass} 
\be\label{mass}
\mathcal{N}(\psi(\cdot, t)):=\mathcal{N}(t) := \int_{\mathbb{R}^d} |\psi(\bx,t)|^2 d\bx \equiv \mathcal{N}(0),
\ee
and the {\sl energy}
\be\label{energy}
\mathcal{E}(\psi(\cdot, t))=: \mathcal{E}(t)=\int_{\mathbb{R}^d} \Big[\fl{1}{2} \bar{\psi}\big(-\nabla^2+m^2\big)^{s}\psi+V(\bx) |\psi|^2
+\fl{\beta}{2} |\psi|^4 +\fl{\lambda}{2} \Phi |\psi|^2 
-\Omega \bar{\psi} L_z \psi  
\Big]
\equiv \mathcal{E}(0).
\ee
Here, $\bar{\psi}$ is the complex conjugate of  $\psi$.
The ground states $\phi_g(\bx)$ of the FNLSE (\ref{SFSE})  are defined by
\be
\label{gs_def}
\phi_g(\bx)={\rm arg}\min_{\phi\in S}\mathcal{E}(\phi),  \qquad S=\{\phi\in \mathbb{C} |\; \|\phi\|_{2}=1,   \mathcal{E}(\phi)<\infty\},
\ee
where $\|\phi\|_{2}$ is the $L^2(\mathbb{R}^{d})$-norm of $\phi$.

%

The FNLSE \eqref{SFSE} brings together a wide range of Schr\"{o}dinger-type PDEs.
 When $s=1$ and $m=0$, FNLSE reduces to the standard nonlinear 
 Schr\"{o}dinger equation (NLSE, also known as GPE). 
 Both the ground states and dynamics  properties of NLSE have been extensively  studied theoretically and
   numerically. One can refer e.g. to  \cite{AABES2008, ABB2013, AD2,AD3, ADBookChapter, BC2013,   BJMZ2013, BJTZ2015, BMTZ2013, BTZ2015, 
       Besse2004,CMS2008,DK2010} and references therein.  
 For $s=1/2$ and $\Phi$ taken as the Coulomb potential, \eqref{SFSE} reduces to 
the semi-relativistic Hartree equation that models boson stars \cite{ES2007}.  Properties of the ground states have been partially investigated
 in \cite{ CN2011, Len2009, LY1987} for $V(\bx)\equiv 0$. The Cauchy problem of  generalized semi-relativistic 
 Hartree equation (with $s\in [\fl{1}{2},1]$)  has also been widely studied
 in \cite{AMS2008,CHHO2013,CO2006, ES2007,FL2007,Len2007}. 
%
%
To the best of our knowledge, there are neither theoretical nor numerical studies on the ground state properties for  $s>0$ other than the
cases $s=1/2$ and $1$. When $s\in(0,1)$   (which would correspond to a \textit{subdispersion} effect in analogy to the subdiffusion process 
characterizing heat-like equations \cite{HLSBookChapter}) and $m=\lambda=\Og=0$, it reduces to the FNLSE that is originally derived by Laskin \cite{Laskin0}. 
Later, he proved the hermitian character of the fractional Hamiltonian, derived the energy spectra of a hydrogen-like atom and computed 
a fractional oscillator \cite{Laskin0,Laskin4,Laskin5}.
Since then, the FNLSE has attracted an increasing attention. 
 For example, for stationary FNLSE with bounded/unbounded potential and various generalised nonlinearities other than 
$|\psi|^2\psi$, the existence of solutions (such as the bound/ground state solutions and radially symmetric solutions), 
and their corresponding properties have been investigated. Moreover, 
the global and/or local well-posedness for the initial value problem (\ref{SFSE})-(\ref{NonLocalPot}) with
$V(\bx)\equiv 0$  and $\lambda=0$ were also studied. We refer to 
\cite{CHHO2014, Feng2013, GH2011, HS2015, Secchi2013, Secchi2014, SS2014, ShangZ2014} and references therein for more details. 
%
%
For $s>1$ (that we call \textit{superdispersion} hereafter),  there are a few
Schr\"{o}dinger-type equations, while it is quite common for the superdiffusion equations \cite{HLSBookChapter}.
We consider here this case for some possible eventual  physical applications.

Generally speaking, it is  difficult to obtain analytical solutions of the FNLSE due to the nonlocal fractional dispersive interaction. 
For example, even for the simplest case with a box potential,  
there is still a controversy over the eigenpair solutions \cite{Bayin2012,Bayin2013,GX2006,HS2013,JXHS2010} .
Therefore, being able to develop some  accurate numerical methods is crucial and would
 provide
a powerful tool to understand fractional quantum mechanics in view of applications.  
Nevertheless, there are few numerical studies so far. Amore \textit{et al.} \cite{AFHS2010} proposed a collocation method and
Wang \textit{et al.} \cite{WH2015} developed an  energy conserving
  Crank-Nicolson  finite difference (FD) scheme when $\Og=\lambda=0$.
Similar FD  schemes were also proposed for coupled equations \cite{WXY2013, WXY2014}.  
As  is well-known, the Crank-Nicolson scheme is nonlinearly implicit and hence requires heavy inner iterations. 
 Worse still, the nonlocal nature of the fractional Laplacian naturally leads to  dense matrix representation
 that hinders efficient computations.  
Recently, the time-splitting Fourier pseudo-spectral method was adapted 
to study the dynamics when $\Og=\lambda=0$ \cite{KZ2014, KSM2014}.  Decoherence properties and  
 finite time blow-up results were studied respectively in \cite{KZ2014} and  \cite{KSM2014}.  
 When the nonlocal nonlinear interaction ($\lambda\ne0$) is taken into account, 
Bao and Dong \cite{BD2011}  proposed a sine pseudo-spectral method
to compute the ground states and dynamics of the three-dimensional semi-relativistic Hartree equation ( $\mu=1$ in (\ref{Kernel})).
In \cite{BD2011}, the Coulomb potential  $\Phi$ (\ref{NonLocalPot}) is reformulated to satisfy the following Poisson equation
\be
\label{3d-poisson}
-\triangle \Phi=|\psi|^2,  \; \; \bx \in \mathbb R^3,\quad \mbox{ with } \quad  \lim_{|\bx|\rightarrow\infty}\Phi(\bx)=0.
\ee
Similar ideas were also applied to nonlocal DDI in NLSE \cite{ BC2013, BCW2010}.
However, due to the slow decay property of $\Phi$ at the far-field, a quite large computational 
domain is necessary to guarantee a satisfactory accuracy.   
Up to now, most existing numerical methods are proposed for non-rotating FNLSE with $s\le1$. 
As far as we know, there were neither theoretical nor numerical methods for the generalized FNLSE (\ref{SFSE}) 
for  both subdispersion $s<1$ and superdispersion  $s>1$, with $\lambda\Og\ne0$. 
The difficulties to develop an accurate and efficient scheme lie in the evaluation of 
the nonlocal interaction $\Phi$ \eqref{NonLocalPot} and 
proper treatment of the rotation term $L_z\psi$.

To compute the nonlocal interaction, Jiang \textit{et al.} \cite{JGB2014} recently  
proposed an accurate NonUniform Fast Fourier Transform (NUFFT)-based algorithm in the Fourier domain by adopting 
the polar/spherical coordinates near the singularity. The method requires $O(N\log N)$ arithmetic operations ($N$ being
 the total number of  grid points)  and  is more accurate than the standard PDE approach (\ref{3d-poisson}). This solver has
been recently integrated within the gradient flow algorithm and time-splitting scheme for computing the ground state
and dynamics of NLSE \cite{BJTZ2015, BTZ2015}.
However, this solver is not ideal because of the large pre-factor in $O(N\log N)$, and it is rather slow for 3D problems.  
Very recently,  by approximating the kernel $\mathcal{U}(\bx)$  with the summation of a finite number of Gaussians, 
Zhang \textit{et al.} \cite{EMZ2015} proposed a Gaussian-sum (GauSum)-based method to evaluate $\Phi$ in the physical space.  
The algorithm also achieves a spectral accuracy, requires $O(N\log N)$ operations and obtains a speed-up factor 
around 3-5 compared with the NUFFT-based algorithm.
Concerning the rotation term, Antoine and Duboscq \cite{AD1,ADBookChapter} proposed 
a robust preconditioned Krylov subspace spectral solver for the ground state computation of the NLSE 
with  large $\Og$ and $\beta$. For the dynamics of the NLSE with a rotation term, Bao \textit{et al.}  \cite{BMTZ2013}
developed  a rotating Lagrangian coordinates transformation method  to reformulate the rotating term into 
a time-dependent trapping potential in the rotating Lagrangian coordinates, which allows for
the implementation of high-order time-splitting schemes for the new NLSE \cite{Besse2015}.  

\

The main objectives of this paper are  threefold.
\begin{enumerate}
\item Investigate theoretically the existence of the ground states of the general FNLSE (\ref{SFSE}) 
with respect to the fractional order $s$ and the rotation speed $\Og$. Develop the dynamical laws for the centre of mass as
well as other standard dynamical quantities for general  $s$ 
and arbitrary $\Og$, 
and compare them with the ones derived in \cite{KZ2014}.
\item Develop some efficient and accurate numerical methods for computing the ground states and 
dynamics of the general FNLSE (\ref{SFSE}) by incorporating the 
GauSum solver into the adapted version of the 
gradient flow and time-splitting Fourier pseudo-spectral method.
The preconditioned Krylov subspace iteration  \cite{AD1}   
and the rotating Lagrangian transformation technique \cite{BMTZ2013} will be also integrated into the numerical methods 
for the ground state computation and dynamics simulation, respectively.
\item Apply our numerical methods to study some interesting behavior, such as  the influence of 
the nonlocal dispersion on the ground states and the vortex pattern as well as possible 
dynamical properties such as chaos and decoherence.
\end{enumerate}
 The rest of the paper is organized as follows. 
 In Section 2, we  briefly review the Gaussian Sum method.  The ground state computation, including the ground states properties and
 numerical methods as well as numerical results are presented in Section 3. In Section 4, we derive some 
 dynamical laws for some global physical quantities that are usually considered for the standard NLSE.  We then
  propose an efficient and robust numerical method for the dynamics simulation. Some numerical results
 are also reported. Finally,  a conclusion and some discussions are developed  in Section 5.

\section{Brief review of the Gaussian-Sum (GauSum) method}
\setcounter{equation}{0}
\setcounter{table}{0}
\setcounter{figure}{0}

With the strong confining potential, the density is smooth and decays exponentially fast. Therefore,
we can reasonably truncate the whole space to a bounded domain, e.g., a square box ${\textbf B}_L := [-L,L]^d$.
The density is then rescaled to be compactly supported in a unit box ${\textbf B}_1$,
which is now  the computational domain.
One of the key ideas is to use a GauSum approximation of the kernel  $\mathcal U$ (see $U_{\textrm{GS}}$ in
\eqref{GS-Gene})  to reformulate the potential into two integrals, namely, the  \textit{long-range regular integral} 
and the \textit{short-range singular integral}. To be precise, we can reformulate the potential \eqref{NonLocalPot} as follows
\bea\label{key_form00}
\Phi(\bx) &\approx& \int_{{\textbf B}_1}{\mathcal U}(\bx-\by)\; \rho(\by) {d} \by =\int_{{\textbf B}_2}{\mathcal U}(\by)\; \rho(\bx-\by) {d} \by \\
  &=& \int_{{\textbf B}_2}U_{\textrm{GS}}(\by)\; \rho(\bx-\by) {d} \by +  \int_{ {\mathcal B}_{\delta} } \big( {\mathcal U}(\by)-U_{\textrm{GS}}(\by) \big)\; \rho(\bx-\by) {d}\by + I_{\delta} \\
\label{key_form2} &:= & I_1(\bx)  + I_2(\bx) + I_{\delta} , \quad \qquad \bx \in {\textbf B}_1,
\eea
where 
\bea\label{remInt}
 I_{\delta}=  \int_{{\textbf B}_{2}\setminus {\mathcal B}_{\delta} } \big( {\mathcal U}(\by)-U_{GS}(\by) \big)\; \rho(\bx-\by) {d}\by,
\eea
with $\mathcal B_{\delta}:= \{\bx \big | |\bx|\leq \delta\}$ being a small neighbourhood of the origin with radius $\delta \sim 10^{-4}-10^{-3}$
and
 $ U_{\textrm{GS}}$ is given explicitly as follows
\bea
\label{GS-Gene} U_{\textrm{GS}}(\by)= U_{\textrm{GS}}(|\by|):=\sum_{q= 0}^Q w_q\, e^{-\tau_q^2 |\by|^2}, \quad Q\in \mathbb N^{+},
\eea
with weights and nodes $\{(w_q, \tau_q)\}_{q=0}^{Q}$. Here, $U_{\textrm{GS}}$ designates an accurate approximation of  $\mathcal U$,  up to $\varepsilon_0 \sim 10^{-14}-10^{-16}$,
within the interval $[\delta,2]$, i.e.
\bea
\| {\mathcal U}(r)-U_{\textrm{GS}}(r)\|_{L^{\infty}([\delta,2])} \leq \varepsilon_{0}.
\eea
For \eqref{remInt}, we have $|I_\delta|\leq C \varepsilon_{0}\, \delta^{d} \, \|\rho\|_{L^{\infty}}$. 
Thus the remainder integral $I_\delta$ is negligible and is 
omitted here. Note that the GauSum approximation can be numerically computed with \textit{sinc quadrature} and we refer to \cite{EMZ2015} for more details.

To compute the \textit{regular integral} $I_1$,  plugging the explicit GauSum approximation \eqref{GS-Gene} into $I_1(\bx)$ yields
\bea
 I_1(\bx)  = \sum_{q = 0}^Q w_q \int_{{\textbf B}_2} 
 e^{-\tau_q^2 |\by|^2} \rho(\bx-\by) {d} \by, \qquad \bx \in \bau.
\eea
For $\bx \in \bau$ and $\by \in {\textbf B}_2$, we have $\bx-\by \in {\textbf  B}_3$ and we can approximate the density on ${\textbf B}_3$ by finite Fourier series.
More specifically, the density $\rho$ is well approximated by Fourier series after zero-padding to ${\textbf B}_3$ as follows
\bea\label{FourSeriB3}
\rho(\bz) \approx \sum_{\bk} \widehat{\rho}_\bk\;  \prod_{j = 1}^d  e^{\frac{\;2\pi i \;k_j}  {b_j-a_j} (z^{(j)} - a_j)}, \quad \quad \bz = (z^{(1)},\hdots,z^{(d)}) \in {\textbf B}_3,
\eea
where $a_j=-3,b_j =3, j = 1,\hdots, d$ and $\bk \in \mathbb Z^d$.
After some careful calculations, we have
\bea
I_1(\bx) & = & 
 \sum_{\bk}   \widehat{\rho}_\bk \left( \sum_{q= 0}^Q w_q    G_\bk^q \right ) \prod_{j = 1}^d  e^{\frac{\;2\pi i  \;\;k_j}  {b_j-a_j} (x^{(j)} - a_j)},
\eea
where 
\bea\label{tensor}
G_\bk^q&=& \prod_{j = 1}^d  \int_{-2}^2 e^{-\tau_q^2 |y^{(j)}|^2}\, e^{\frac{-2\pi i  k_j \; y^{(j)}}  {b_j-a_j} } {d } y^{(j)} = \prod_{j = 1}^d  \int_{0}^2 2\,e^{-\tau_q^2 |y^{(j)}|^2}\, \cos(\tfrac{2\pi  k_j \; y^{(j)}}{b_j-a_j}) {d } y^{(j)},
\eea
can be pre-computed once for all if the potential is computed on the same grid.

For the \textit{near-field correction integral} $I_2$,  within the small ball $\mathcal B_{\delta}$, 
the density function $\rho_\bx(\by) :=\rho(\bx-\by)$  is approximated by a low-order Taylor expansion as follows
\bea\label{taylor}
\rho_{\bx}(\by) \approx \mathrm P_\bx(\by)= 
\rho_{\bx}(\textbf{0}) + \sum_{j=1}^d \frac{\partial \rho_{\bx}(\textbf{0})}{\partial y_j} y_j
+ \frac{1}{2}\sum_{j,k=1}^d \frac{\partial^2 \rho_{\bx}(\textbf{0})}{\partial y_j \partial y_k} y_j\, y_k +
\frac{1}{6}\sum_{j,k,\ell=1}^d \frac{\partial^3 \rho_{\bx}(\textbf{0})}{\partial y_j \partial y_k \partial y_{\ell}} y_j\, y_k\, y_{\ell}.\eea 
Next, we integrate in spherical/polar coordinates. The computation boils down to a multiplication 
of the Laplacian $\Delta \rho$ since the contributions of the odd derivatives in \eqref{taylor} and 
off-diagonal components of the Hessian vanish. The derivatives of $\rho$ are computed by using the Fourier series approximation of the density.

The GauSum method achieves a spectral accuracy and is essentially as efficient as FFT algorithms within $O(N \log N)$ arithmetic operations.
The algorithm was implemented for the Coulomb-type kernels in \cite{EMZ2015}. The evaluations of 2D and 3D DDIs boil down 
to the Coulomb potentials with some modified densities. More explicitly, 
the 2D and 3D DDIs can be reformulated as follows
\bea
\Phi(\bx)& =& -\fl{3}{2} \left(\partial_{\bn_{\perp}\bn_{\perp}}-n_3^2 
\nabla_{\perp}^2\right)\left( \fl{1}{2\pi|\bx|} \right) \ast \rho  
= \left( \fl{1}{2\pi|\bx|} \right) \ast \left(-\fl{3}{2} (\partial_{\bn_{\perp}\bn_{\perp}}\rho 
-n_3^2 \nabla_{\perp}^2\rho )\right), \quad \bx \in \mathbb R^{2} , \quad \\
\Phi(\bx) &=& -(\bn \cdot \bn ) \rho(\bx) + \partial_{\bn \bn}\left(\frac{1}{4\pi|\bx|} \right) \ast \rho
=-(\bn \cdot \bn ) \rho(\bx) + \frac{1}{4\pi|\bx|} \ast \left(\partial_{\bn\bn}\rho \right) , \quad \bx \in \mathbb R^{3}.
\eea
Then, we need to substitute the modified densities, i.e. $\-\fl{3}{2} (\partial_{\bn_{\perp}\bn_{\perp}}\rho -n_3^2 \nabla_{\perp}^2\rho )$ and $\left(\partial_{\bn\bn}\rho \right)$ for $\rho$ in  \eqref{key_form00} for the  2D and 3D cases, respectively.

\section{Ground state computation: properties, numerical scheme and simulations}
\setcounter{equation}{0}
\setcounter{table}{0}

In this section,  we first prove some results related to the existence/non-existence of the ground states (subsection \ref{ThmOnGs}).
We next propose in subsection \ref{NumericsGrSt} an efficient and accurate numerical method for computing the ground states
by combining  the normalized gradient flow which is discretized  by the semi-implicit backward Euler Fourier pseudo-spectral method
and the Gaussian-Sum nonlocal interaction solver. We shall refer to this new method as {\sl GF-GauSum} hereafter.
Finally, subsection \ref{NumericsSectionFGPE} reports some simulations of the ground states to show some special features related to FNLSEs.

\subsection{Existence and nonexistence of the ground states}
\label{ThmOnGs}

To simplify the presentation, we divide the energy functional $\mathcal{E}(\phi)$  (\ref{energy})
 into five parts, i.e. the kinetic, potential, rotating,
local   and nonlocal interactions energy parts
\bea
\mathcal{E}(\phi(\bx))=  \mathcal{E}_{\rm kin}(\phi)+ \mathcal{E}_{\rm pot}(\phi) + 
\mathcal{E}_{\rm rot}(\phi) +  \mathcal{E}_{\rm int}(\phi) +
 \mathcal{E}_{\rm non}(\phi),
\eea
where \beas
&&\mathcal{E}_{\rm kin}(\phi) : =\fl{1}{2}\big\langle\big(-\nabla^2+m^2\big)^{s}\phi, \phi\big\rangle , 
\quad \quad \mathcal{E}_{\rm pot}(\phi) := \langle V(\bx)\phi, \phi \rangle, \\
&&\mathcal{E}_{\rm rot}(\phi) =  -\Og  \langle L_z\phi, \phi \rangle,\quad  \,\mathcal{E}_{\rm int}(\phi) := \fl{\beta}{2} \langle |\phi|^2, |\phi|^2  \rangle,
\quad \,  \mathcal{E}_{\rm non}(\phi): = \fl{\lambda}{2}\langle \Phi, |\phi|^2  \rangle,
\eeas
with $\langle f, g \rangle=\int_{\mathbb{R}^d}f\,\bar{g}\, d\bx$. 
We first prove some properties of the energy functional $\mathcal{E}(\phi(\bx))$ for any $\phi \in S$.

\begin{lemma}
\label{lemma_energy}
If  the convolution   kernel $\mathcal{U}(\bx)$ in (\ref{Kernel}) is chosen as the Coulomb-type interaction
and $V$ is the harmonic potential defined by (\ref{harm_poten}),  we have
the following properties

\begin{itemize}
\item[(i)]  For any positive $\vep >0$,  we have for $\phi \in S$ 
\be
\label{estimate_nonlocal}
\Big| \big\langle \Phi, \rho  \big\rangle  \Big|=\Big| \big\langle 
\mathcal{U}(\bx)\ast \rho , \rho \big\rangle  \Big| \le \vep \|\nabla\phi\|^2_2+\,C_\vep,
\ee
where $C_\vep$ is a real-valued constant that depends only on $d,\mu$ and $\vep$.
\item[(ii)] When $s>1$, for any $ m\ge 0$ and   $\phi\in S$, we have 
\bea
\nn
&& \int_{\mathbb{R}^d} \Big[ \fl{1}{8} \bar{\phi}\big(-\nabla^2+m^2\big)^s\phi +
\Big( V(\bx) -\fl{\gm_r^2 |\bx|^2}{2}  \Big) |\phi|^2 +\fl{\beta}{2} |\phi|^4 
\Big ]d\bx
+C_1 \le \mathcal{E}(\phi) \\[0.5em]
\label{ineqaulity_energy}
 &&\qquad\qquad \le \int_{\mathbb{R}^d} \Big[ \fl{7}{8} \bar{\phi}\big(-\nabla^2+m^2\big)^s\phi +
\Big( V(\bx) +\fl{\gm_r^2 |\bx|^2 }{2}  \Big) |\phi|^2 +\fl{\beta}{2} |\phi|^4  
\Big ]d\bx
+C_2, \qquad\qquad\quad
\eea
where $\gm_r=\min\{\gm_x, \gm_y\}$, $C_1$ and $C_2$ are two constants that only depend on $\Og$, $s$,  $\gm_r$, $d$ and $\mu$.    
\end{itemize}
\end{lemma}

\noindent
{\bf Proof.}
(i)  Using  the Hardy-Littlewood-Sobolev (HLS) inequality, we have for the Coulomb-type interaction
\be
\label{HLS}
\Big| \big\langle \Phi, \rho  \big\rangle  \Big|
=\Big| \big\langle \mathcal{U}(\bx)\ast \rho , \rho \big\rangle  \Big| 
=\fl{1}{2^{d-1}\pi} \int_{\mathbb{R}^d}\int_{\mathbb{R}^d} \fl{\rho(\bx)\rho(\by)}{|\bx-\by|^\mu}d\bx d\by 
\le c_{d,\mu} \|\rho\|^2_{p}   = c_{d,\mu} \|\phi\|^4_{2 p},
\ee
where $1<p=\fl{2d}{2d-\mu}\le \fl{2d}{d+1}<2$ and the constant
$c_{d,\mu}$ depends only on $d$ and $\mu$. 
For the 3D case,   let us introduce  $\sigma=\fl{3-p}{2p}$.  By the H\"{o}lder's inequality, Young's inequality
and the embedding theorem, we obtain    
\be
\label{3d_Cou_est}
 \Big| \big\langle \Phi, \rho  \big\rangle  \Big|
 \le
 c_{d,\mu} \Big( \|\phi\|_2^\sigma\; \|\phi\|_6^{1-\sigma} \Big)^{4} 
 =c_{d,\mu} \Big( \|\phi\|_6^2 \Big)^{2(1-\sigma)} 
 \le~\widetilde{\vep}\, \|\phi\|^2_6 +C_\vep
 \le \vep \|\nabla\phi\|^2_2 + C_\vep.
\ee
Similarly, for the 2D case,  let $q=\fl{4p}{2-p}>2p>2$ and $\sigma=\fl{q-2p}{p (q-2)}$. Then, one gets
\be
\label{2d_Cou_est}
 \Big| \big\langle \Phi, \rho  \big\rangle  \Big|
 \le
 c_{d,\mu} \Big( \|\phi\|_2^\sigma\; \|\phi\|_q^{1-\sigma} \Big)^{4} 
 =c_{d,\mu} \Big( \|\phi\|_q^2 \Big)^{2(1-\sigma)}
 \le \widetilde{\vep}\,\|\phi\|^2_q + \widetilde{C_\vep} 
 \le \vep \|\nabla\phi\|^2_2 +C_\vep .
\ee

\medskip

\noindent
(ii) Let $\gm_r=\min\{\gm_x, \gm_y\}$.
By Young's inequality and Plancherel's formula, we have
\bea
\nn
&&\Big|\,\Og\int_{\mathbb{R}^d}\bar{\phi}L_z\phi \;d\bx\,\Big| 
\le  \int_{\mathbb{R}^d} \Big[ \big| (\gm_r x\bar{\phi})\, (\Og\p_y\phi/\gm_r) \,\big| 
+ \big| (\gm_r y\bar{\phi})\, (\Og\p_x\phi/\gm_r) \,\big|  d\bx\Big]  \qquad\qquad\\[0.5em]
\nn  
&&\le \fl{\gm_r^2}{2} \int_{\mathbb{R}^d}|\bx|^2|\phi|^2d\bx 
+\fl{\Og^2}{2\gm_r^2}  \int_{\mathbb{R}^d} \big|\nabla\phi\big|^2\,d\bx
\;=  \fl{\gm_r^2}{2} \int_{\mathbb{R}^d}|\bx|^2 |\phi|^2d\bx + 
\fl{\Og^2}{2\gm_r^2 (2\pi)^d}\int_{\mathbb{R}^d} |\bk|^2|\widehat{\phi}|^2 d\bk
\qquad\qquad \\[0.5em] \nn
&&\le  \fl{\gm_r^2}{2} \int_{\mathbb{R}^d}|\bx|^2 |\phi|^2d\bx 
+ \fl{\Og^2}{2\gm_r^2 (2\pi)^d}\int_{\mathbb{R}^d} \Big[  \fl{\gm_r^2 (|\bk|^2+m^2)^{s} }{2 \Og^2} + 
\Big(\fl{ \gm_r^2}{2\Og^2}\Big)^{\fl{1}{1-s}}\,\Big] |\widehat{\phi}|^2 \,d\bk 
-\fl{\Og^2 m^2}{2\gm^2_r}
\\[0.5em]
&& \label{proof_ineq}\le \fl{\gm_r^2}{2} \int_{\mathbb{R}^d}|\bx|^2|\phi|^2d\bx 
+ \fl{1}{4} \int_{\mathbb{R}^d} \bar{\phi}\big( -\nabla^2+m^2\big)^{s} \phi\, d\bx+C.
\eea 
Similarly, for the Coulomb-type nonlocal interaction, we obtain 
\bea
\nn
\Big| \fl{\lambda}{2} \big\langle \Phi, \rho  \big\rangle  \Big|
& \le &\widetilde{\vep} \|\nabla\phi\|^2_2+\widetilde{C_\vep}
 =\fl{\widetilde{\vep}}{(2\pi)^d}\int_{\mathbb{R}^d} |\bk|^2 \,|\widehat{\phi}|^2 \ d\bk+\widetilde{C_\vep}
 \le \fl{1}{(2\pi)^d}\int_{\mathbb{R}^d}
 \Big[ \fl{1}{8}(|\bk|^2+m^2)^s + 
C \Big]\,|\widehat{\phi}|^2  \, d\bk \\[0.5em]
 \label{Cou_inq}
 &=&\fl{1}{8}  \int_{\mathbb{R}^d} \bar{\phi}\big( -\nabla^2+m^2\big)^{s} \phi\, d\bx+ C.
\eea
Therefore, the inequality (\ref{ineqaulity_energy}) follows  from (\ref{proof_ineq}) and (\ref{Cou_inq}).
\hfill $\square$\\
\begin{thm}
\label{existence_gs}
If $V(\bx)$ is a trapping harmonic potential defined in (\ref{harm_poten}), then the following properties hold.
\begin{itemize}
\item[(i)] If $s>1$ and  $\beta\ge0$, then  there exists a ground state of the FNLSE for all $\Omega>0$ if one of the following conditions holds:
\begin{itemize}
\item[(A)]  $\mathcal{U}(\bx)$ reads as either Coulomb-type.

\item[(B)]  For 3D DDI:   $-\beta/2\le\lambda\le\beta.$

\item[(C)]  For 2D DDI:  (c1) $\lambda=0$. (c2)  $\lambda>0$ and $n_3=0$. (c3)  $\lambda<0$ and $n_3^2 \ge\fl{1}{2}$.
\end{itemize}
\end{itemize}
\begin{itemize}
\item[(ii)] If $\Og=0$, $\beta>0$ and $\lambda>0$, then the ground state of the FNLSE exists for all $s>0$. 

\item[(iii)] If $0< s <1$,   there exists no ground state if  one of the following conditions holds

\begin{itemize}
\item[(A)] $\forall \;\Og>0$, $\mathcal{U}(\bx)$ is a Coulomb-type interaction or a 3D DDI.

\item[(B)] $\mathcal{U}(\bx)$ is the 2D DDI, $\forall \;\Og>\Og_0=c|\lambda|^{\fl{2}{5}}$ with 
$c=\big(\fl{(2\pi^2+1)^4\gm^6}{48 e \pi^9}\big)^{\fl{1}{5}}(\approx 0.54 \textrm{ for } \gamma=1).$
Here, $\gm=\max\{\gm_x, \gm_y\}$
\end{itemize}
\end{itemize}
\end{thm}

\noindent
{\bf Proof:}
(i) For the Coulomb-type interaction, it is clear by Lemma \ref{lemma_energy}  that the energy functional
$\mathcal{E}$ is bounded below, coercive and weakly lower semi-continuous on $S$. Hence, (A) follows.
For the  DDI, the proof is similar as those for the non-fractional case \cite{BAC2012, BCW2010}
by noticing (\ref{proof_ineq}).  Similar arguments lead to (ii).

%
%

\smallskip

\noindent
(iii) Denote $\gm=\max\{\gm_x, \gm_y\}$.  In 2D, we choose the function
\be
\label{nonexist_fun}
\phi_n(\bx)=\mathcal{F}^{-1}(\widehat{\phi_n})(\bx), \quad  {\rm with} \quad
\widehat{\phi_n}(\bk)=\mathcal{F}(\phi_n)(\bk)=
\left(4\pi \vep^{n+1}\right)^{1/2} (n!)^{-1/2} {\rm exp}(-\vep|\bk|^2/2) |\bk|^ne^{i n\theta}.
\ee
By  Plancherel's formula,  it is easy to check that  $\|\phi_n\|_2=\fl{1}{2\pi}\|\widehat{\phi_n}\|_2=1$, and thus $\phi_n \in S$.
Let $\rho_n=|\phi_n |^2$. By Young's inequality and Cauchy-Schwarz inequality, we obtain
\bea
\nn
\mathcal{E}_1(\phi_n)&=:&\mathcal{E}_{\rm kin}(\phi_n)+
\mathcal{E}_{\rm pot}(\phi_n) + 
\mathcal{E}_{\rm rot}(\phi_n)\\[0.5em]
\nn
&\le& \fl{1}{4\pi^2}\bigg[ \fl{1}{2}\langle(|\bk|^2+m^2\big)^{s}\widehat{\phi}_n,  \widehat{\phi}_n\rangle
-\fl{\gm^2}{2}\langle\Delta\widehat{\phi}_n,\widehat{\phi}_n\rangle
-i\Omega\langle \widehat{J}_{z_k}\widehat{\phi}_n,\widehat{\phi}_n \rangle \bigg] \\[0.5em]
&\le&\fl{e^{m^2 \vep}}{2\vep^{s}}\fl{\Gamma(n+1+s)}{ \Gamma(n+1)}
+\bigg(  \fl{\vep \gm^2}{2}-\Og \bigg) n+ \fl{\vep \gm^2}{2},
\\[1em]
\nn
|\mathcal{E}_{\rm int}(\phi_n)|
&=& \fl{|\beta|}{2} \| \rho_n\|_2^2
=\fl{|\beta|}{2(2\pi)^6} \| \widehat{\phi_n}\ast \widehat{\phi_n^*}\|^2_2 
\le \fl{|\beta|}{2(2\pi)^4} \left(\fl{1}{(2\pi)^2}\| \widehat{\phi_n} \|^2_{2}\right)  \| 
\widehat{\bar{\phi}_n} \|^2_{1}  \\[0.5em]
\label{2d-inter}
&=&\fl{|\beta|}{2(2\pi)^4}  \| \widehat{\phi_n} \|^2_{1}
=\fl{|\beta| }{2\pi \vep}\; 
\fl{ 2^{n} \big(\Gamma(n/2+1)\big)^2}{\Gamma(n+1)}.
\eea
Furthermore, we compute the nonlocal interaction energy $\mathcal{E}_{\rm non}(\phi_n)$. 
For the Coulomb-type interaction,  by using the HLS inequality (\ref{HLS})
and the H\"{o}lder's inequality, we obtain
\bea
\label{nonlocal_2dCou_est}
\big|\mathcal{E}_{\rm non}(\phi_n)\big|
&\le& \fl{c_{\mu} |\lambda|}{2} \|\rho_n\|^2_{p} 
\le \widetilde{c}_{\mu}  \|\rho_n\|_1^{2-\mu}  \|\rho_n\|_{2}^{\mu} 
\le   \widetilde{c}_{\mu} \left[ \fl{1}{\pi\vep }\; 
\fl{ 2^{n} \big(\Gamma(n/2+1)\big)^2}{\Gamma(n+1)}
 \right]^{\mu/2},
\eea
where $p=\fl{4}{4-\mu}$,  $0<\mu \le1$, and $ \widetilde{c}_{\mu}$  depends only on $\mu$ and $\lambda$.  
Together  with the Stirling's formula
\be
\label{Stirling}
\Gamma(x+1) \thicksim \sqrt{2\pi x}\left(\fl{x}{e}\right)^x, \textrm{ when} \quad x \rightarrow \infty,
\ee
one gets 
\be
|\mathcal{E}_{\rm int}(\phi_n)|\thicksim  \sqrt{n}, \qquad
|\mathcal{E}_{\rm non}(\phi_n)|\thicksim  \sqrt{n},
\qquad  n  \rightarrow \infty.
\ee 
Let $\vep<\fl{2\Og}{\gm^2}$, $\forall\, \Og>0$ and $s<1$, we then prove that
\bea 
\limsup_{n \rightarrow \infty}
 \mathcal{E}(\phi_{n}) 
 \nn
&\le& \limsup_{n \rightarrow \infty}
\big[ \mathcal{E}_1(\phi_{n})  + |\mathcal{E}_{\rm int}(\phi_{n})|
+\mathcal{E}_{\rm non}(\phi_{n}) \big]   \\
&\le&
\limsup_{n \rightarrow \infty}
\left[ 
\fl{c_2 n^s}{2\vep^s}
 + \bigg(  \fl{\vep \gm^2}{2}-\Og \bigg) n +c_1 \sqrt{n}
+c_0  \right]  
= -\infty,
\eea
which implies the nonexistence of the ground states.

For the 2D DDI,  we have
\bea
\label{nonlocal_2dDDI_est}
\big|\mathcal{E}_{\rm non}(\phi_n)\big|
&=&\fl{|\lambda|}{8\pi^2} \big|\big \langle \widehat{\mathcal{U}}_{\rm ker} \widehat{\rho_n}, \widehat{\rho_n} \big\rangle \big|
\le \fl{|\lambda|}{8\pi^2} \Big\langle \fl{3 (|\bk\cdot\bn_{\perp}|^2+n_3^2 |\bk|^2)}{2 |\bk|} \widehat{\rho_n}, \widehat{\rho_n}     \Big\rangle
\le \fl{ 3 |\lambda|}{16 \pi^2} \big\langle |\bk|\, \widehat{\rho_n}, \widehat{\rho_n}     \big\rangle.
\eea
By the generalized  Minkowski inequality, we prove that
\bea\nn
4\pi^2\Big( \big\langle |\bk|\, \widehat{\rho_n}, \widehat{\rho_n}  \big\rangle \Big)^\fl{1}{2}
&=&\left[ \int_{\mathbb{R}^2} |\bk| \Big| \int_{\mathbb{R}^2} \widehat{\phi_n}(\bk-\bm{\xi}) 
\widehat{\bar{\phi}_n}(\bm{\xi}) d\bm{\xi} \Big|^2 d\bk \right]^\fl{1}{2} 
\le \int_{\mathbb{R}^2}   \left[ \int_{\mathbb{R}^2} |\bk| \big| \widehat{\phi_n}(\bk-\bm{\xi}) \big|^2 
 \big| \widehat{\bar{\phi}_n}(\bm{\xi}) \big|^2 d\bk \right]^\fl{1}{2} d\bm{\xi}  
\\[0.5em] \nn
&=& \int_{\mathbb{R}^2}  \big| \widehat{\bar{\phi}_n}(\bm{\xi})\big|  \left[ \int_{\mathbb{R}^2} |\bk-\bm{\xi}| \big| \widehat{\phi_n}(\bk) \big|^2 
       d\bk \right]^\fl{1}{2} d\bm{\xi} 
 \le \int_{\mathbb{R}^2}  \big| \widehat{\bar{\phi}_n}(\bm{\xi})\big|  \left[ \int_{\mathbb{R}^2} (|\bk|+|\bm{\xi}|) \big| \widehat{\phi_n}(\bk) \big|^2 
       d\bk \right]^\fl{1}{2} d\bm{\xi}  
\\[0.5em] \nn       
&\le&  \int_{\mathbb{R}^2}  \big| \widehat{\bar{\phi}_n}(\bm{\xi}) \big| 
\left[ \sqrt{|\bm{\xi}|} + \Big(\int_{\mathbb{R}^2}  |\bk| \big| \widehat{\phi_n}(\bk) \big|^2 
       d\bk \Big)^\fl{1}{2} \right] d\bm{\xi} 
 = \left(\big\|\sqrt{|\bk|}\,\widehat{\phi_n}\big\|_1 + \big\|\widehat{\phi_n}\big\|_1\,\|\sqrt{|\bk|}\,\widehat{\phi_n}\|_2\right)  
 \\[0.5em] 
 &=&\fl{\pi^{\fl{1}{2}}2^{\fl{n}{2}+\fl{9}{4}}}{\vep^{\fl{3}{4}}}\, \fl{\Gamma(\fl{n}{2}+\fl{5}{4})}{\sqrt{ \Gamma(n+1)} }
 +\fl{\pi^{\fl{5}{2}} 2^{\fl{n}{2}+3}}{\vep^{\fl{3}{4}}} \,\fl{ \sqrt{\Gamma(n+\fl{3}{2})} \, \Gamma(\fl{n}{2}+1)}{\Gamma(n+1)}. 
\eea
Again, by the Stirling's formula (\ref{Stirling}), we show that
\be
|\mathcal{E}_{\rm non}(\phi_n)| \precsim  \fl{c_3 |\lambda| }{\vep^{3/2}}\,n,\qquad 
\qquad  n  \rightarrow \infty,
\ee
where $c_3=\fl{3\sqrt{2} (2\pi^2+1)^2}{32\pi \sqrt[4]{e \pi}}.$
Let us set $\vep=\Big(3c_3|\lambda|/\gm^2 \Big)^{\fl{2}{5}} $    and 
$\Og>\Og_0=\big(\fl{(2\pi^2+1)^4\gm^6}{48 e \pi^9}\big)^{\fl{1}{5}} |\lambda|^{\fl{2}{5}}$. It follows that 
\bea 
\limsup_{n \rightarrow \infty}
 \mathcal{E}(\phi_{n}) 
 \nn
&\le&
\limsup_{n \rightarrow \infty}
\left[ 
\fl{c_2 n^s}{2\vep^s}
 + \bigg(  \fl{\vep \gm^2}{2}+\fl{c_3}{\vep^{\fl{3}{2}}}-\Og \bigg) n +c_1 \sqrt{n}
+c_0  \right]  
= -\infty,
\eea
leading to  the nonexistence of the ground states.

In 3D, we choose the sequence
\be
\phi^{\textrm{3D}}_n(\bx)=\phi_n(x, y) \phi(z),
\ee
where $\phi(z) =\left(\fl{\gm_z}{\pi}\right)^{1/4} \exp\{-\fl{\gm_z z^2}{2} \}$ and
 $\phi_n(x, y)=\mathcal{F}^{-1}(\widehat{\phi}_n(\bk))$, with $\widehat{\phi}_n(\bk)$ reading as (\ref{nonexist_fun}).
Then, the argument proceeds similarly as those in 2D for the 3D Coulomb potential. As for the 3D DDI, noticing that
\be
|\mathcal{E}_{non}(\phi^{\textrm{3D}}_n)| \le \fl{3 |\lambda|}{2}\|\phi^{\textrm{3D}}_n\|_4^4=\fl{3 |\lambda |\sqrt{\gm_z} }{2\sqrt{2\pi}}\, \|\phi_n\|_4^4
=\fl{3 |\lambda |\sqrt{\gm_z} }{2\sqrt{2\pi}}\, \|\rho_n\|_2^2,
\ee
the left argument proceeds similarly as those in 2D from  (\ref{2d-inter}).

\hfill $\square$

\begin{remark}
For the 2D DDI, one open question concerns the plausible fact that (iii)(B) in Theorem (\ref{existence_gs}) maybe hold for $\forall$ $\Og_0>0$.
The proof presented here does not seem to be directly applicable for this conjecture.
\end{remark}

\begin{remark}
It might be interesting to understand the existence/non-existence and the uniqueness of the ground states  for the more general 
FNLSE
\be
\label{GenFSE}
i\p_t\psi=\left[\fl{1}{2}(-\nabla^2+m^2)^s+\frac{1}{2}\gm_r^2 |\bx|^p+\beta |\psi|^q +\lambda\Phi-\Og L_z\right]\psi,
\ee
where the constants $\beta$ and $\lambda$ can be positive or negative and the powers $p$ and $q$ are real-valued
positive constants.  
We leave it as an open problem for some future studies.
\end{remark}

\subsection{Numerical method}\label{NumericsGrSt}


For a constant time step $\Delta t$, we introduce the discrete times $t_n=n\Delta t$ for $n=0, 1, 2,\ldots$
The gradient flow with discrete normalization (GFDN) method reads as
\bea \label{gf-eq1}
&&\partial_t \phi(\bx,t)= -\left[\frac{1}{2} (-\nabla^2+m^2)^{s}  +
V(\bx)+\beta |\phi|^2 + \lambda \,\Phi(\bx,t) -\Og L_z\right]
\phi(\bx,t),\\
\label{gf-colb1}
&& \Phi(\bx,t) = \left(\mathcal{U}\ast |\phi|^2\right)(\bx,t), \;\;\;\quad\quad \qquad \qquad \qquad \qquad\quad
\bx \in {\mathbb{ R}}^d, \quad t_n\le t< t_{n+1},\\
\label{gffg3}
&&\phi(\bx,t_{n+1})=\frac{\phi(\bx,t_{n+1}^-)}{\|\phi(\bx,t_{n+1}^-)\|_{2}}, \;\qquad \quad \qquad \qquad \qquad \qquad
\bx \in {\mathbb{ R}}^d, \quad n\ge0,
\eea
with the initial data
\be
\phi(\bx,0)=\phi_0(\bx), \qquad \bx \in {\mathbb{ R}}^d, \qquad {\rm with} \qquad \|\phi_0\|_{2}=1.
\ee
Let $\phi^n(\bx)$ and $\Phi^n(\bx)$ be the approximations of $\phi(\bx,t_n)$ and
$\Phi(\bx,t_n)$, respectively.
The above GFDN is usually discretized in time \textit{via} the semi-implicit backward Euler method \cite{AD1,BJTZ2015,
BTZ2015,ZD2011}
\bea 
\label{gf-eq2}
&&\frac{\phi^{(1)}(\bx)-\phi^{n}(\bx)}{\Delta t}=- \left[\frac{1}{2} (-\nabla^2+m^2)^{s}
+V(\bx)+ \beta |\phi^n|^2+\lambda \,\Phi^n(\bx) -\Og L_z\right]
\phi^{(1)}(\bx),\\
\label{gf-colb2}
&& \Phi^n(\bx) = \left(\mathcal{U}\ast |\phi^n|^2\right)(\bx),
\quad\qquad\qquad\qquad \qquad \quad \qquad  
\bx \in {\mathbb{ R}}^d, \\
\label{gffg33}
&&\phi^{n+1}(\bx)=\frac{\phi^{(1)}(\bx)}{\|\phi^{(1)}(\bx)\|_{2}}, \quad \qquad\qquad \qquad \qquad \qquad \qquad
\bx \in {\mathbb{ R}}^d, \quad n\ge0.
\eea
The ground states  decay exponentially fast due to the trapping potential.
Therefore, in practical computations,
we first truncate the whole space to a bounded rectangular domain and impose periodic boundary conditions.
Then, we discretize the equation (\ref{gf-eq2}) \textit{via} the Fourier pseudo-spectral method in space and 
evaluate the nonlocal interaction $\Phi^n(\bx)$ by the GauSum solver.  
The full discretized  scheme of  system (\ref{gf-eq2})-(\ref{gffg33})
can then be solved by a fixed-point iteration or a preconditioned Krylov subspace solver
 with a similar preconditioner as those in \cite{AD1}. Let us define the operators 
\bea
&& A^{\rm BE,n}:=\fl{ I}{\Delta t}+\frac{1}{2} (-\nabla^2+m^2)^{s}
	+V(\bx)+ \beta |\phi^n|^2+\lambda \,\Phi^n(\bx) -\Og L_z, \\
&&{ P}^{\rm BE}_{\Delta }=\left[ \fl{ I}{\Delta t}+ \frac{1}{2} (-\nabla^2+m^2)^{s} \right]^{-1},
	\qquad { A}^{\rm BE,n}_{{\rm TF}}=V(\bx)+ \beta |\phi^n|^2+\lambda \,\Phi^n(\bx) -\Og L_z, \\
&&	{ P}^{\rm BE,n}_{\rm TF}=\left[ \fl{  I}{\Delta t}+V(\bx)+ \beta |\phi^n|^2+\lambda \,\Phi^n(\bx)  \right]^{-1},
	\qquad  { A}^{\rm BE,n}_{\Delta, \Og}= \frac{1}{2} (-\nabla^2+m^2)^{s}-\Og L_z.
\eea
Moreover,  we denote by $\mathbb{I}$, $\mathbb{A}^{\rm BE,n}$, $\mathbb{P}^{\rm BE}_\Delta$, $\mathbb{A}^{\rm BE,n}_{{\rm TF}}$, 
$\mathbb{P}^{\rm BE,n}_{\rm TF}$,  $\mathbb{A}^{\rm BE,n}_{\Delta,\Og}$ 
 the discretized versions of the above operators, and by $\phi^{(1)}$ and $\phi^n$
 the discretization of
 $\bm{\phi}^{(1)}$ and $\bm{\phi}^n$
 through the Fourier pseudo-spectral approximation.
   Then, the finite-dimensional linear system corresponding to  \eqref{gf-eq2}-\eqref{gffg33} reads as
\be
\mathbb{A}^{\rm BE,n}\bm{\phi}^{(1)}=\bm{b}^n:=\bm{\phi}^n/\Delta t.
\ee
Two preconditioned versions of the linear system are the following 
\be
\label{precond_type1}
\Big(\mathbb{I}+\mathbb{P}^{\rm BE}_{\Delta}
\mathbb{A}^{\rm BE,n}_{ {\rm TF}}\Big)\bm{\phi}^{(1)}
	=\mathbb{P}^{\rm BE}_{\Delta}\bm{b}^{n},
\qquad {\rm or} \qquad
\Big(\mathbb{I}+\mathbb{P}^{\rm BE,n}_{\rm TF}
\mathbb{A}^{\rm BE,n}_{\Delta,\Og}\Big)\bm{\phi}^{(1)}
	=\mathbb{P}^{\rm BE,n}_{\rm TF}\bm{b}^{n}.
\ee
We refer the reader to \cite{AD1} for more details and omit them here for brevity. 
Like in the standard case \cite{AD1,ADBookChapter}, the most efficient solver uses the first preconditioned system (left) in (\ref{precond_type1})
based on $\mathbb{P}^{\rm BE}_{\Delta}$. In particular, the acceleration of the convergence of the Krylov 
subspace solver (BiCGStab) is visible when $\Omega$, 
$\beta$ and $\lambda$ are large. In practice, we use this preconditioned solver in subsection \ref{NumericsSectionFGPE}.

\subsection{Numerical results}\label{NumericsSectionFGPE}

In this subsection, we report some numerical results concerning  the ground states of (\ref{SFSE})-(\ref{NonLocalPot}) computed by  the
{\sl GF-GauSum} solver built in the previous subsection.
To this end, unless stated, we fix $m=0$ and  $d=2$. We carry
 out the computation  on the domain ${\textbf B}=[-32, 32]\times[-32, 32]$
that is discretized with uniform mesh sizes $h_x=h_y=\fl{1}{8}$. We use a constant time step $\Delta t=10^{-3}$.
 The trapping potential  $V(\bx)$ is chosen as (\ref{harm_poten}) with $\gm_x=\gm_y=1$. 
The nonlocal interaction is of Coulomb-type with $\mu=1$. The initial guess $\phi_0(\bx)$ is chosen as
\be\label{ini_guess}
\phi_0(\bx)=\fl{(1-\Og)\phi_{\rm ho}(\bx) + \Og \phi^v_{\rm ho}(\bx)}{\|(1-\Og)\phi_{\rm ho}(\bx) + \Og \phi^v_{\rm ho}(\bx)\|},\;\;
{\rm with}\; \;
\phi_{\rm ho}(\bx)=\fl{1}{\sqrt{\pi}}\,e^{-\fl{|\bx|^2}{2}},\; 
\phi^v_{\rm ho}(\bx)=\fl{x+i y}{\sqrt{\pi}}\,e^{-\fl{|\bx|^2}{2}}, \; \bx\in {\textbf B}.
\ee
The ground state $\phi_g(\bx)$ is reached when the stopping criterion holds:
 $\| \phi^n(\bx)-\phi^{n+1}(\bx)\|_{\infty} \le \vep_0\, \Delta t $. 
In the computations,  we choose the   accuracy parameter $\vep_0= 10^{-9}$.

\begin{figure}[h!]
\centerline{
\psfig{figure=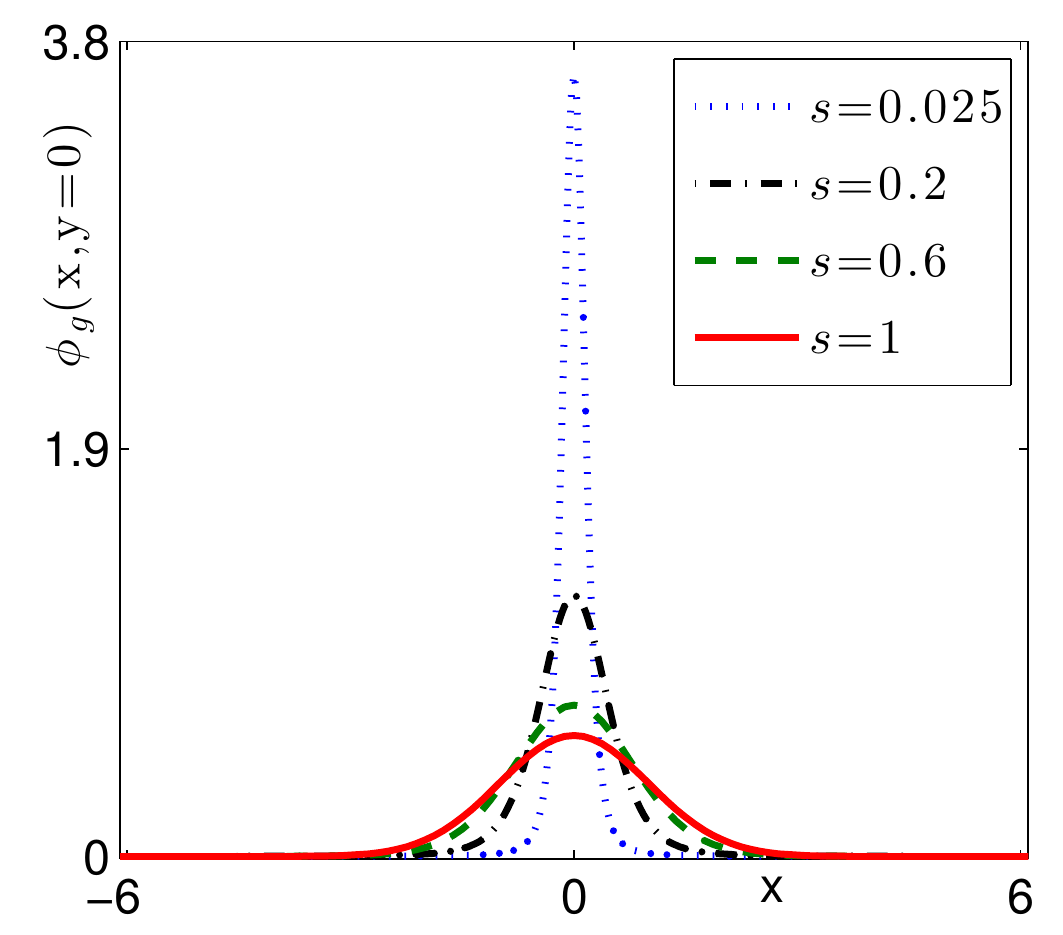,height=3.9cm,width=4.0cm,angle=0}
\psfig{figure=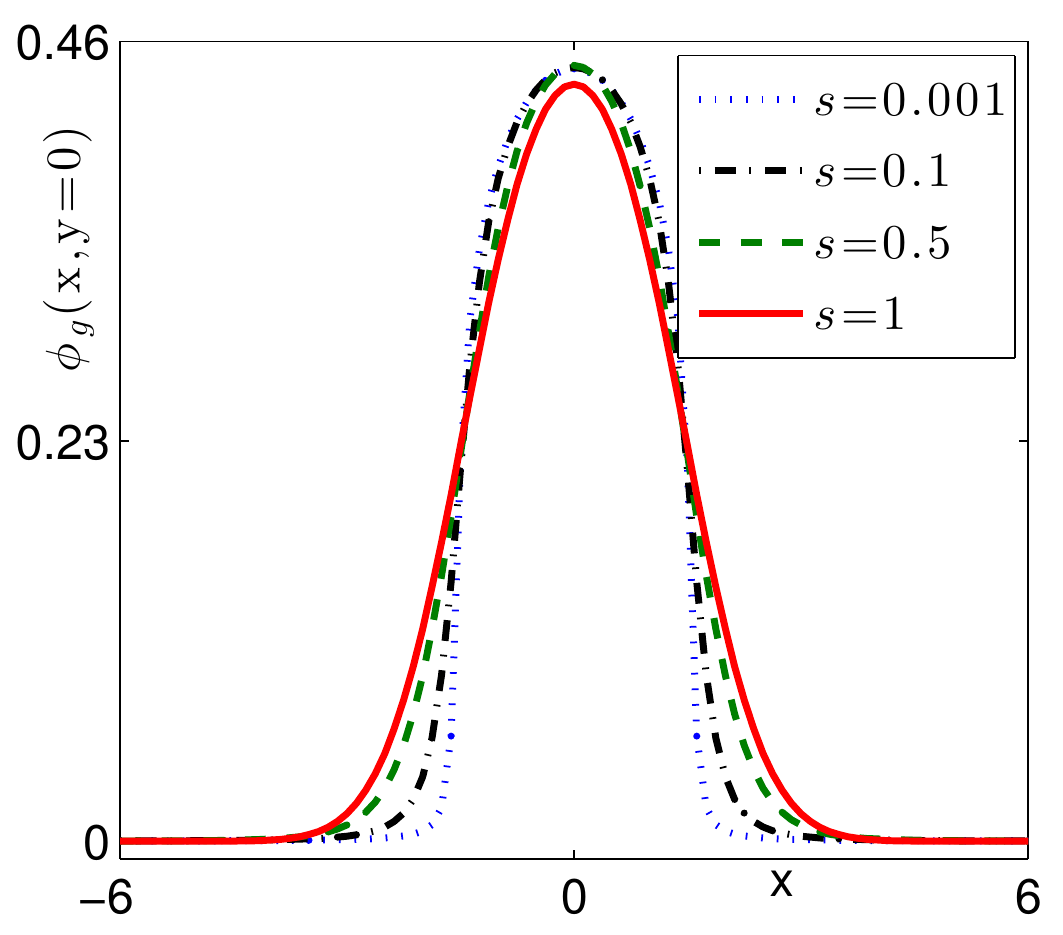,height=3.9cm,width=4.0cm,angle=0}
\psfig{figure=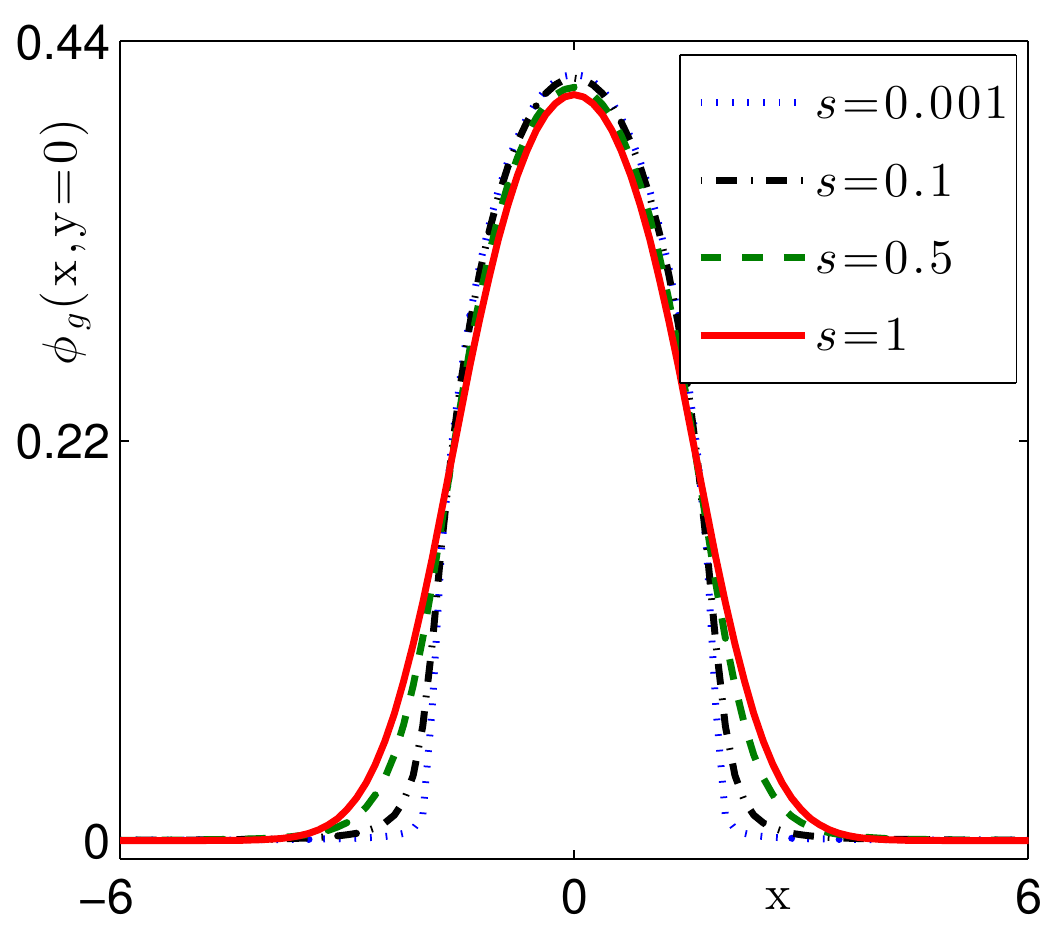,height=3.9cm,width=4.0cm,angle=0}
\psfig{figure=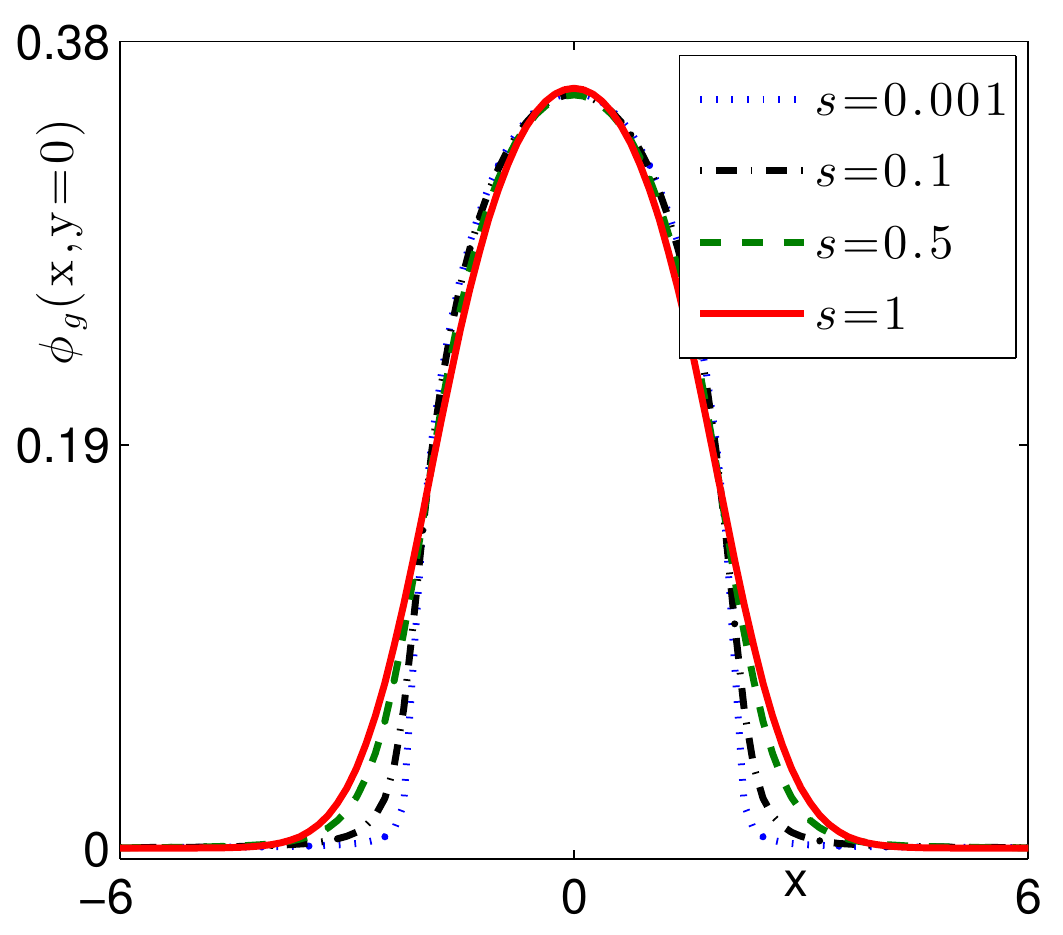,height=3.9cm,width=4.0cm,angle=0}
}
\centerline{
\psfig{figure=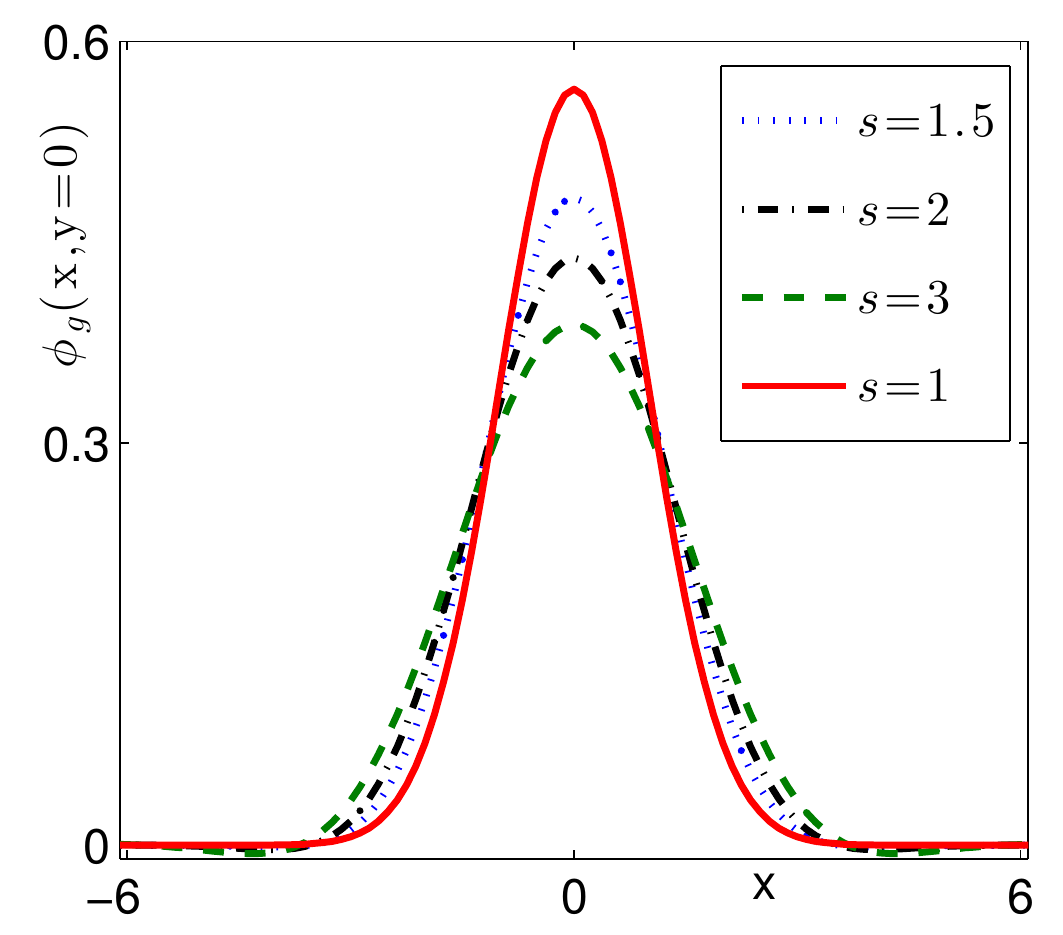,height=3.9cm,width=4.cm,angle=0}
\psfig{figure=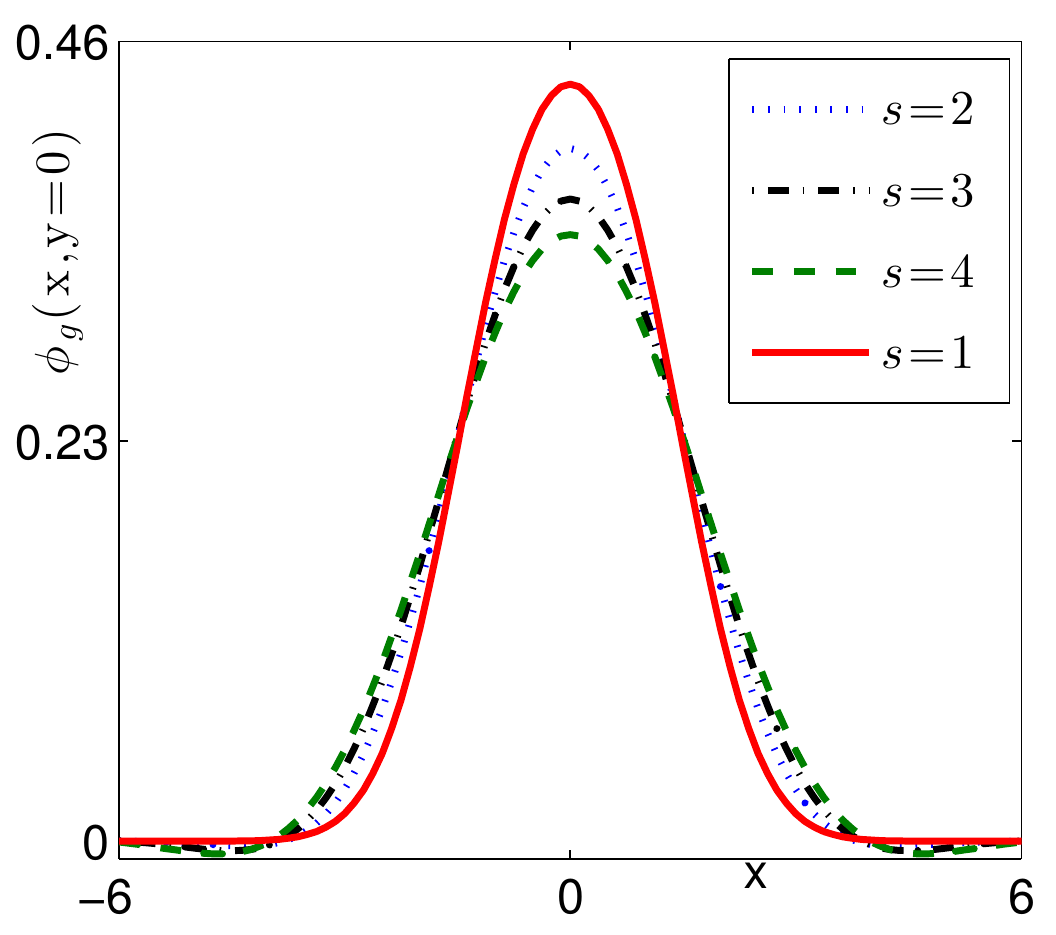,height=3.9cm,width=4.cm,angle=0}
\psfig{figure=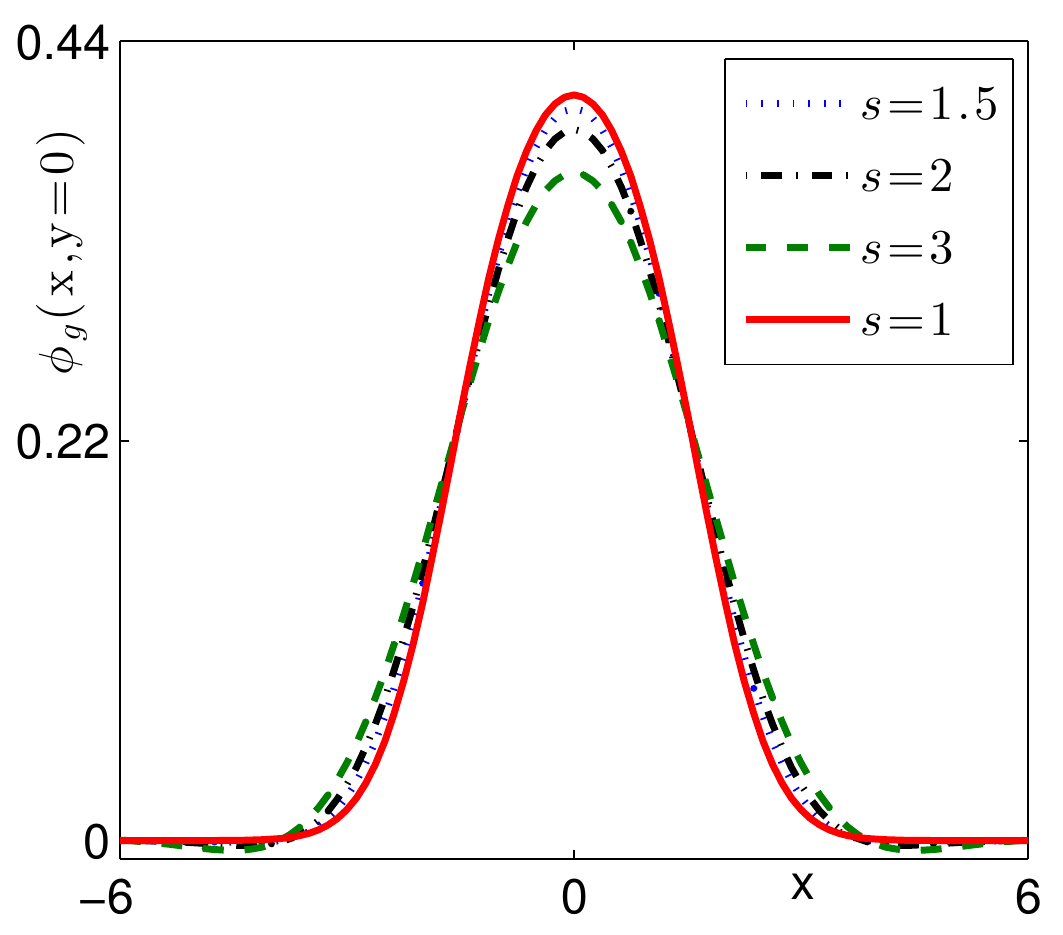,height=3.9cm,width=4.cm,angle=0}
\psfig{figure=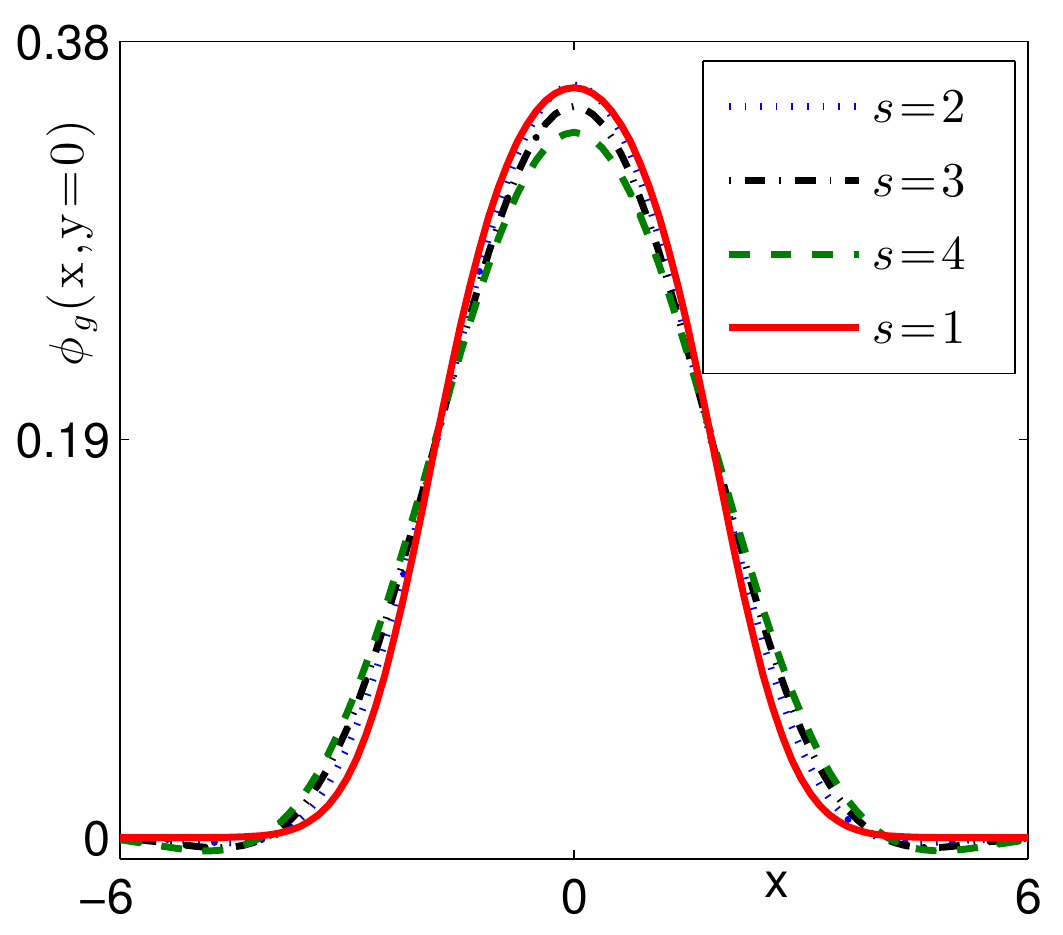,height=3.9cm,width=4.cm,angle=0}
}
\caption{ Slice plots of $\phi_g(x,0)$ for {\bf Cases} {\bf I}--{\bf IV} (from left to right) for subdispersion
$s\le 1$ (top row)  and superdispersion $s\ge 1$ (bottom)  in example \ref{GS: NonRot}.}
\label{fig:ex1}
\end{figure}

\begin{exmp}\label{GS: NonRot} {\bf Non-rotating FNLSE.}
Here, we impose $\Og=0$. We study the ground states of the following four cases:
\begin{itemize}
\item {\bf Case I}.    Linear case, i.e. $\beta=\lambda=0$.
\item {\bf Case II}.   Purely long-range interaction, i.e. $\beta=0$ and $\lambda=10$.
\item {\bf Case III}.   Purely short-range interaction, i.e. $\beta=10$ and
$\lambda=0$.
\item  {\bf Case IV}. Both long-range and short-range interactions, i.e. $\lambda=\beta=10$.
\end{itemize}
\end{exmp}
Figure \ref{fig:ex1} shows the slice plots of the ground states along the $x$-axis, i.e. $\phi_g(x,0)$, for different fractional orders $s$
of the FNLSE.

\begin{figure}[h!]
\centerline{
\quad\psfig{figure=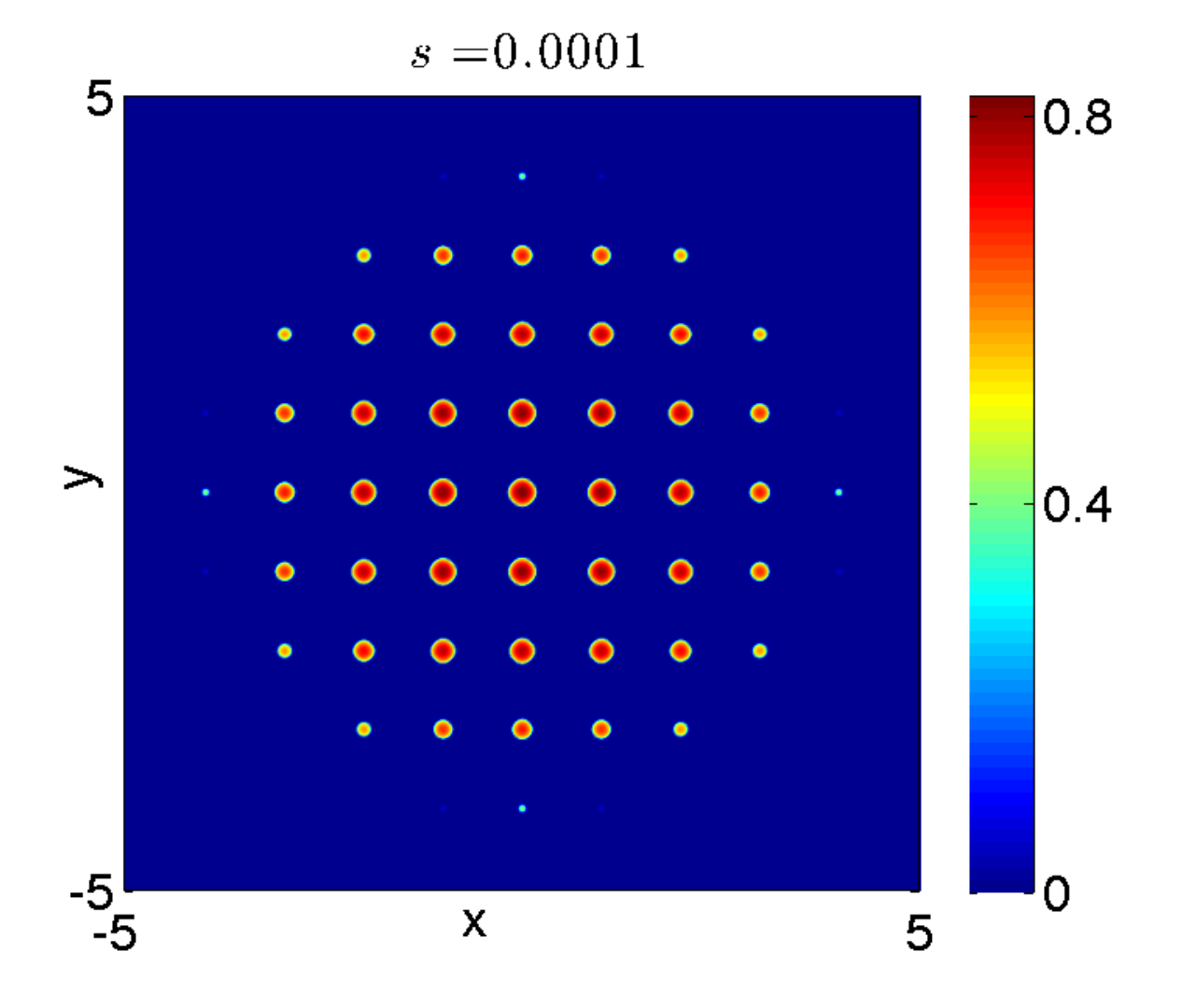,height=3.6cm,width=4.2cm,angle=0}
\psfig{figure=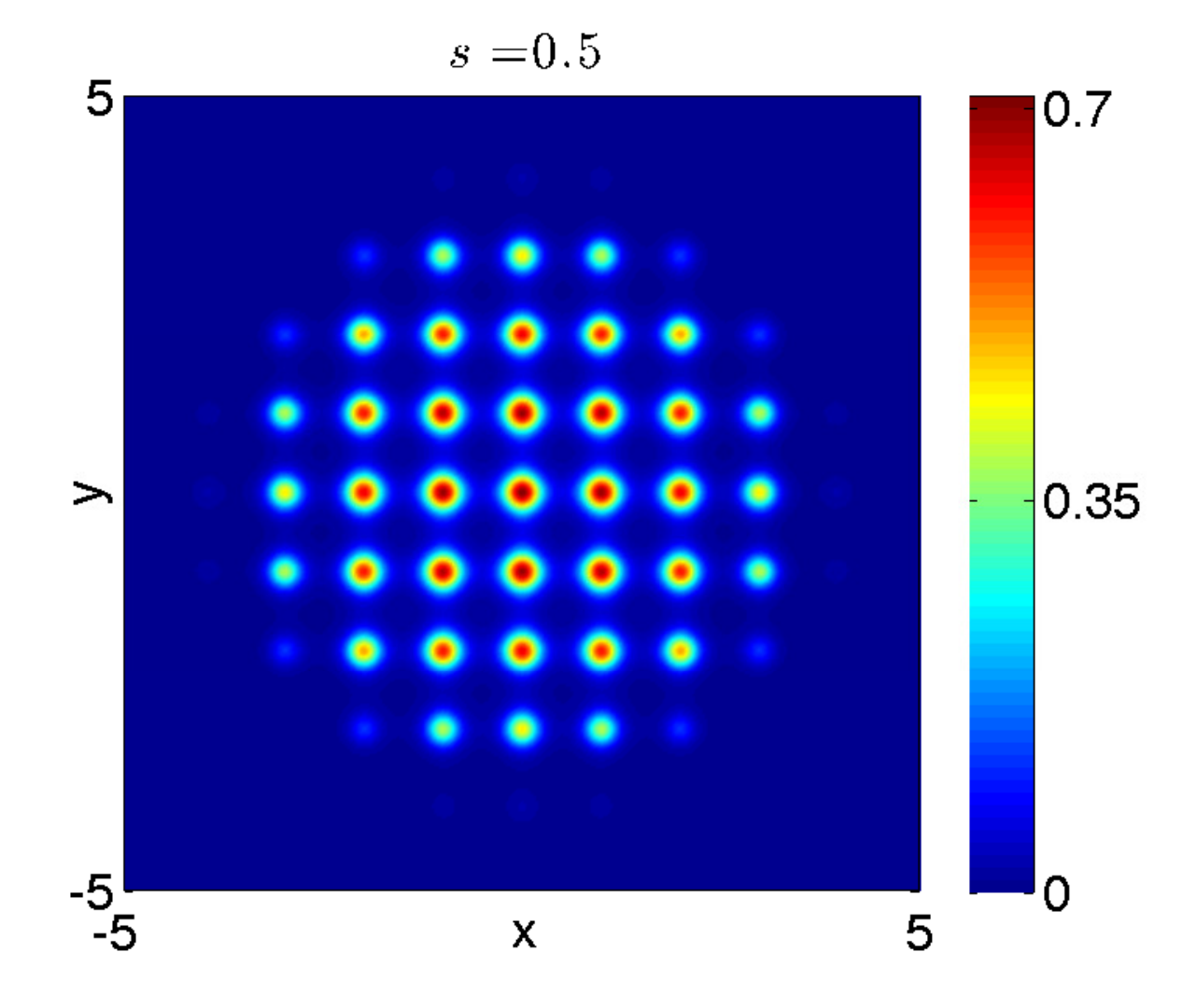,height=3.6cm,width=4.5cm,angle=0}
\psfig{figure=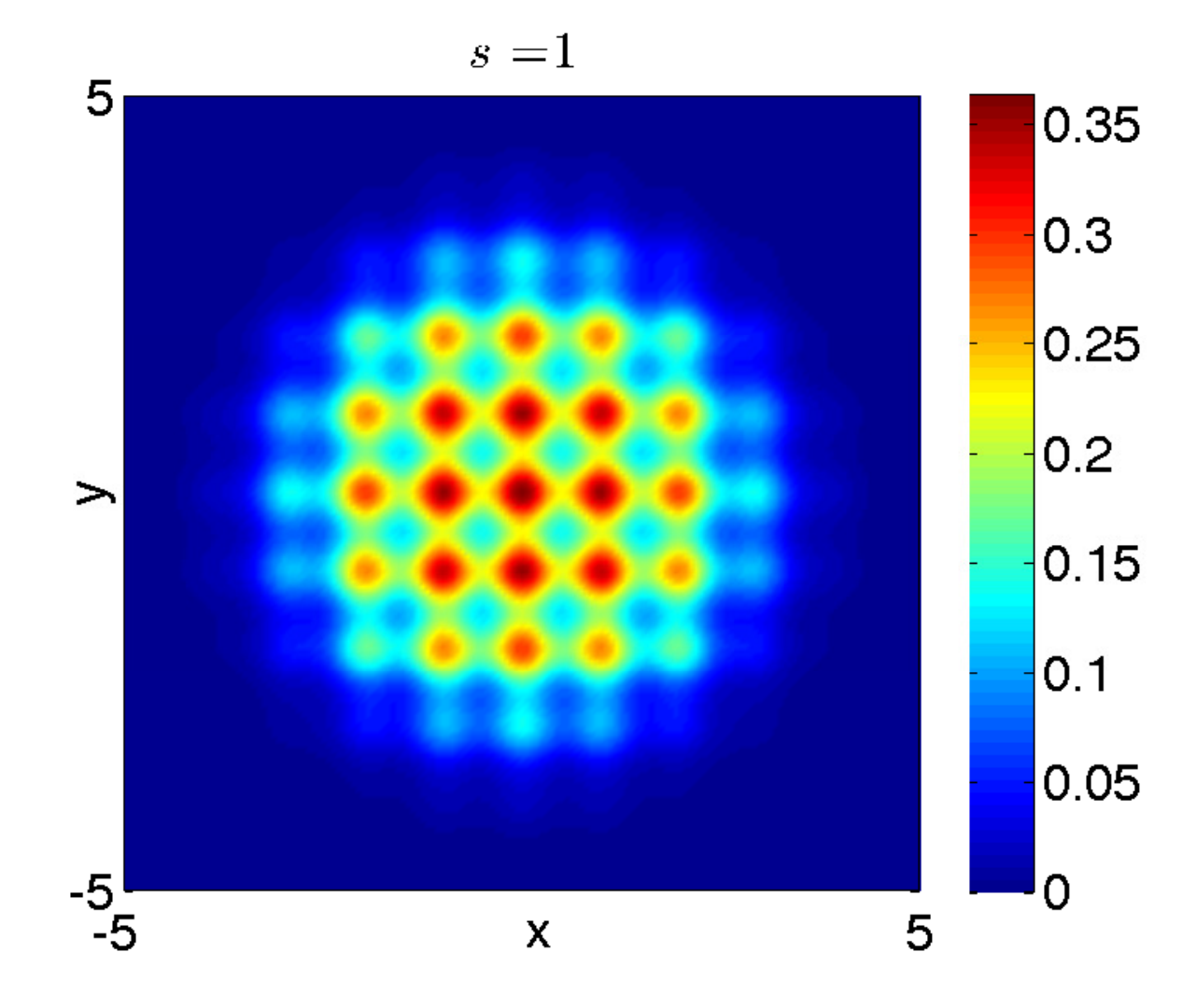,height=3.6cm,width=4.2cm,angle=0}
}
\centerline{
\,\psfig{figure=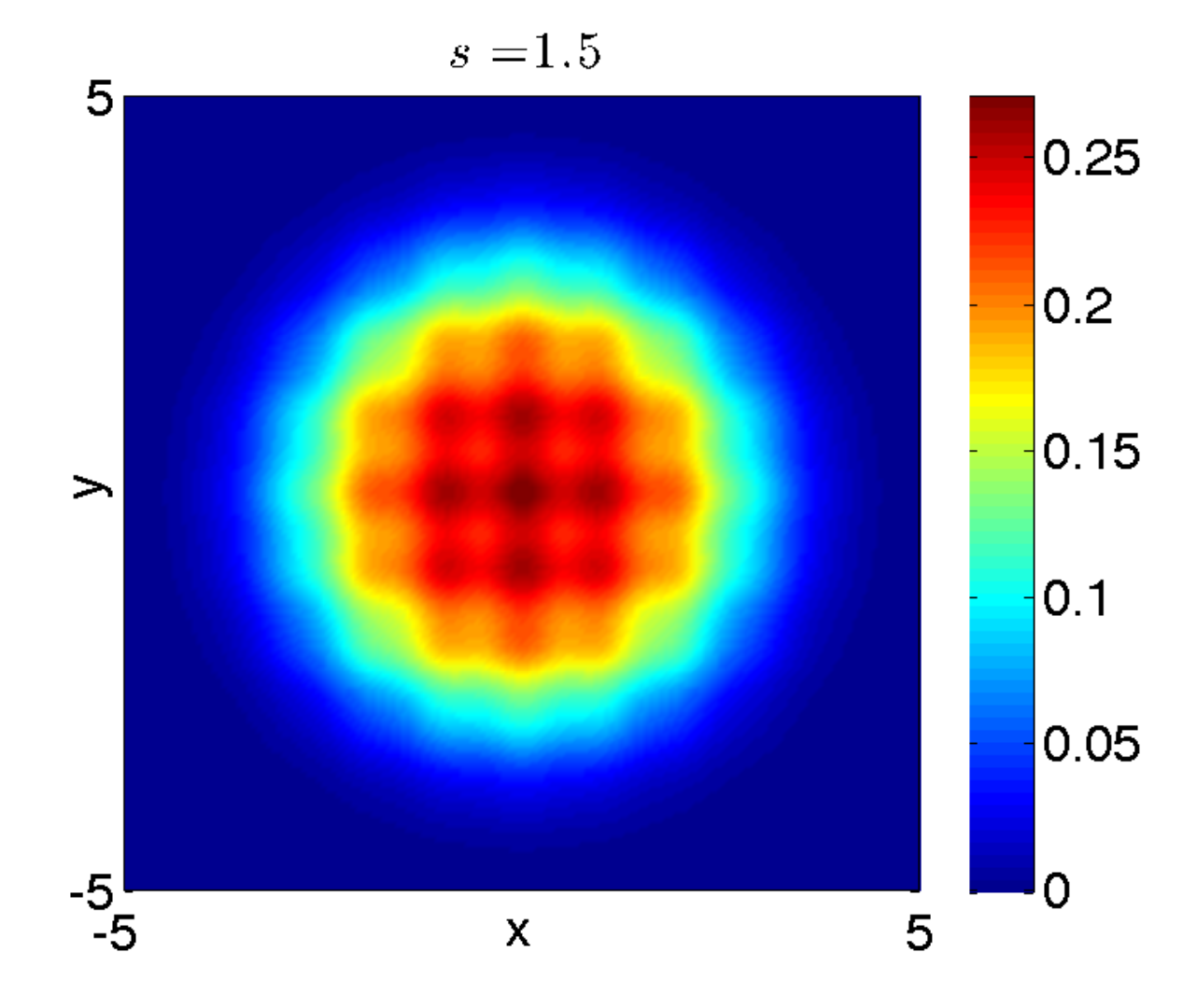,height=3.6cm,width=4.2cm,angle=0}\;\;
\psfig{figure=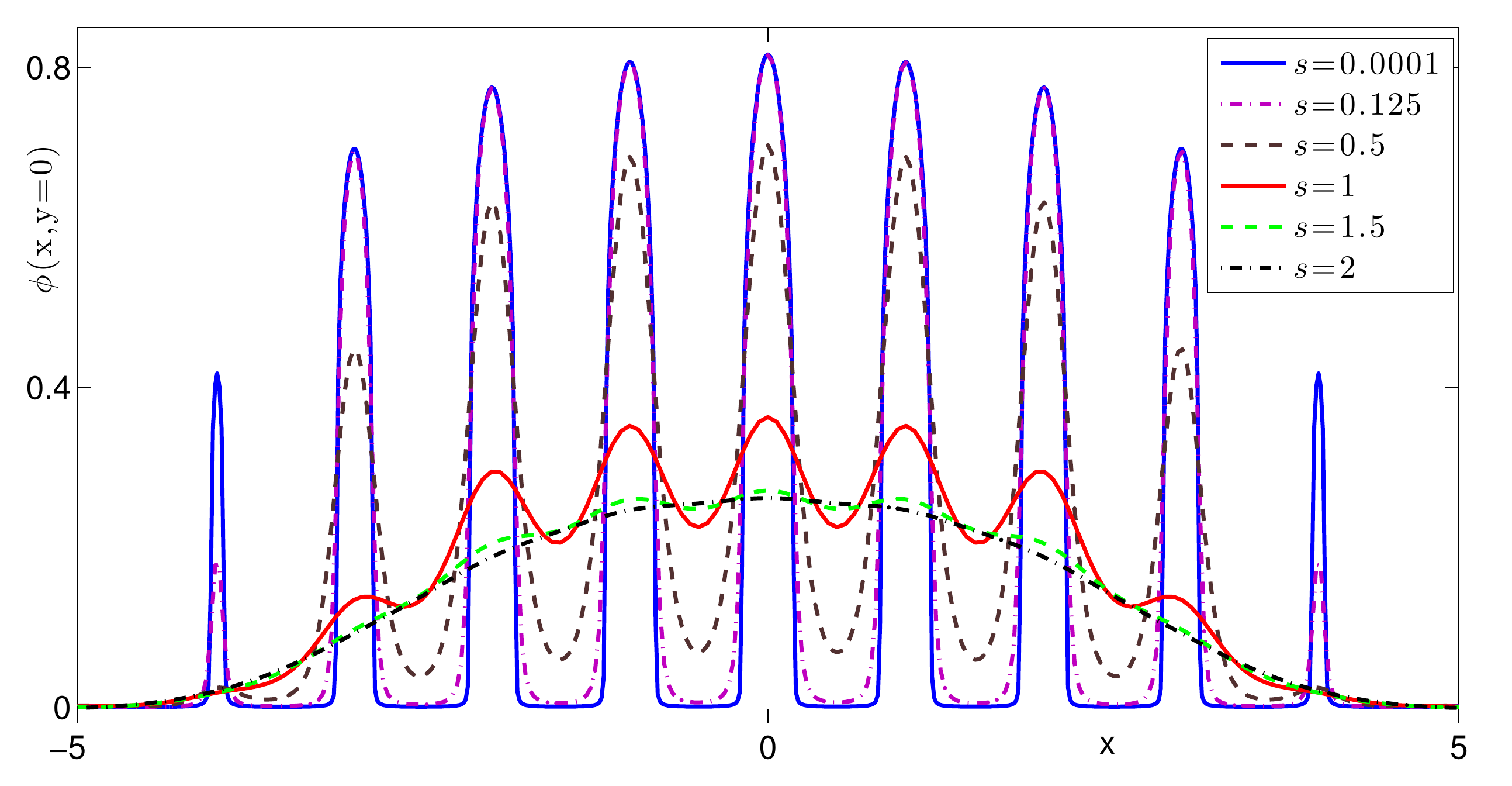,height=3.6cm,width=8.4cm,angle=0}
}
\caption{Contour plots of the density of the ground state $\phi_g (x,y)$  and
 the  slice plot of $\phi_g (x,y=0)$  in example \ref{GS: NonRot2}.}
\label{fig:ex2_den}
\end{figure}

\begin{exmp}\label{GS: NonRot2}
{\bf Non-rotating FNLSE with harmonic + optical lattice potential.}
 Here, we choose $\Og=0$. We consider the ground states of the
FNLSE in a harmonic plus optical lattice potential  with different parameters.  
To this end, we  let $\lambda=64$ and $\beta=0$ and choose the potential as
\[
V(x,y)=\fl{x^2+y^2}{2} + 10 \big(\sin^2(\pi x) +\sin^2(\pi y) \big).
\]
\end{exmp}
The spatial mesh sizes  are chosen as $h_x=h_y=\fl{1}{32}$ in this case. 
Figure \ref{fig:ex2_den} shows the contour plot of the ground state density $\rho_g:=|\phi_g(\bx)|^2$ 
and the slice plot of $\phi_g(x,0)$ with  different fractional orders $s$.

From Figures \ref{fig:ex1}--\ref{fig:ex2_den} and additional results not shown here, we can conclude that 
(i) The ground states become more peaked and
narrower as the fractional order $s$ tends smaller, which corresponds to subdispersion.  
(ii) A large fractional order helps in smoothing out the density profile (cf. Fig. \ref{fig:ex2_den}) for the superdispersion case.
(iii) The repulsive local/nonlocal interactions suppress 
the ``focus" or ``homogenization'' effect as the dispersive order $s$ tends smaller or larger. In other words, the repulsive nonlinear interaction 
helps to stabilize the ground states.
(iv) When $\beta$ and/or $\lambda$ are/is large, the nonlinear interaction dominates and the dispersive effect can be neglected. 
%

\begin{exmp}\label{GS: RotFNLS}
{\bf  Rotating FGPE. }
In this example, we present the ground states of the rotating FGPE with only local nonlinear interaction, i.e. $\lambda=0$
and $\beta=100$.
\end{exmp}

\begin{figure}[h!]
\centerline{
\psfig{figure=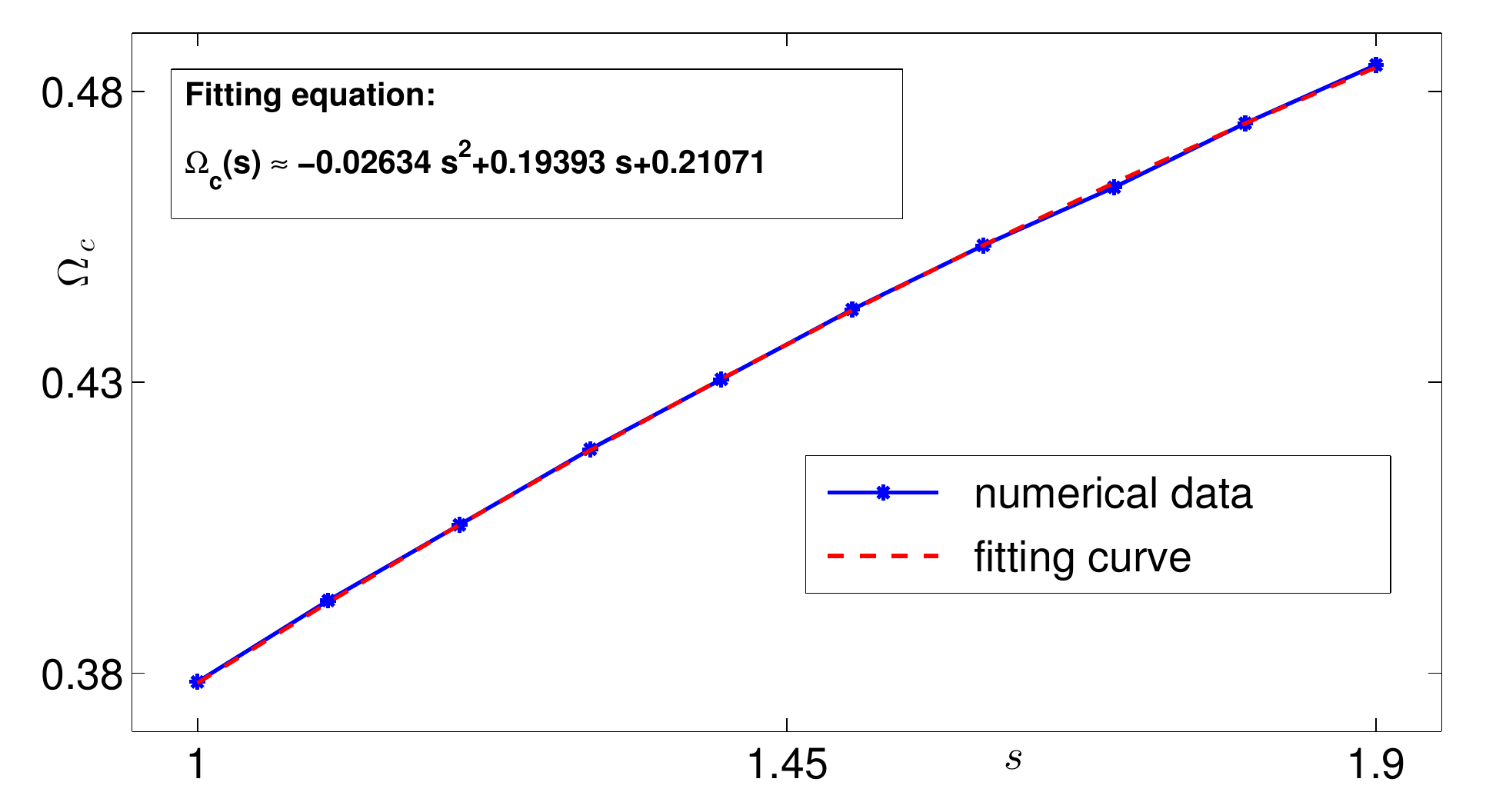,height=6cm,width=12cm,angle=0}}
\caption{Critical rotating frequency vs. the fractional order $s$ in example \ref{GS: RotFNLS}. }
\label{fig:ex3_cri_og}
\end{figure}

\begin{figure}[h!]
\centerline{
\psfig{figure=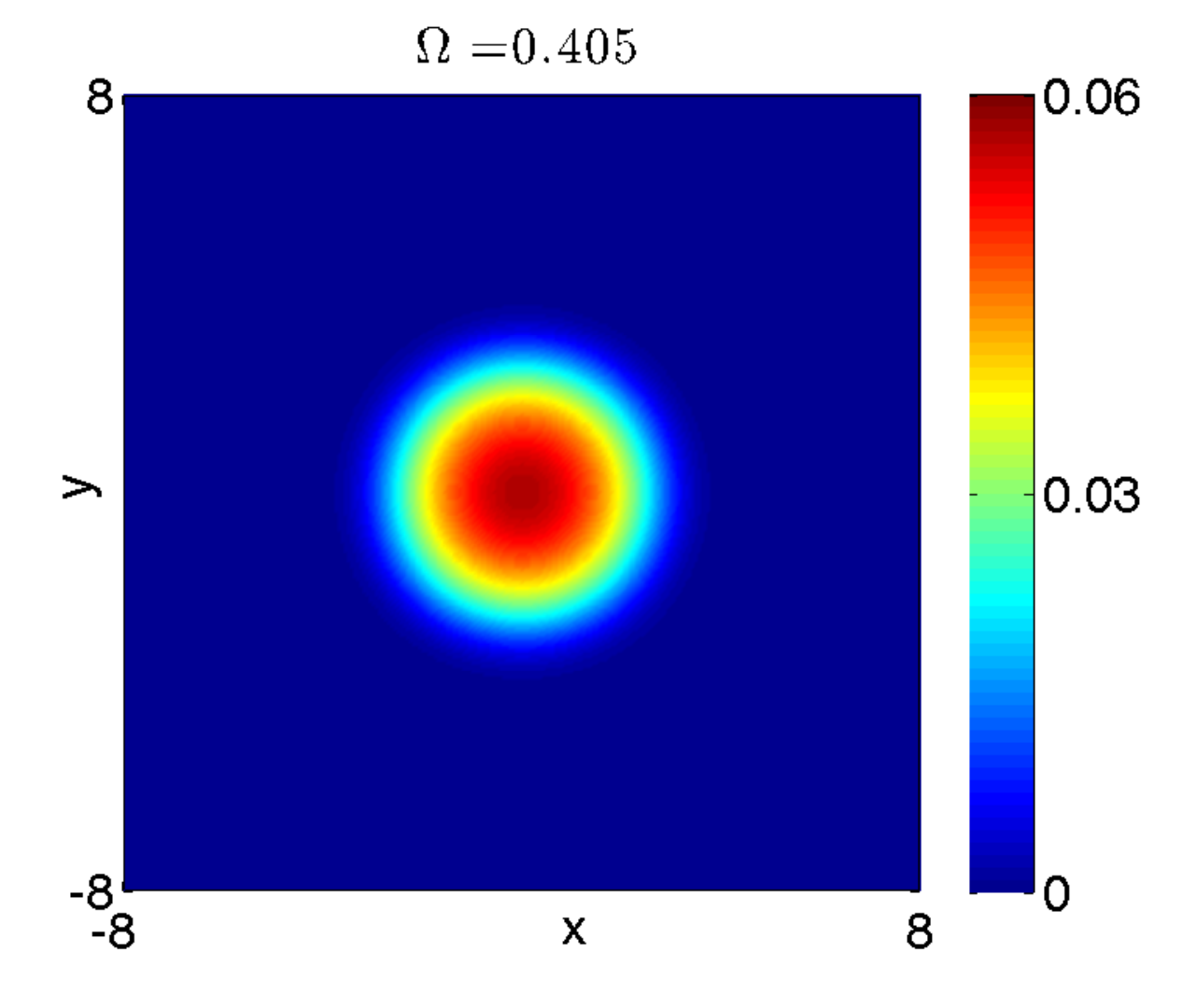,height=3.6cm,width=4.2cm,angle=0}
\psfig{figure=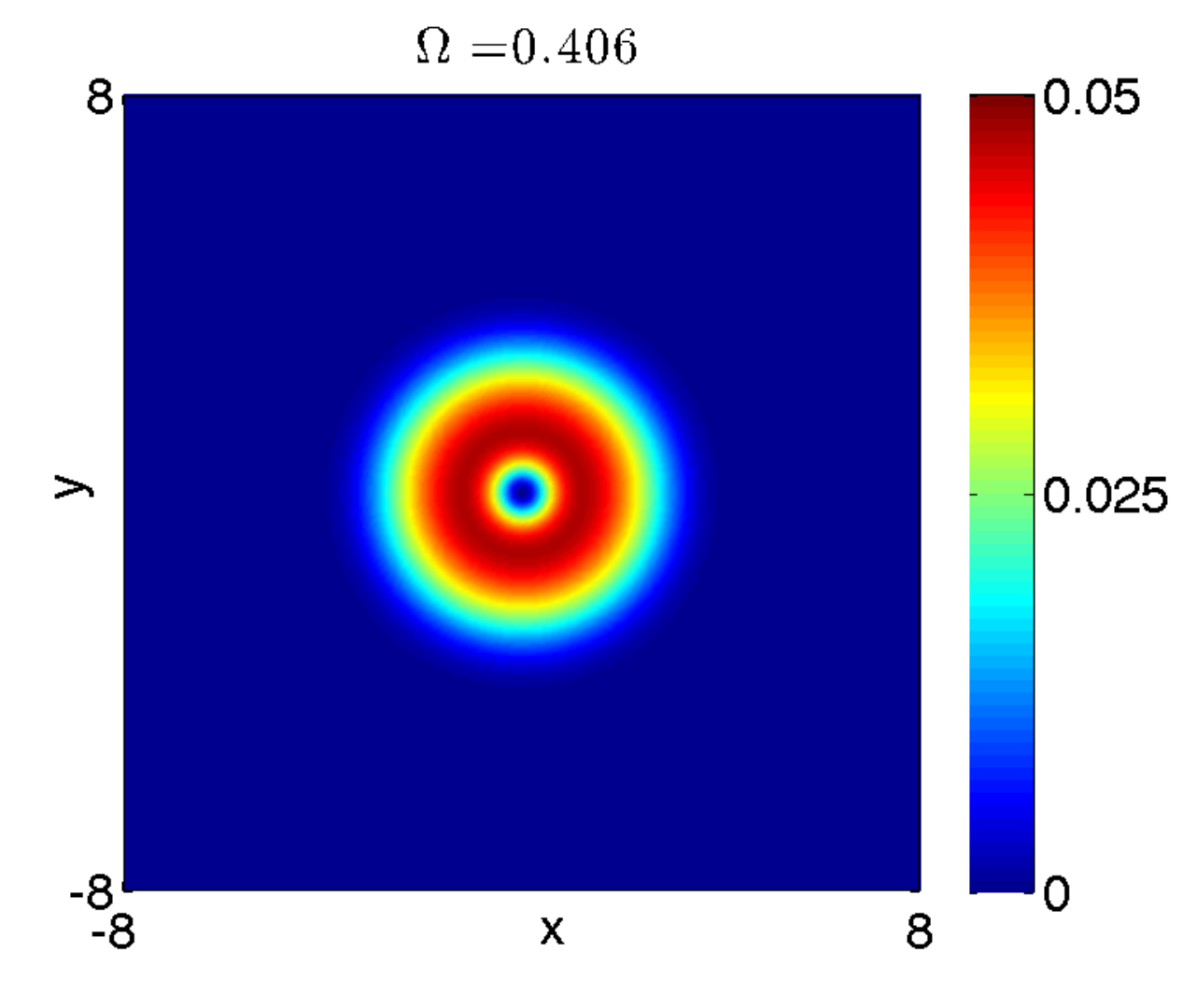,height=3.6cm,width=4.2cm,angle=0}
\psfig{figure=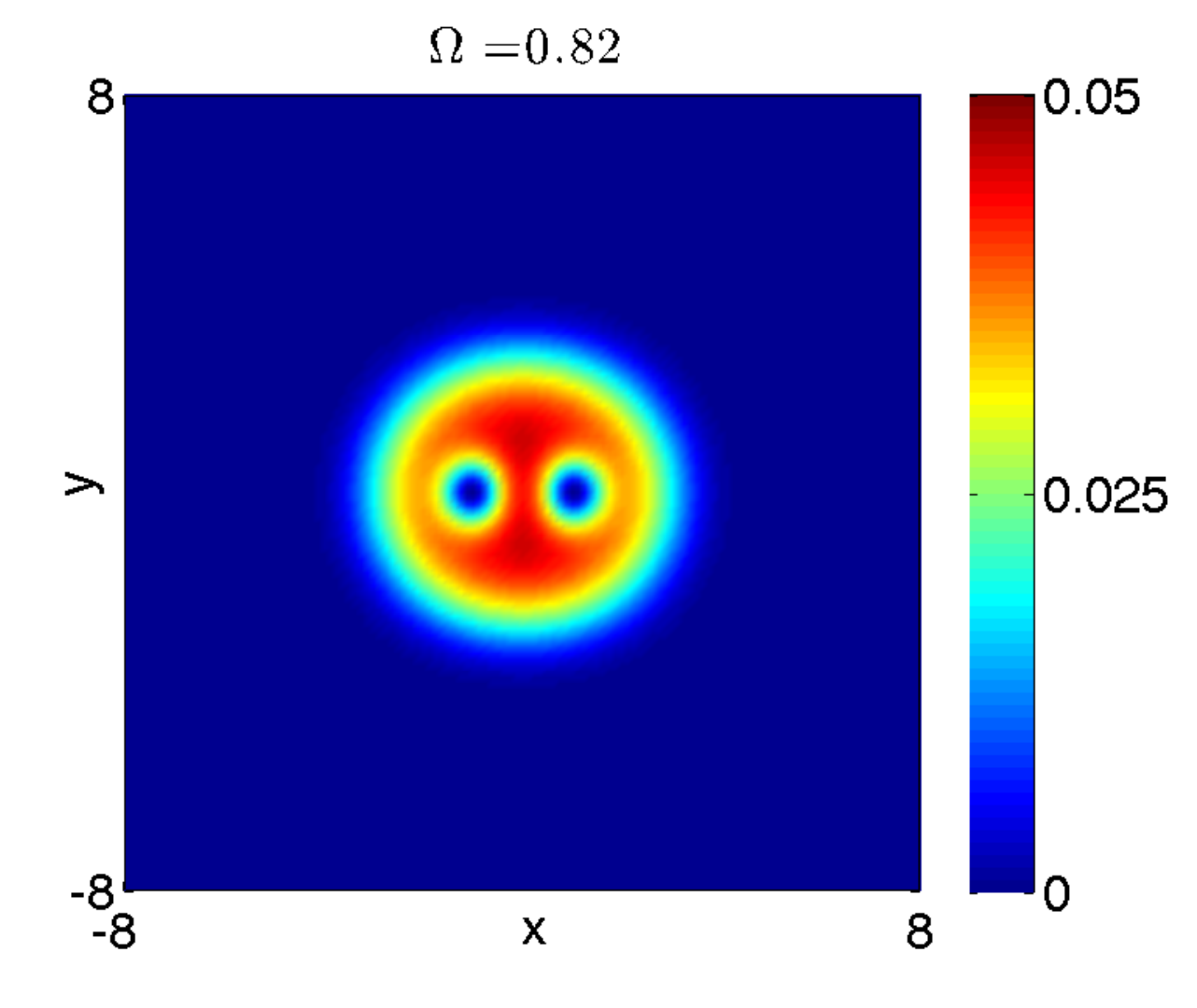,height=3.6cm,width=4.2cm,angle=0}
\psfig{figure=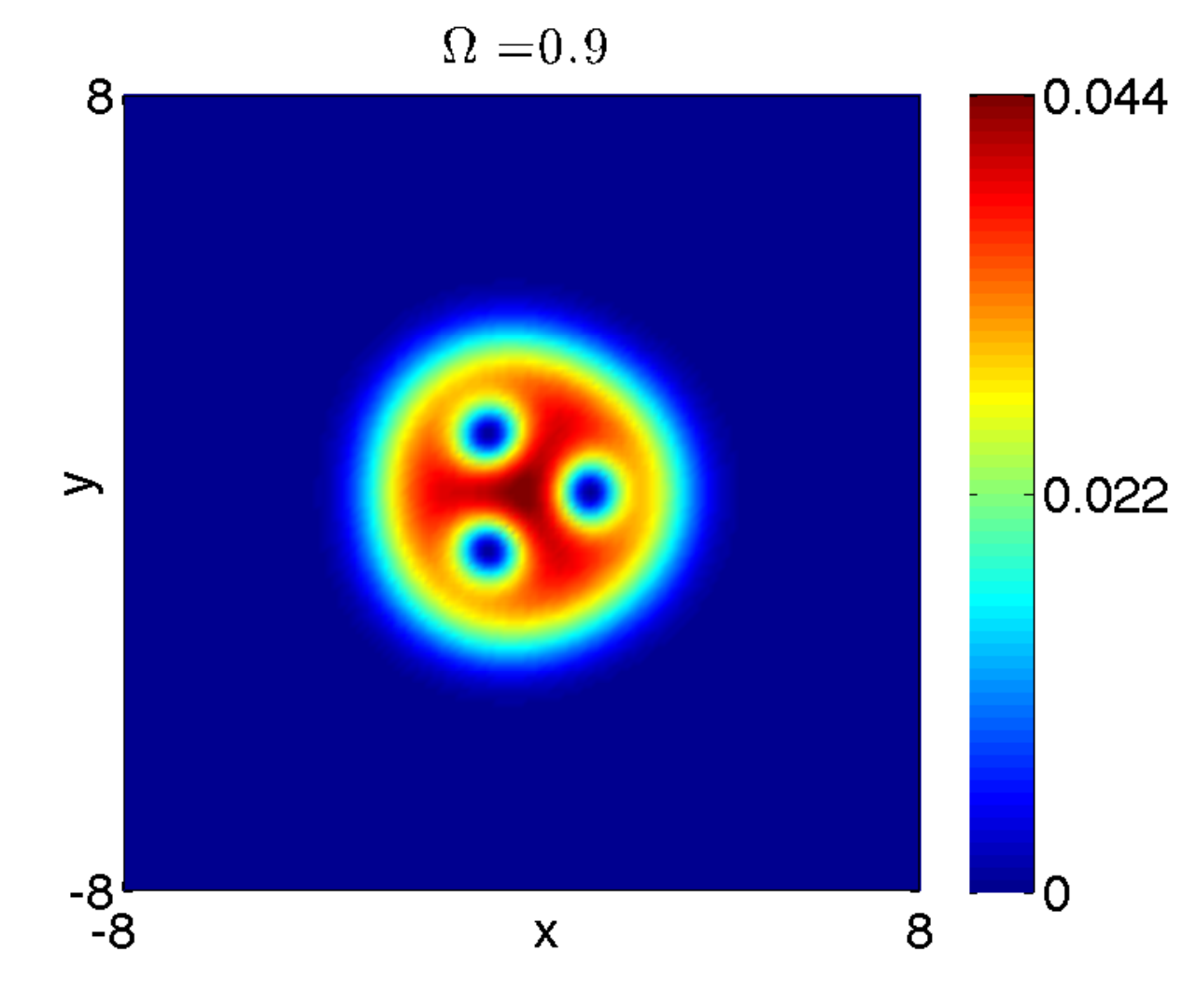,height=3.6cm,width=4.2cm,angle=0}
}
\centerline{
\psfig{figure=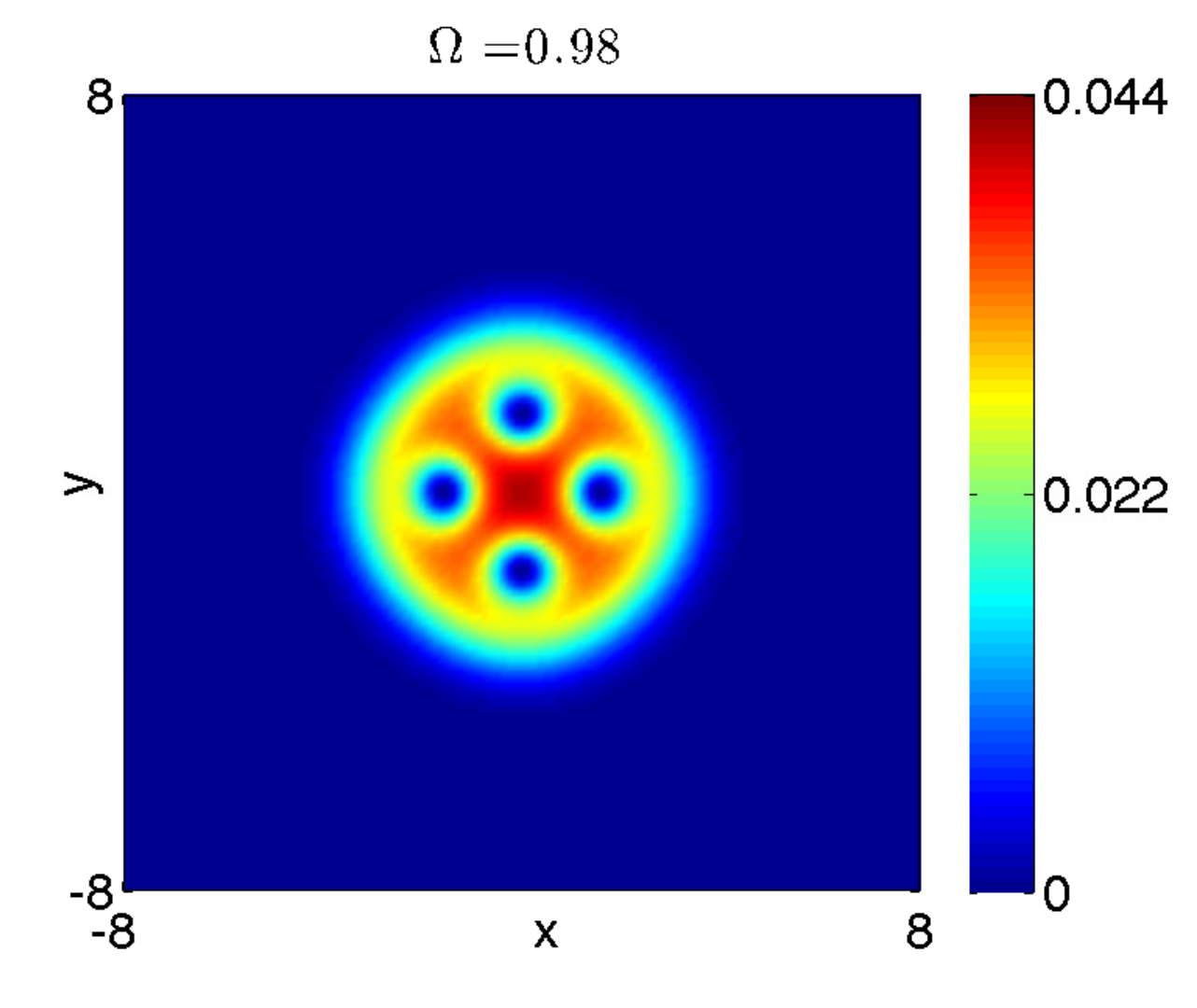,height=3.6cm,width=4.2cm,angle=0}
\psfig{figure=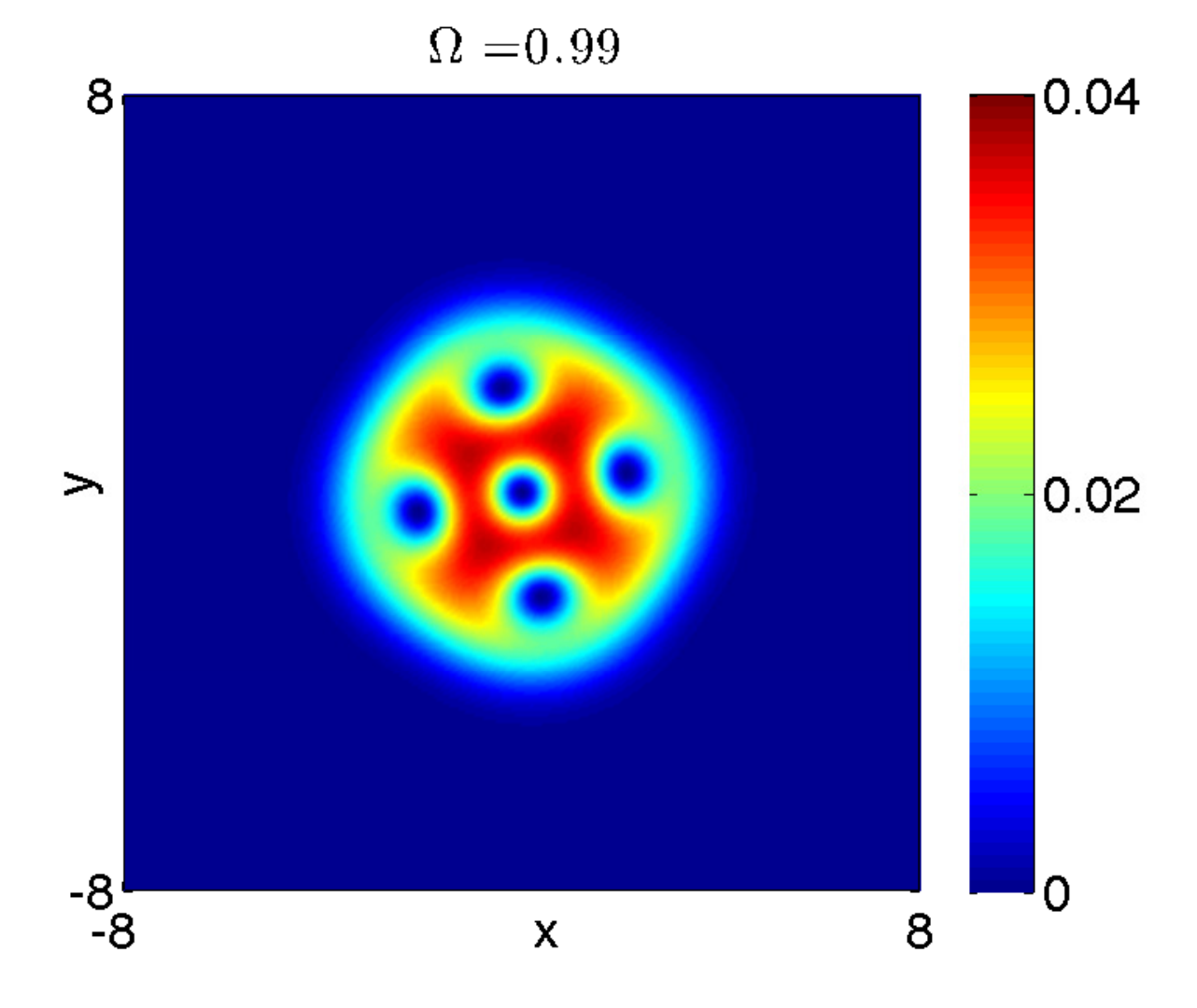,height=3.6cm,width=4.2cm,angle=0}
\psfig{figure=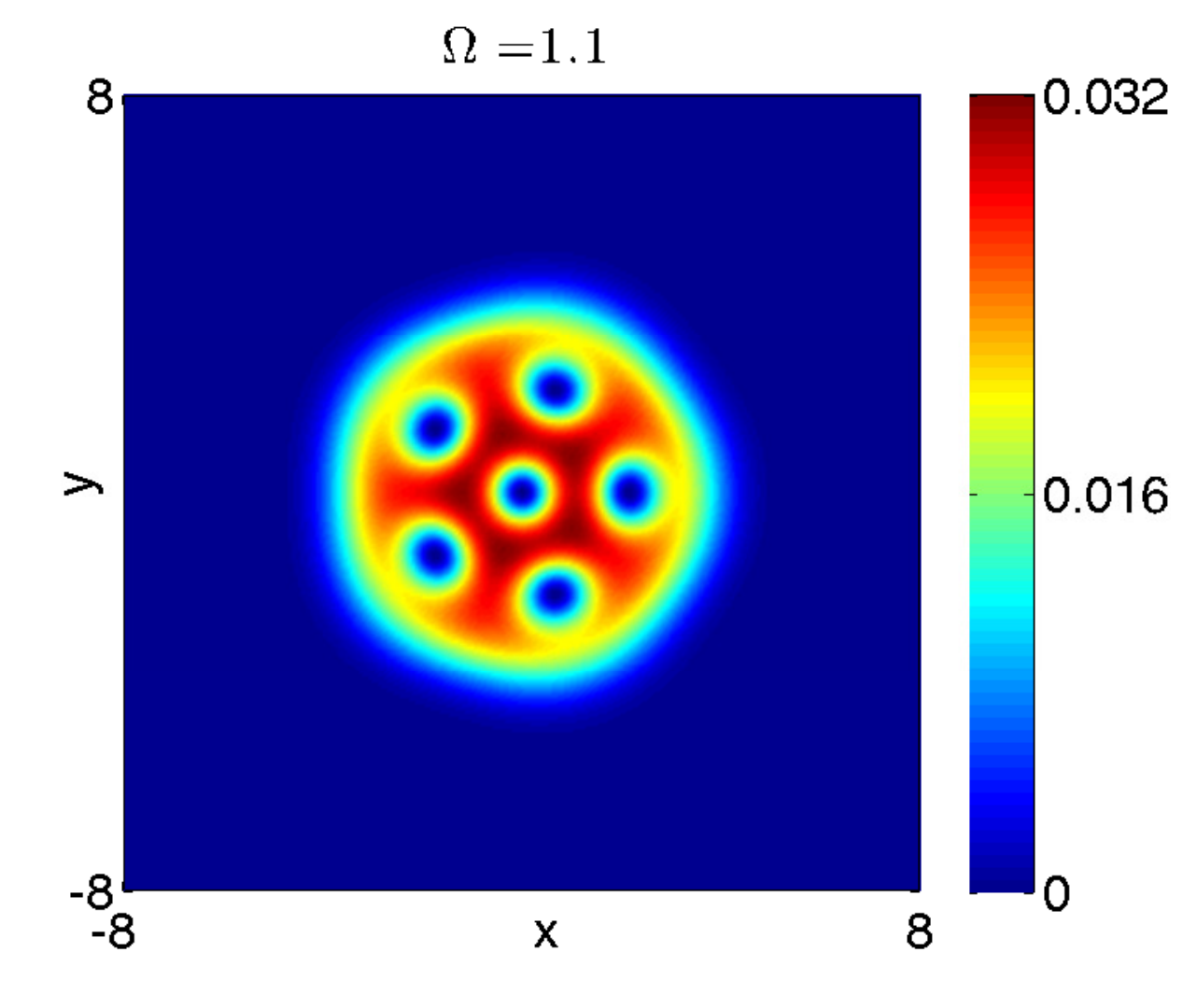,
height=3.6cm,width=4.2cm,angle=0}
\psfig{figure=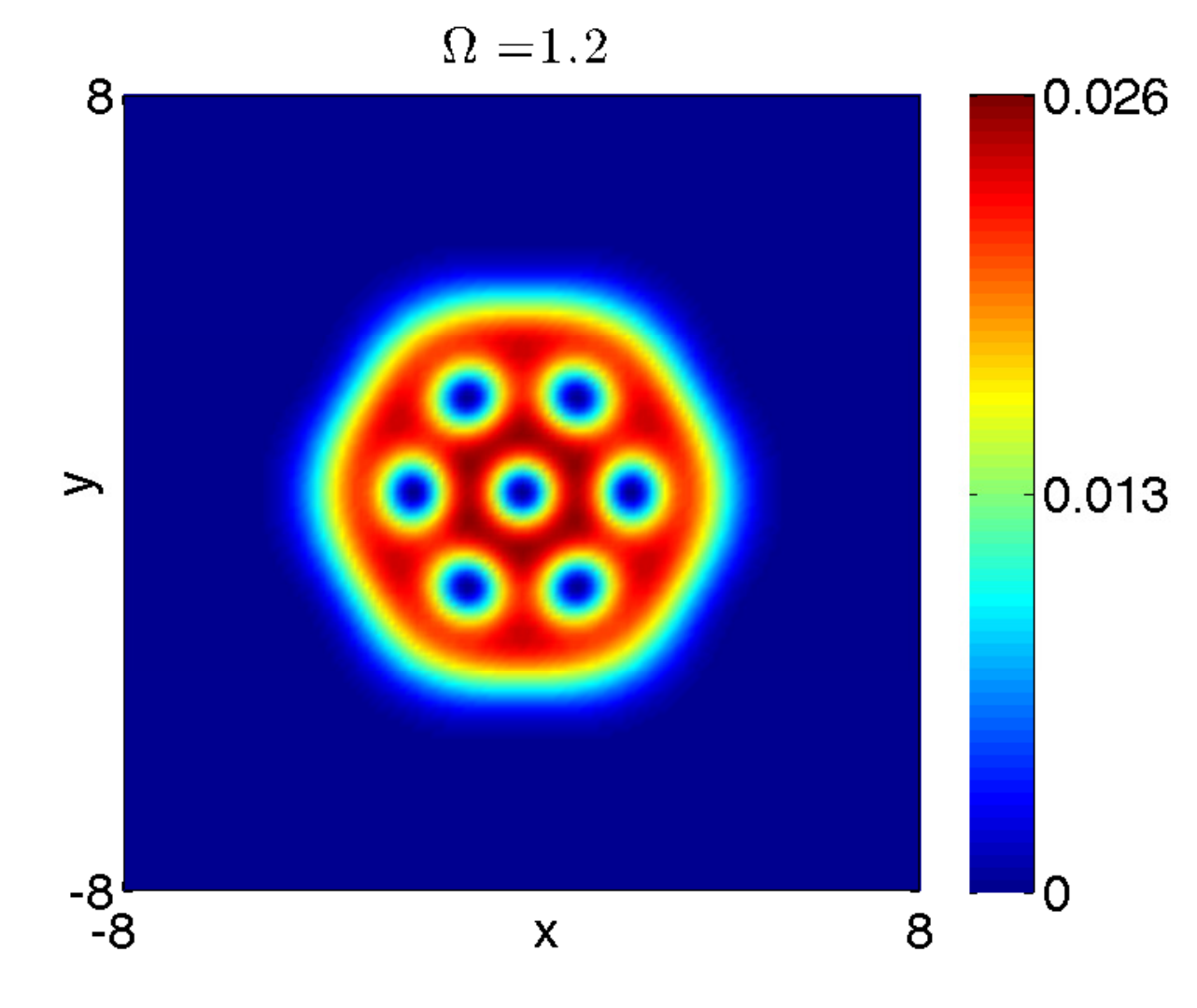,
height=3.6cm,width=4.2cm,angle=0}
}
\centerline{
\psfig{figure=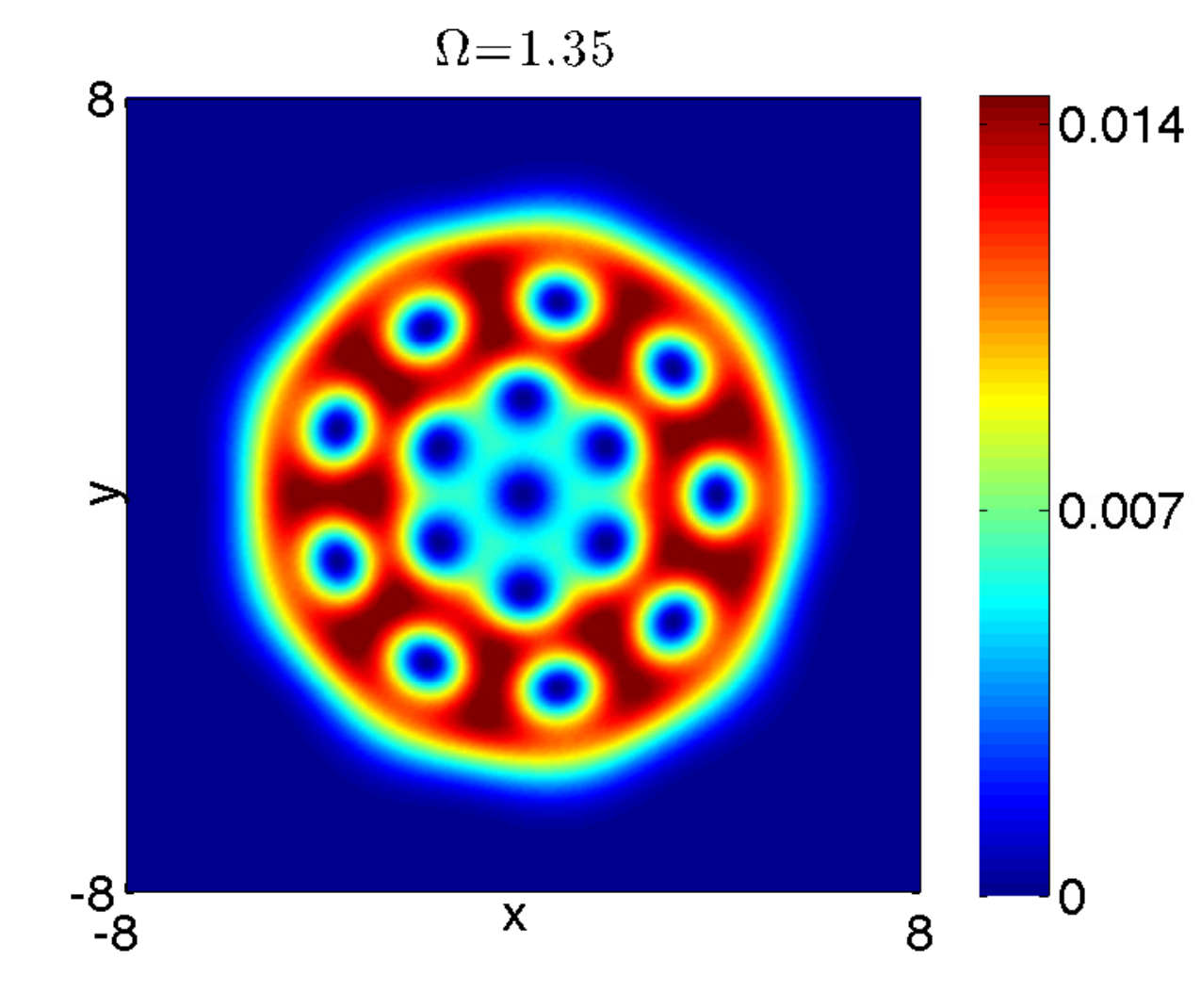,
height=3.6cm,width=4.2cm,angle=0}
\psfig{figure=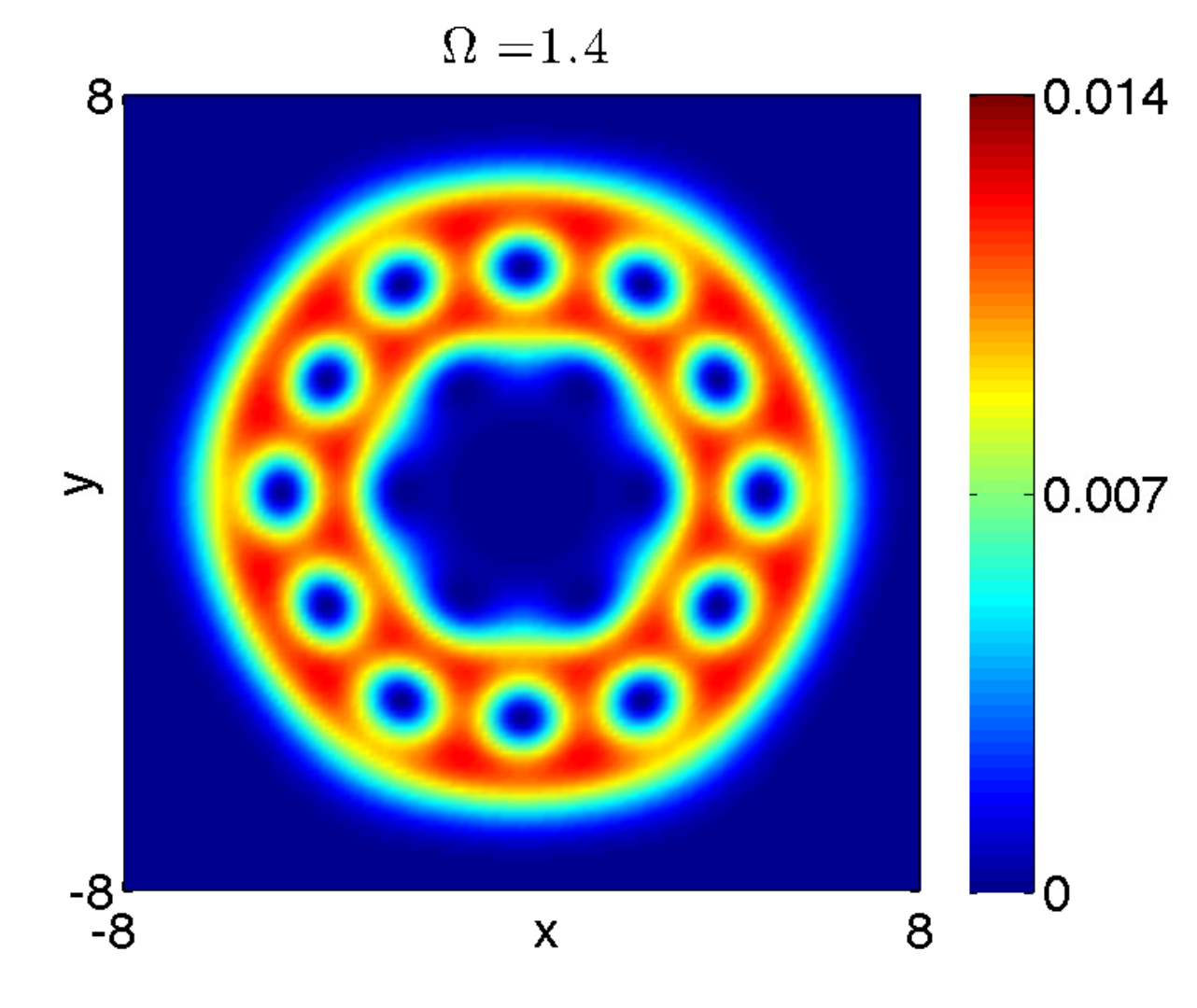,
height=3.6cm,width=4.2cm,angle=0}
\psfig{figure=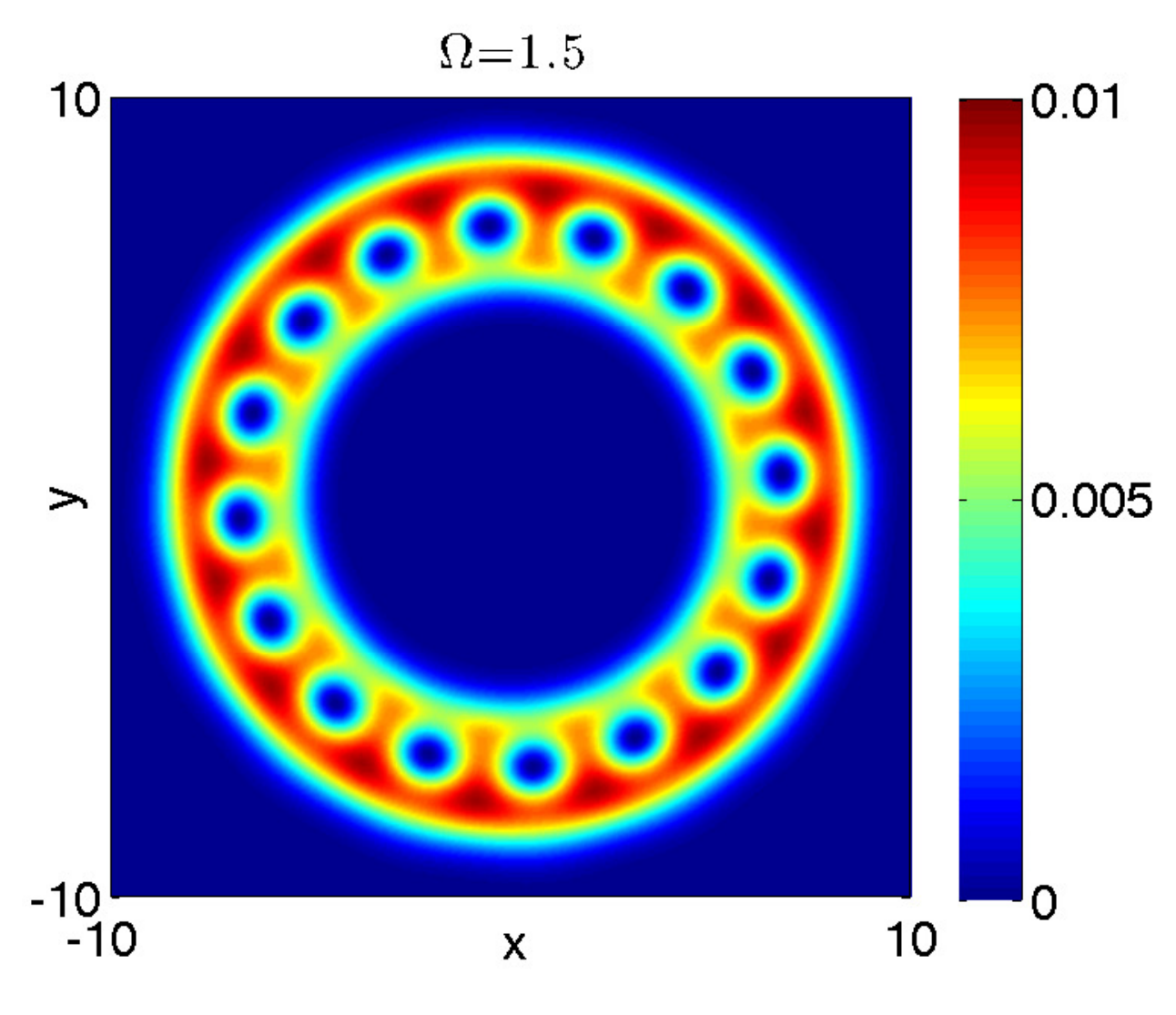,
height=3.6cm,width=4.2cm,angle=0}
\psfig{figure=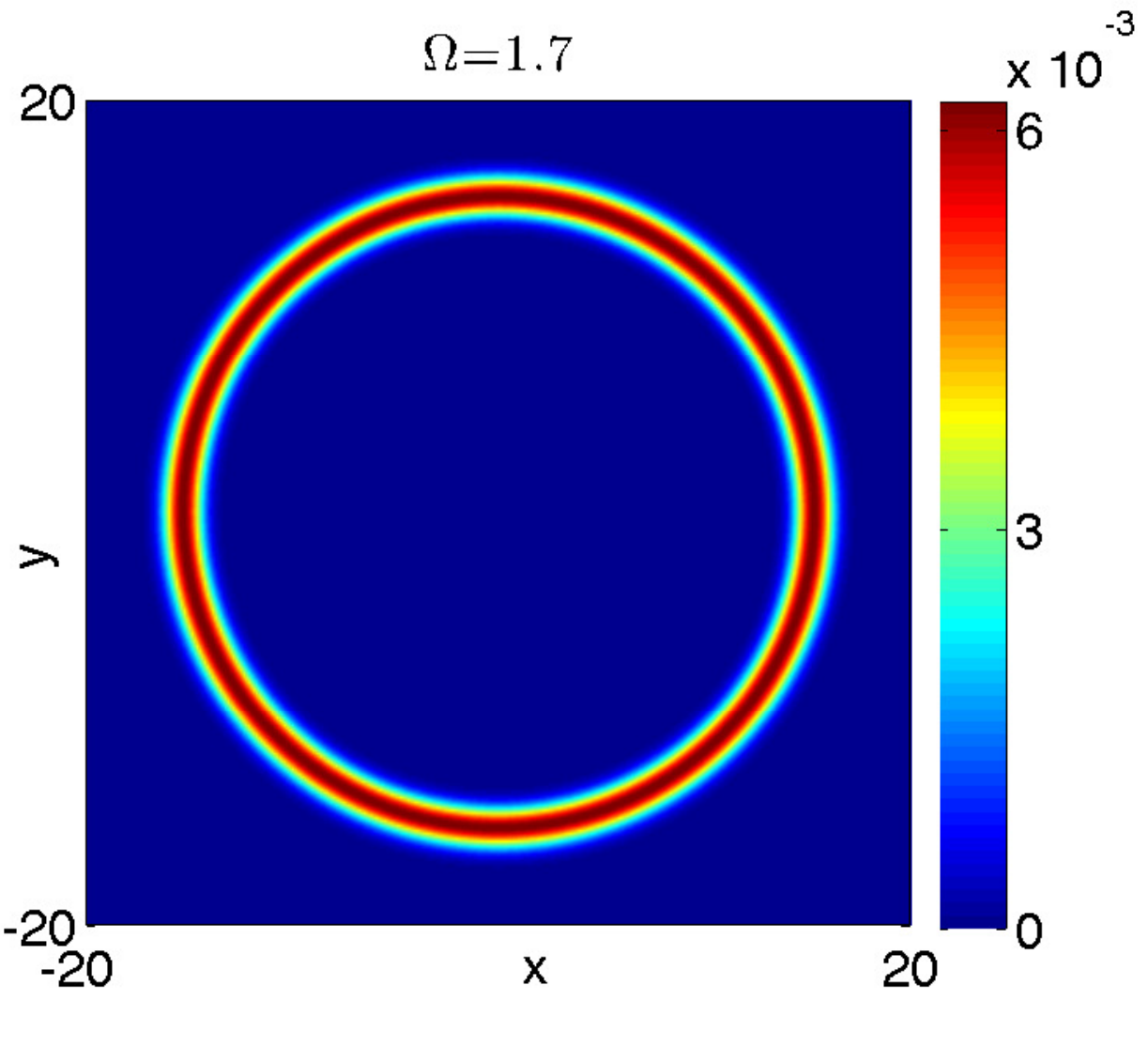,
height=3.6cm,width=4.0cm,angle=0}
}
\caption{Contour plots of the density $|\phi_g(\bx)|^2$ in example \ref{GS: RotFNLS} (superdispersion). }
\label{fig:ex3_diff_Dens_vs_og}
\end{figure}

We propose to numerically study the dependence of the first critical rotating velocity  $\Og_c$ to create a vortex with respect
to the fractional dispersive order $s$.  
Figure \ref{fig:ex3_cri_og} shows this relation  derived by a linear regression
\be
\Og_c(s)\approx-0.02634\;  s^2+0.19393\; s + 0.21071.
\ee
Figure  \ref{fig:ex3_diff_Dens_vs_og} displays the contour plots of the ground state density $\rho_g$ 
for different values of $\Og$ but with  $s=1.2$ (superdispersion).
From Figure \ref{fig:ex3_cri_og}-\ref{fig:ex3_diff_Dens_vs_og} and additional results not shown here, we can conclude that
(i) The first critical rotating velocity  $\Og_c$ depends almost linearly on $s$.
(ii) For the superdispersion case, i.e. $s>1$, the ground states exist for all velocities $\Og$. 
As $\Og$ increases, the ground states will undergo three phase transitions (similar to the non-fractional GPE with 
quartic order trapping potential), i.e., from Gaussian-type to one-vortex profile, from vortex lattice to vortex-lattice with a hole at the
center and then to a giant vortex. 
It would be interesting to study how these critical rotating frequencies for the 
transitions depend on $s$ and  how they compare with those in the case of the standard GPE.

%

\section{Dynamics computation: properties, numerical scheme and simulations}
\setcounter{equation}{0}
\setcounter{table}{0}
In this section, we first present analogous dynamical laws for 
some commonly used quantities in the classical rotating GPE.
Then, we extend the rotating Lagrangian coordinates transform
proposed for the standard GPE \cite{BMTZ2013} to the FGPE. In the rotating Lagrangian coordinates, 
the rotation term vanishes, giving rise to a time-dependent potential.
Based on the new FNLSE, we propose a time-splitting Fourier pseudo-spectral method
incorporated with the GauSum solver to simulate the dynamics.

\subsection{Dynamical properties}
\label{sec: dyn_prop}

Here we study the dynamical properties of  the mass, energy, angular momentum expectation and center of 
mass \cite{BMTZ2013}. The dynamical laws can be used as benchmarks to test the numerical methods and are briefly listed here.  
For details, one can either refer to  appendices or to \cite{Tang2013} for   analogous proofs to their non-fractional counterparts.

\noindent\textbf{Mass and energy.} 
The FNLSE (\ref{SFSE})-(\ref{NonLocalPot}) conserves  the mass (\ref{mass}) and the energy (\ref{energy}), i.e.
\be\label{MassEnerg}
\mathcal{N}(t)=\mathcal{N}(t=0),\qquad \qquad \mathcal{E}(t)=\mathcal{E}(t=0).
\ee

\noindent
{\bf Proof:}
It is straightfoward to prove in a similar way as their non-fractional counterparts \cite{Tang2013}  by using the
Plancherel's formula.  
\hfill $\square$

\

\noindent\textbf{Angular momentum expectation.}   The {\it angular momentum expectation}
is defined  as
\be\label{AME}
\langle L_z\rangle(t) =  \int_{{\Bbb R}^d} \bar{\psi}(\bx,t)  L_z \psi(\bx,t)\,d\bx,\qquad t\geq 0.
\ee

\begin{lemma}
\label{lawsDyn_AME}
The angular momentum expectation  $\langle L_z\rangle(t) $ satisfies the following equation 
\be\label{AME_Law}
\fl{d }{d t}\langle L_z\rangle(t)=\int_{\mathbb{R}^{d}} |\psi|^2 (y\p_x-x\p_y) \big(V(\bx)+\lambda\Phi(\bx,t)\big) d\bx.
\ee
This implies that the angular momentum expectation is conserved, i.e.
\be
\langle L_z\rangle(t) =\langle L_z\rangle(0),\qquad t\ge0, 
\ee
when $V(\bx)$ is radially/cylindrically symmetric in 2D/3D and one of the  following conditions holds: (i) $\lambda=0$,
(ii) $\lambda\ne0$, $\Phi(\bx)$ is the Coulomb potential or (iii)  $\lambda\ne0$, $\Phi(\bx)$
is the dipole potential with dipole axis $\bn=(0, 0, 1)^T$, i.e. is parallel to the $z$-axis.
\end{lemma}

\noindent
{\bf Proof.}
Details of the proof are given in   \ref{Appendix_A}. 
\hfill $\square$

\

\noindent\textbf{Center of mass.} The center of mass is defined by 
\bea
\label{CoM}
\bx_c(t) = \int_{{\mathbb R}^d} \bx\, |\psi(\bx,t)|^2 d\bx= \langle \bx\psi, \psi  \rangle.
\eea

\begin{lemma}
\label{lawsDyn_COM}
The center of mass $\bx_c(t)$ satisfies the following equations, for $0<s\le 1$ (subdispersion),
\bea
\label{CoMLaw1}
&&\dot{ \bx}_c -\Og J \bx_c
	= i\big\langle G\ast\psi, \nabla\psi \big\rangle,  \\[0.5em]
&&\ddot{\bx}_c -2\Og J \dot{\bx}_c+\Og^2J^2\bx_c
	= 2 {\rm Re}\Big(  \big \langle G\ast (\mathcal{V}\psi),   \nabla\psi \big \rangle\Big).
\label{CoMLaw2}
\eea
Here, we set $\mathcal{V}(\bx,|\psi|)=V(\bx)+\beta|\psi|^2+\lambda\Phi(\bx,t)$, and
\bea
\label{matrix1}
{J} = \left(\begin{array}{cc}
0 & 1  \\
-1 & 0
\end{array}\right), \quad
\mbox{for\quad $d = 2$},
\qquad
{J} = \left(\begin{array}{ccc}
0 & 1 & 0 \\
-1 & 0 & 0 \\
0 & 0 & 0
\end{array}\right), \quad
\mbox{for\quad $d = 3$}.
\eea
The convolution kernel $G(\bx)$ reads as
\be\label{Gx_theor}
G(\bx)=\left\{
\begin{array}{cl}
\delta(\bx),	& \qquad	s=1, \\[0.5em]
\fl{2^{s-d/2}s}{\Gamma(1-s )\, \pi^{d/2}}
\left(\fl{m}{|\bx|}\right)^{\fl{d}{2}+s-1}K_{\fl{d}{2}+s-1}\Big(m |\bx|\Big), 
&	\qquad 0<s<1,
\end{array}
\right.
\ee
where $\delta(\bx)$ is the Dirac delta function and $K_{v}(z)$, the modified Bessel function of the second-kind
and order $v$, is given explicitly as follows 
\be\label{Bes2}
K_{v}(z)=\fl{(2 z)^v \Gamma(v+\fl{1}{2})}{\sqrt{\pi}} \int_{0}^{\infty}\fl{\cos(t)}{(t^2+z^2)^{v+\fl{1}{2}}} dt.
\ee
\end{lemma}

\noindent
{\bf Proof.}
A detailed proof is reported in   \ref{Appendix_B}. 
\hfill $\square$

\begin{remark}
If $s=1$, $V(\bx)$ is the harmonic potential  (\ref{harm_poten}) and $\Phi(\bx)$  
is the Coulomb potential or DDI with $\bn=(0, 0, 1)^T$, 
then (\ref{CoMLaw2}) reduces to  \cite{BMTZ2013, KZ2014}
\be
\label{classical_CoM}
\ddot{\bx}_c -2\Og J \dot{\bx}_c+ ( \Og^2 J^2 + \Lambda_d)\bx_c={\bf 0},
\ee
where 
\bea
\label{matrix2}
\Lambda_d = \left(\begin{array}{cc}
\gm^2_x & 0  \\
0 & \gm^2_y
\end{array}\right), \quad
\mbox{for\quad $d = 2$},
\qquad
\Lambda_d= \left(\begin{array}{ccc}
\Lambda_2 & {\bf 0}  \\
 {\bf 0}  & \gm^2_z
\end{array}\right), \quad
\mbox{for\quad $d = 3$}.
\eea
In \cite{KZ2014}, the authors derived a dynamical law for the center of mass for the FNLSE with $s\in(\fl{1}{2}, 1]$ 
and for a harmonic trapping potential. Compared with their results, the dynamical laws (\ref{CoMLaw1})-(\ref{CoMLaw2})
are simpler and hold for a general potential $V(\bx)$ as well as for the full subdispersion case, i.e. $\forall$  $s\in(0, 1].$
It is also interesting to explore similar equations for the superdispersion case $s>1$.

\end{remark}

\begin{remark}
We also remark here that it might be interesting to derive the dynamical laws for the condensate width $\delta_v$ which is defined as
\be\label{ConWidth}
\delta_v(t)=\int_{\mathbb{R}^d}v^2 |\psi(\bx)|^2d\bx,\qquad v=x,y\ {\rm in\ 2D} \ {\rm and} \ v=x,y,z\ {\rm in\ 3D}.
\ee
The derivation and proof is  feasible but tedious. One can refer to \cite{Tang2013} for the analogous details.
%
%
\end{remark}

\subsection{Numerical  method}
In this subsection, we first introduce a coordinates transformation and 
reformulate the rotating FGPE (\ref{SFSE})-(\ref{NonLocalPot})  in
the new coordinates,  eliminating hence the rotation term.

\subsubsection{Rotating Lagrangian coordinates transformation}

For any time $t\geq 0$, let ${ A}(t)$ be the orthogonal rotational matrix  defined as \cite{BMTZ2013}
\bea\label{Amatrix}
{ A}(t)=\left(\begin{array}{cc}
\cos(\Omega t) & \sin(\Omega t) \\
-\sin(\Omega t) & \cos(\Omega t)
 \end{array}\right),  \quad  \mbox{if \ \ $d = 2$,}  \quad  \
{A}(t)=\left(\begin{array}{ccc}
         \cos(\Omega t) & \sin(\Omega t) & 0 \\
         -\sin(\Omega t) & \cos(\Omega t) & 0 \\
         0 & 0  & 1
         \end{array}\right), \quad\ \mbox{if \ \ $d = 3$.} \quad
\eea
It is easy to check that $A^{-1}(t) =A^T(t)$ for any $t\ge0$ and ${A}(0) = { I}$, where
 $I$ is the identity matrix. For any $t\ge0$, we introduce the {\it rotating Lagrangian coordinates}
$\tbx$ as \cite{AMS2012,BMTZ2013,GPV2001}
\bea\label{transform}
\tbx={A}^{-1}(t) \bx=A^T(t)\bx \quad \Leftrightarrow \quad \bx= {A}(t){\tbx},   \qquad \bx\in {\mathbb R}^d,
\eea
and we denote by $\phi:=\phi(\tbx, t)$ the wave function in the new coordinates
\bea\label{transform79}
\phi(\tbx, t):=\psi(\bx, t)= \psi\left({A}(t){\tbx},t\right), \qquad \bx\in {\mathbb R}^d, \quad t\geq0.
\eea
By some simple calculations, one can easily obtain 
\bea
\label{operator1}
&& \p_t\phi(\tbx,t) =\p_t\psi(\bx, t)  + \nabla\psi(\bx,t)\cdot\left(\dot{A}(t) \tbx\right) = \p_t\psi(\bx,t)
- \Og(x\p_y-y\p_x)\psi(\bx,t),\\
\label{operator2}
&& (-\nabla^2+m^2)^s\psi(\bx,t)=(-\nabla^2+m^2)^s\phi(\tbx,t).
\eea
Plugging them back into  (\ref{SFSE})-(\ref{NonLocalPot})  gives the  following
FNLSE in the rotating Lagrangian coordinates
\bea\label{FSE1_Rot}
\ i\p_t \phi(\widetilde{\bx}, t)& =& \left[\fl{1}{2}(-\nabla^2+m^2)^s + \mathcal{W}(\tbx, t) + \bt|\phi|^2
+\lambda \widetilde{\Phi}(\tbx,t) \right]\phi(\widetilde{\bf x},t), \quad \tbx\in{\mathbb R}^d, \quad t>0,\\
\label{NonLoc_Rot}
\widetilde{\Phi}(\tbx,t) &=&\widetilde{\mathcal{U}}*|\phi|^2,\qquad \tbx\in{\mathbb R}^d, \quad t\geq0.
\eea
Here, 
$\mathcal{W}(\tbx,t)=V(A(t)\tbx)$ and
$\widetilde{\mathcal{U}}(\tbx,t)$ reads as 
\be\label{Kernel_Rot}
\widetilde{\mathcal{U}}(\tbx,t)=
\left\{
\begin{array}{ll}
\fl{1}{2^{d-1}\pi |\tbx|^{\mu}}, \qquad 0<\mu<d-1, &	\quad {\rm Coulomb,}\\[0.5em]
-\delta(\tbx)-3\,\partial_{\bbm(t)\bbm(t)}   
\left( \fl{1}{4\pi|\tbx|} \right), &\quad {\rm 3D \ \ DDI,} \\[0.5em]
-\fl{3}{2} \left(\partial_{\bbm_{\perp}(t)\bbm_{\perp}(t)}-m_3^2 \nabla_{\perp}^2\right)\left( \fl{1}{2\pi|\tbx|} \right),
 & \quad {\rm 2D\ \ DDI},
\end{array}
\right.
\ee
with ${\bf m}(t) \in {\mathbb R}^3$  defined as ${\bf m}(t) ={A}^{-1}(t)\bn =:\big((m_1(t),m_2(t),m_3(t)\big)^T$
and ${\bf m}_\perp(t)=\big(m_1(t),m_2(t)\big)^T$.

We can clearly see  that the   rotation term vanishes in the new coordinates (see \eqref{FSE1_Rot}).
Instead, the trapping potential and the dipole axis become time-dependent. The absence of the rotating term allows us
to develop a simple and efficient time-splitting scheme. 

%

\subsubsection{A time-splitting pseudo-spectral method}
Here we shall consider the new equation (\ref{FSE1_Rot})-(\ref{NonLoc_Rot}) which has been reformulated in rotating Lagrangian coordinates.
In a practical computation, we first truncate the problem into 
a bounded computational domain $ {\textbf B} = [L_\tx, R_\tx]\times[L_\ty, R_\ty]\times[L_\tz, R_\tz]$
if $d=3$,  or $ {\textbf B} = [L_\tx, R_\tx]\times[L_\ty, R_\ty]$ if $d=2$. From $t=t_n$
to $t=t_{n+1}:=t_n+\Delta t$, the equation is solved in two steps. One first considers
\be\label{1step}
i\p_t\phi(\tbx, t) = \frac{1}{2} (-\nabla^2+m^2)^{s}\phi(\tbx,t), \qquad  \quad
\tbx\in{\textbf B}, \qquad t_n\leq t\leq t_{n+1},
\ee
with periodic boundary conditions on the boundary $\p {\textbf B}$
for a time step $\Dt t$, then  solves
\bea
\label{2stepA}
 i\p_t \phi({\tbx}, t)&=&\left[ \mathcal{W}(\tbx,t) + \beta |\phi|^{2} +
\lambda\widetilde{\Phi}(\tbx,t) \right]\phi(\tbx,t),  \qquad \tbx\in {\textbf B},
\quad t_n\leq t\leq t_{n+1}, \qquad  \\
\label{2stepB}
\widetilde{\Phi}(\tbx,t)&=&\big(\widetilde{\mathcal{U}}\ast \widetilde{\rho}\big) (\tbx, t), \qquad\qquad\quad
 \tbx\in {\textbf B}, \quad t_n\leq t\leq t_{n+1},
\eea
for the same time step. Here, $\widetilde{\rho}(\tbx,t)=|\phi(\tbx,t)|^2$ if $\tbx\in {\textbf B}$ and  $\widetilde{\rho}(\tbx,t)=0$
otherwise. The linear subproblem (\ref{1step}) is discretized in space by the Fourier pseudo-spectral
method  and integrated in time exactly in the phase space. The nonlinear subproblem
(\ref{2stepA})-(\ref{2stepB})  preserves the density pointwise, i.e. $|\phi(\tbx,t)|^2\equiv|\phi(\tbx,t=t^n)|^2=|\phi^n(\tbx)|^2$, and it  
can be integrated exactly as
\bea
\label{solu_step2}
\phi(\bx,t)&=&\exp\left\{ -i \left[  (t-t_n) \beta|\phi^n(\tbx)|^2+ \lambda\, \bm{\varphi}(\tbx, t)  + P(\bx,t) \right]\right\},\quad \tbx\in {\textbf B}, \quad t_n\leq t\leq t_{n+1}, \\
\label{solu_step2_conv}
\bm{\varphi}(\tbx, t) &=&  \int_{\mathbb{R}^d}  \widetilde{ \mathcal{K}}(\tby,t) \rho(\tbx-\tby,t^n) d\tby,
\eea
where the time-dependent kernel $\widetilde{ \mathcal{K}}(\tbx,t)$ has the form
\be
\widetilde{ \mathcal{K}}(\tbx,t)=\int_{t^n}^{t} \widetilde{\mathcal{U}}(\tbx, \tau)d\tau=
\left\{
\begin{array}{ll}
   (t-t_n)/(2^{d-1}\pi |\tbx|^\mu),      &\quad  {\rm Coloumb},  \\[0.5em]
   -\delta(\tbx) (t-t^n)-3 \widetilde{L}_3(t)(\fl{1}{4\pi|\tbx|}),    &\quad    {\rm 3D\ \ DDI },\\[0.5em]
     -\fl{3}{2} \widetilde{L}_2(t)(\fl{1}{2\pi|\tbx|}),    & \quad   {\rm 2D\ \ DDI }.
\end{array}
\right.
\ee
Here, the differential operators $\widetilde{L}_3(t)=\int_{t^n}^{t}\p_{\bbm(\tau)\bbm(\tau)} d\tau$
and $\widetilde{L}_2(t)=\int_{t^n}^{t}\p_{\bbm_{\perp}(\tau)\bbm_{\perp}(\tau)} d\tau$
can be actually  integrated analytically and have some explicit expressions. One
 refers to {\em section 4.1} in \cite{BMTZ2013}
for more details.  The GauSum solver is then applied to evaluate
the nonlocal nonlinear interaction $\bm{\varphi}(\tbx, t) $ (\ref{solu_step2_conv}). In addition, we have
\be
\label{int_poten}
P(\tbx,t)=\int_{t^n}^{t}\mathcal{W}(\tbx,\tau)d\tau=\int_{t^n}^{t}V(A(\tau)\tbx)d\tau.
\ee
If  $V(\bx)$
is chosen as the harmonic potential (\ref{harm_poten}), then   $P(\tbx,t)$  
can be calculated analytically. For a general potential, a numerical quadrature can be used to approximate the 
integral (\ref{int_poten}). 

To simplify the notations, we  only present the scheme for the 2D case. 
 Let $L$ and $M$ be two even positive integers. We choose $h_\tx=\fl{R_\tx-L_\tx}{L}$ and $h_\ty=\fl{R_\ty-L_\ty}{M}$
 as the spatial mesh sizes in the $\tx$- and $\ty$-directions, respectively.  We define the indices and grid points sets as
\beas
{\mathcal T}_{LM} &=& \left\{(\ell, m)\in\mathbb{N}^2\,|\,0\leq \ell\leq L, \ 0\leq m \leq M \right\}, \\
\widetilde{\mathcal T}_{LM} &=& \left\{(p, q)\in\mathbb{N}^2\,|\,-L/2\leq p\leq L/2-1,
\ -M/2\leq q\leq M/2-1\right\},\\
{\mathcal G}_{\tx\ty} &=& \left\{ (\tx_\ell, \ty_m) =: (L_x + \ell \,h_x, L_y + m \,h_y ),\ (\ell,m) \in  {\mathcal T}_{LM}  \right\}.
\eeas
We introduce the following functions
\[W_{pq}(\tx,\ty)=e^{i \mu_p^\tx(\tx-L_\tx)}\,e^{i \mu_q^\ty(\ty-L_\ty)},
\quad (p, q)\in\widetilde{\mathcal T}_{LM},
\]
with
\[  \mu_p^\tx = \fl{2\pi p }{R_\tx-L_\tx}, \;\quad \mu_q^\ty = \fl{2\pi q}{R_\ty-L_\ty},\quad (p, q)\in \widetilde{\mathcal T}_{LM}.\]
Let $ f_{\ell m}^n$ ($f=\phi$, $\bm{\varphi}$ or $P$) be the approximation of $f(\tx_{\ell}, \ty_m, t_n)$ for
$ (\ell, m) \in {\mathcal T}_{LM}$ and $n\ge0$. We denote by $\bm{\phi}^n$  the solution 
at time $t=t_n$, with components $\{\phi_{\ell m}^n, \ (\ell,m)\in {\mathcal T}_{LM}\}$.
We take the initial data as
$\phi_{\ell m}^0=\phi_0(\tx_{\ell}, \ty_m)$, for $(\ell,m)\in {\mathcal T}_{LM}$. 
A second-order time-splitting Fourier pseudo-spectral (TSFP) method to solve  \eqref{FSE1_Rot}-\eqref{NonLoc_Rot} is
given by
\bea
\phi_{\ell m}^{(1)}
\label{tssp1}
&=&\sum_{p=-L/2}^{L/2-1}\,\sum_{q=-M/2}^{M/2-1}
	e^{-\fl{i\Dt t}{4} \left[ (\mu_p^\tx)^2+ (\mu_q^\ty)^2+m^2\right]^{s}}
	\widehat{(\bm{\phi}^n)}_{pq}\; W_{pq}(\tx_{\ell},\ty_m), \\[0.5em]
\phi_{\ell m}^{(2)}
\label{tssp2}
&=&\phi_{\ell m}^{(1)}\exp\left\{ -i  \left[ \Dt t \beta|\phi^n_{\ell m}|^2+\lambda \bm{\varphi}^{n+1}_{\ell m} + P^{n+1}_{\ell m} \right] \right\} ,  \\[0.5em]
\phi_{\ell m}^{n+1}
\label{tssp3}
&=&\sum_{p=-L/2}^{L/2-1}\,\sum_{q=-M/2}^{M/2-1} e^{-\fl{i\Dt t}{4} \left[ (\mu_p^\tx)^2+ (\mu_q^\ty)^2+m^2\right]^{s}} 
\widehat{(\bm{\phi}^{(2)})}_{pq}\; W_{pq}(\tx_{\ell},\ty_m).
\eea
Here, $\widehat{(\bm{\phi}^n)}_{pq}$ and $ \widehat{(\bm{\phi}^{(2)})}_{pq}$ are
the discrete Fourier series coefficients of
the vectors $\bm{\phi}^n$ and $\bm{\phi}^{(2)}$,
respectively.  This method is referred  to as {\em TS2-GauSum}. 
The {\em TS2-GauSum} method (\ref{tssp1})-(\ref{tssp3}) is explicit, efficient, simple to implement, unconditionally stable and can 
be easily extended to high-order time-splitting schemes.


\subsection{Numerical results}

In this subsection, we present some numerical results for the dynamics of the FNLSE/FGPE solved by {\em TS2-GauSum}.  
To this end, unless stated,  we let $m=0$, $\Og=0$,  $d=2$  and choose the computational domain as 
${\textbf B}=[-16, 16]\times[-16, 16]$. The mesh sizes in space and time are chosen  as $h_x=h_y=\fl{1}{8}$ and $\Delta t=10^{-3}$, 
respectively. The trapping potential  $V(\bx)$ is chosen as (\ref{harm_poten}) with $\gm_x=\gm_y=1.$  The nonlocal 
interaction is of Coulomb-type with $\mu=1.$  The initial data is set to
\be
\label{ini_dyn}
\psi_0(\bx)=\phi^{s}_g(\bx-\bx_0)\, e^{i\,v_0 (0.8x+0.5y)},
\ee
 where $\phi^{s}_g$ is the ground state of the FNLSE with the fractional order $s$. 
Starting from the  ground state $\phi_{g}^{s}(\bx)$, we shift it by $\bx_{0}\in \mathbb R^{2}$ and/or imprint an
 initial momentum as shown above.

\begin{exmp}\label{dy:exmp:1}
{\bf Dynamics of the FNLSE  ($\bx_{0}=(0,0)^{T}$).}
In this example, let $\beta=0$, $\lambda=-1$, $v_0=1$ and $\bx_0=(0, 0)^T$ in (\ref{ini_dyn}).  
We study two cases in (\ref{ini_dyn}): {\bf Case I:}  $s=1$,   {\bf Case II:}  $s\neq 1$.
\end{exmp}

\begin{figure}[h!]
\centerline{
(a)
\psfig{figure=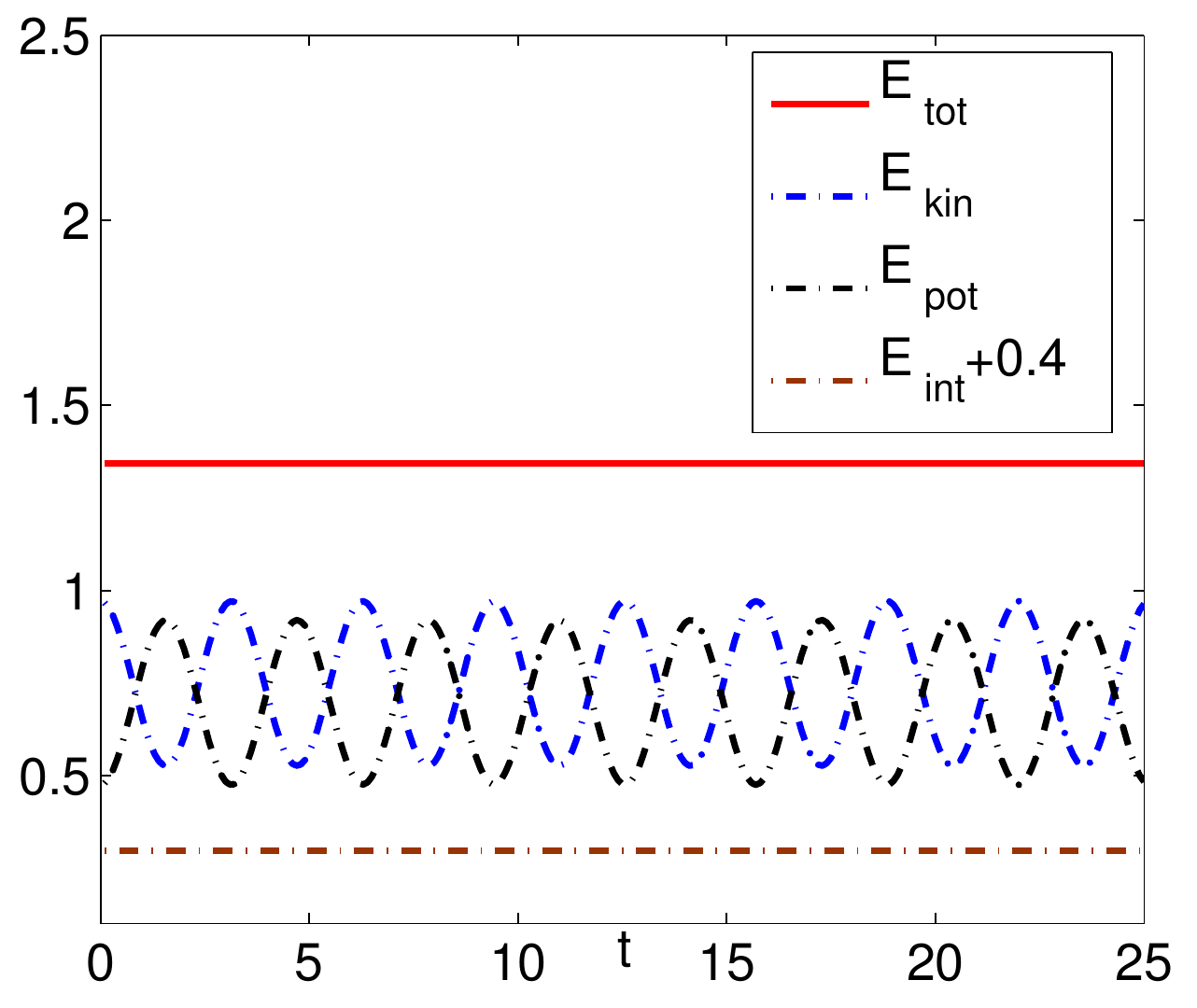,height=4.4cm,width=5.3cm,angle=0}\quad
\psfig{figure=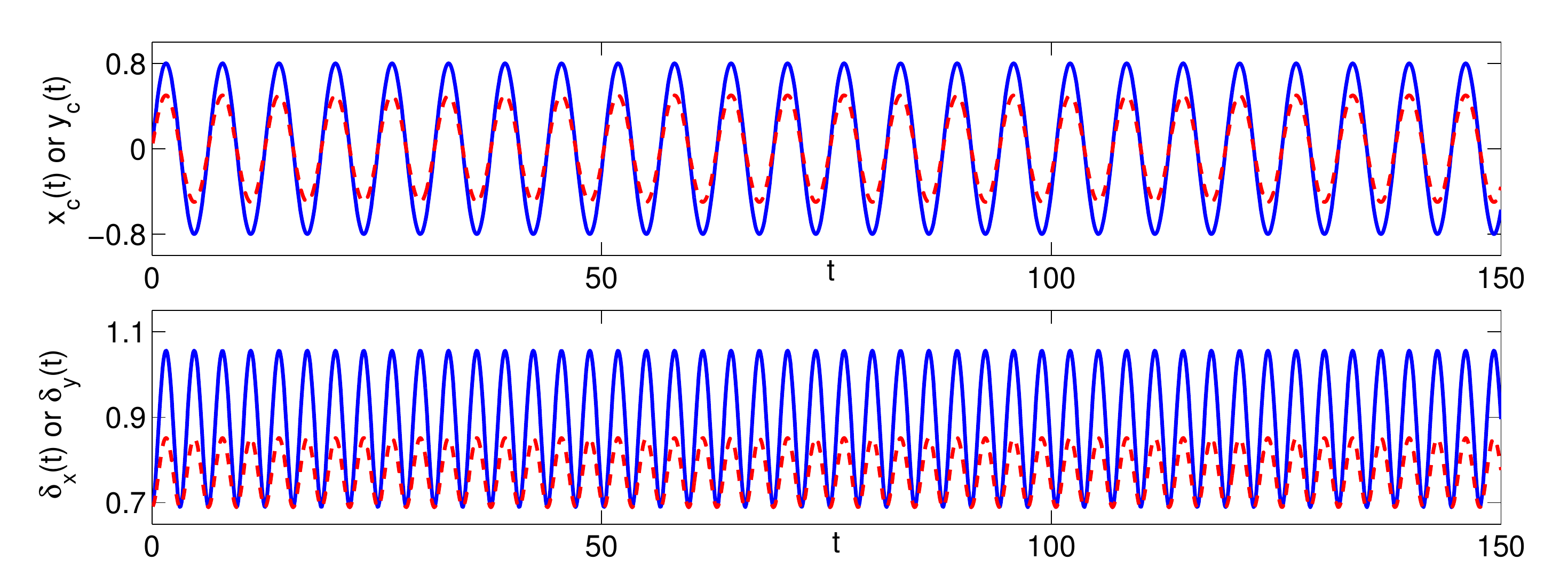,height=4.6cm,width=11cm,angle=0}
}
\centerline{
(b)
\psfig{figure=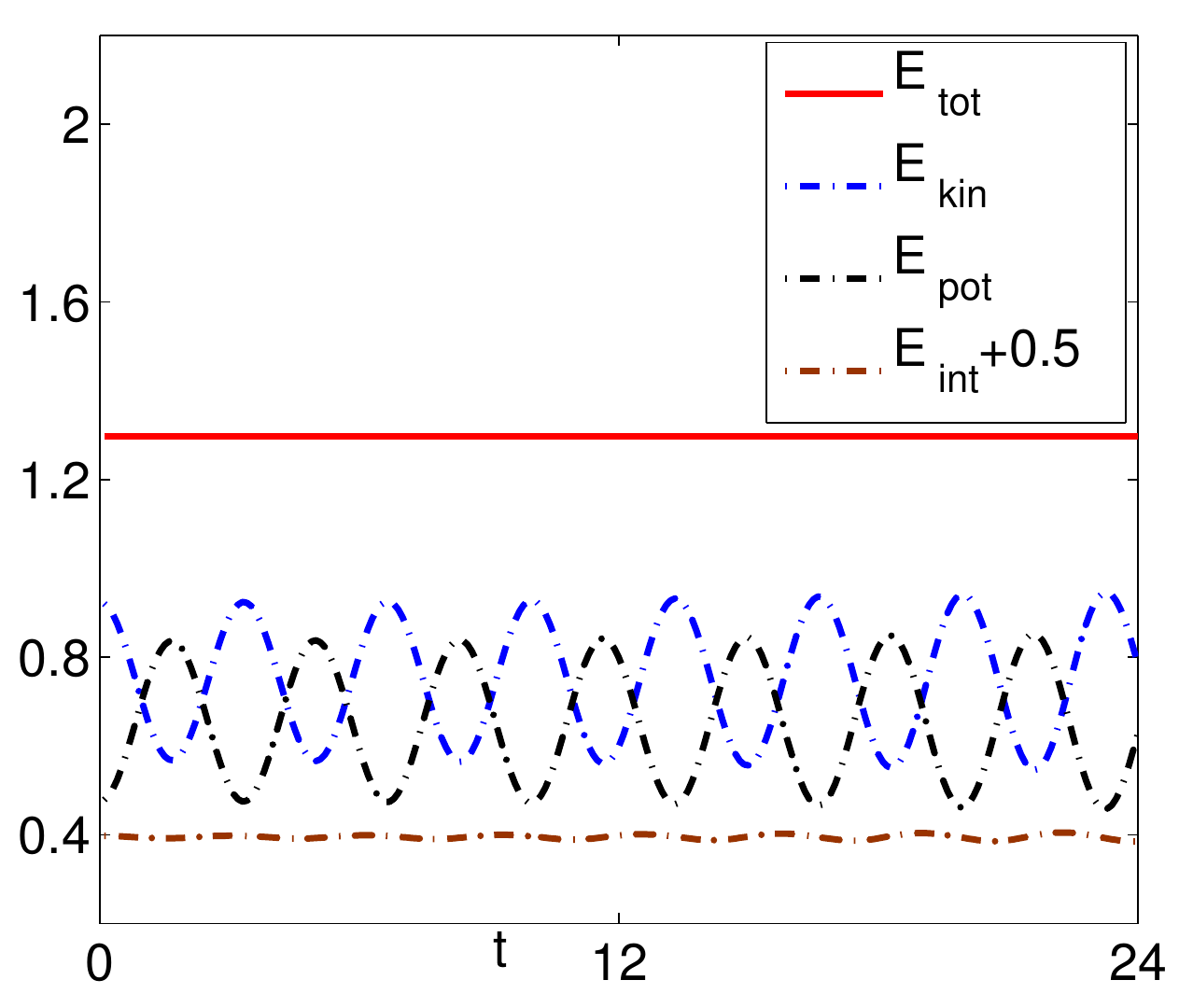,height=4.4cm,width=5.3cm,angle=0}\quad
\psfig{figure=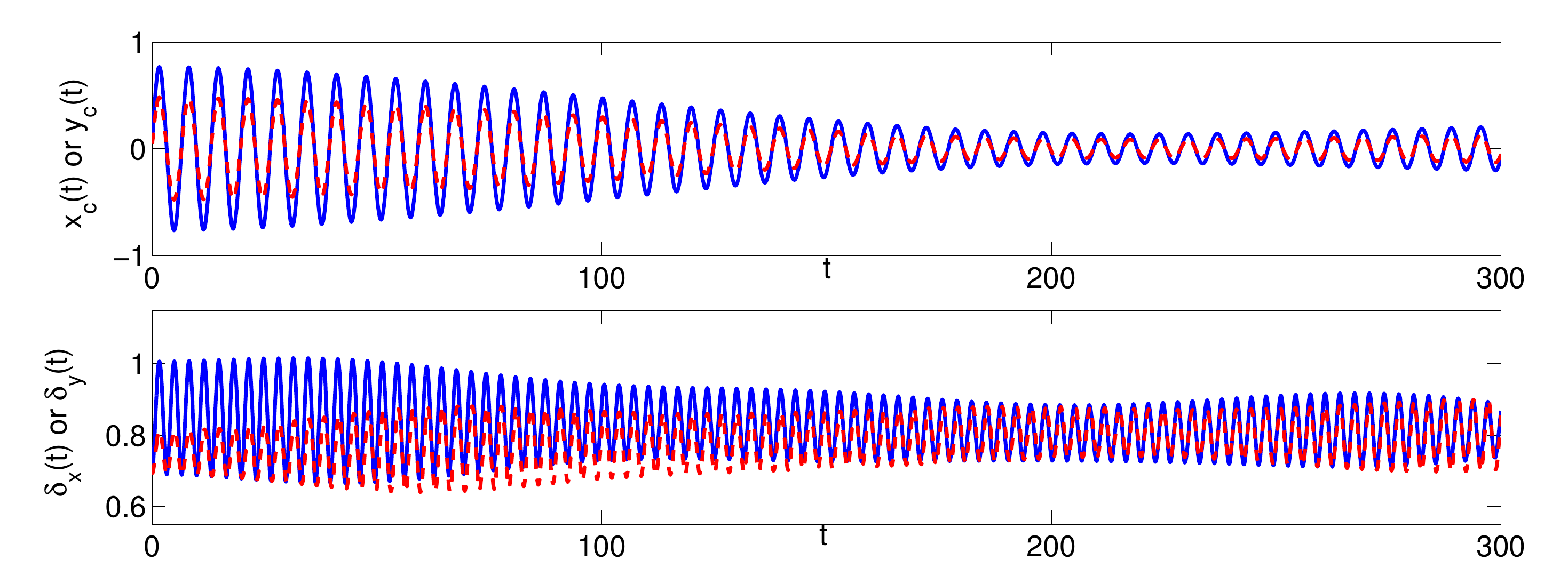,height=4.6cm,width=11cm,angle=0}
}
\centerline{
(c)
\psfig{figure=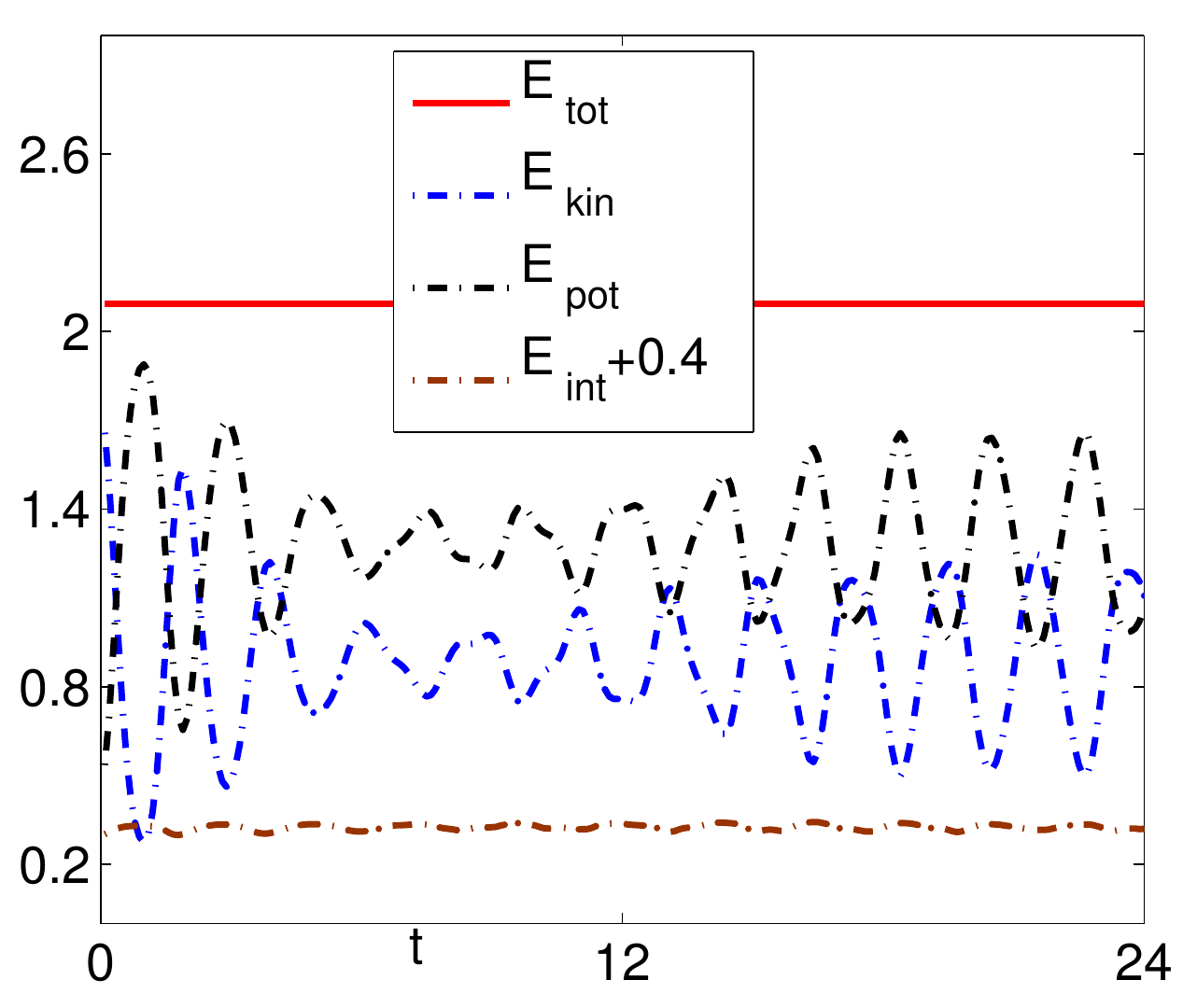,height=4.4cm,width=5.3cm,angle=0}\quad
\psfig{figure=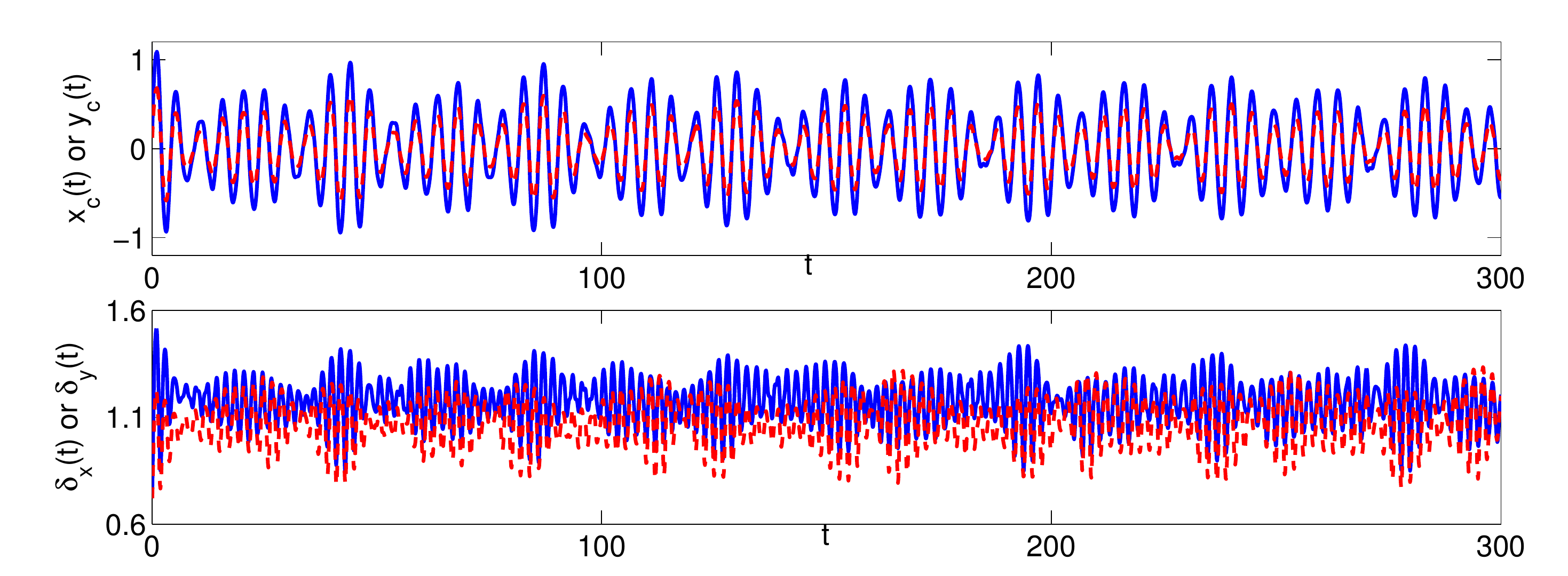,height=4.6cm,width=11cm,angle=0}
}
\caption{Time evolution of the energies, centre of mass and condensate 
width in example \ref{dy:exmp:1} for case I for:   (a)  $s=1$ (standard case),  (b) $s=0.95$ (subdispersion case) and (c)
 $s=1.5$ (superdispersion case).}
\label{fig:ex6}
\end{figure}

\begin{figure}[h!]
\centerline{
(a)
\psfig{figure=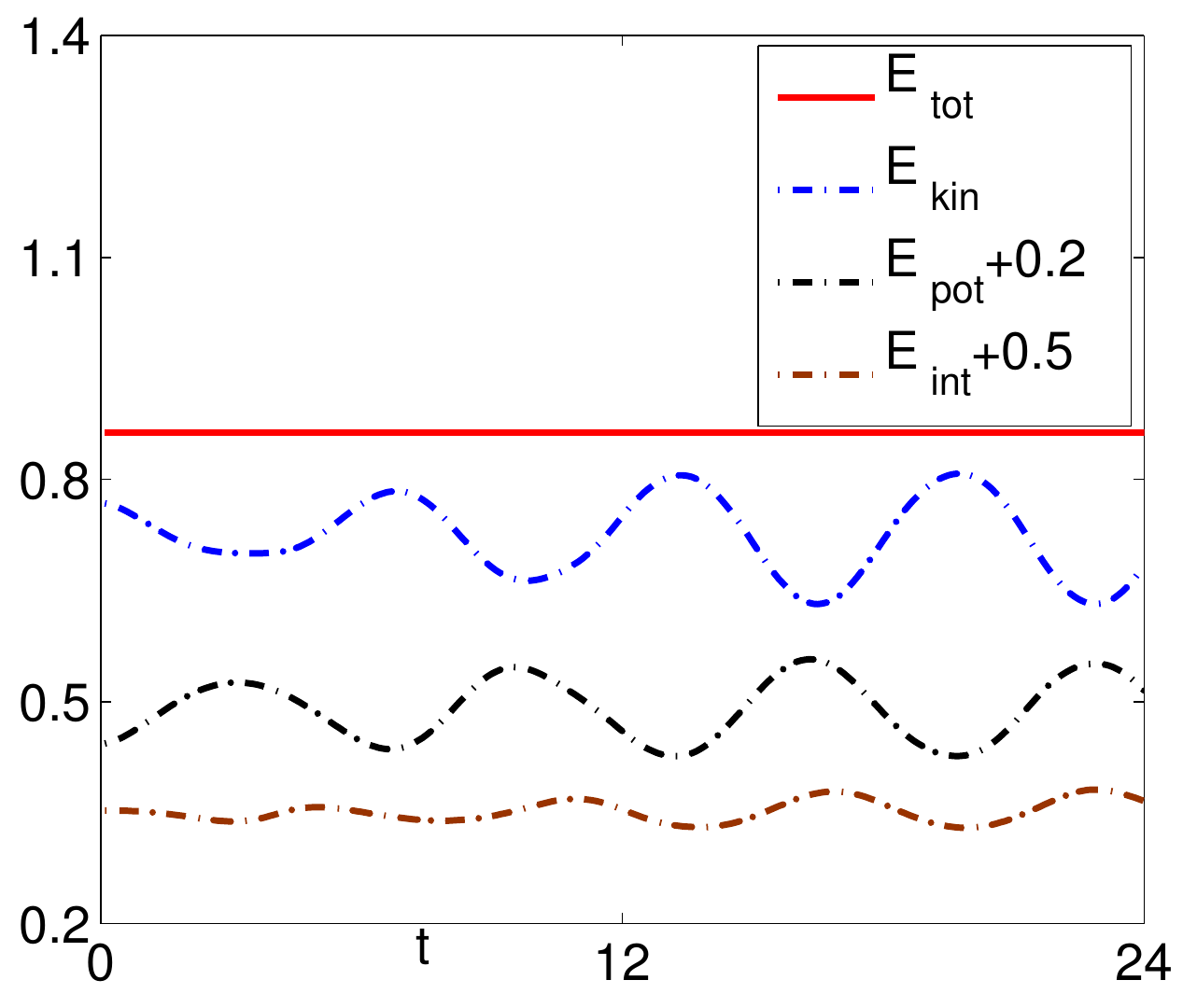,height=4.4cm,width=5.3cm,angle=0}\quad
\psfig{figure=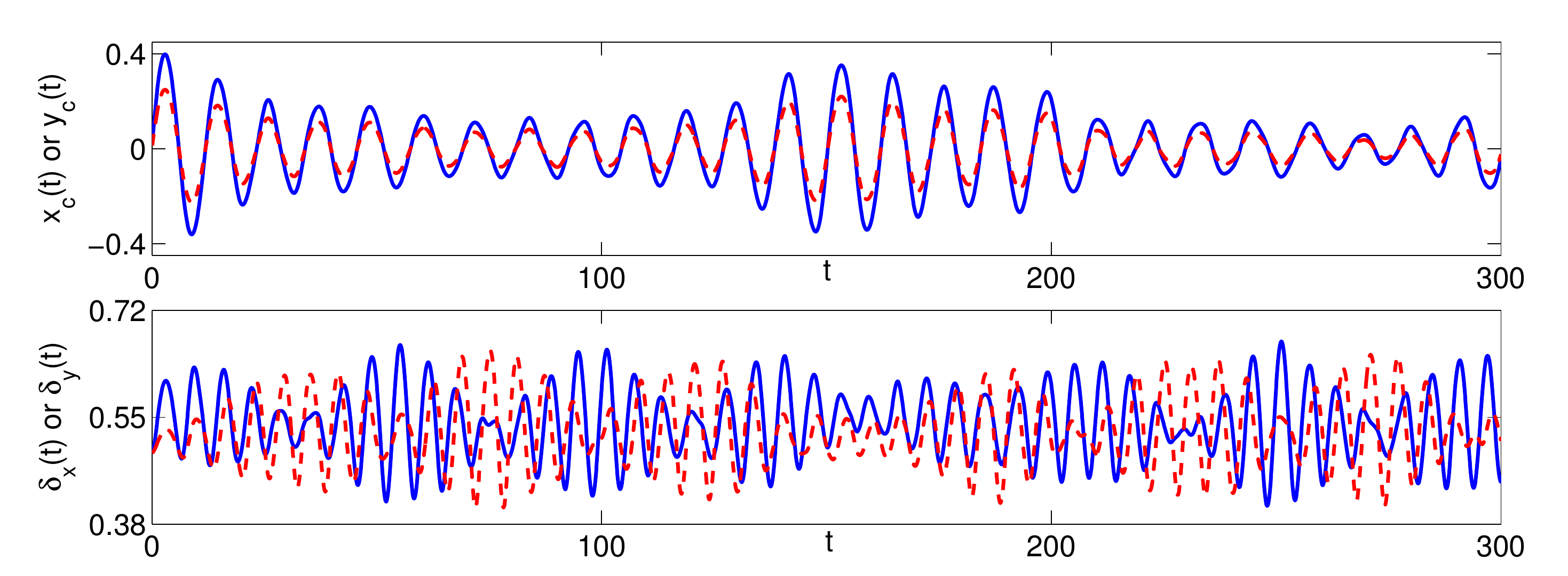,height=4.6cm,width=11cm,angle=0}
}
\centerline{
(b)
\psfig{figure=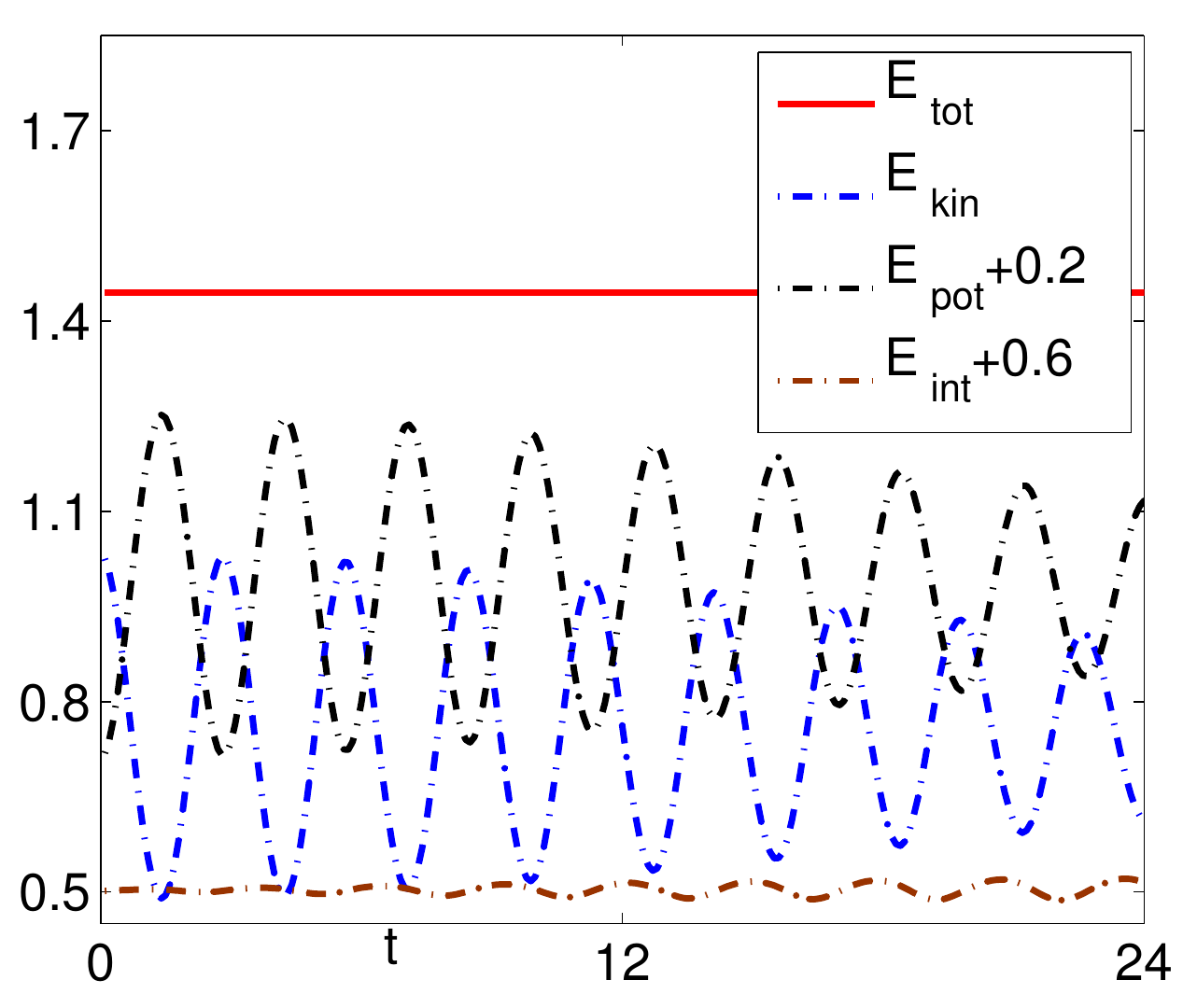,
height=4.4cm,width=5.3cm,angle=0}\quad
\psfig{figure=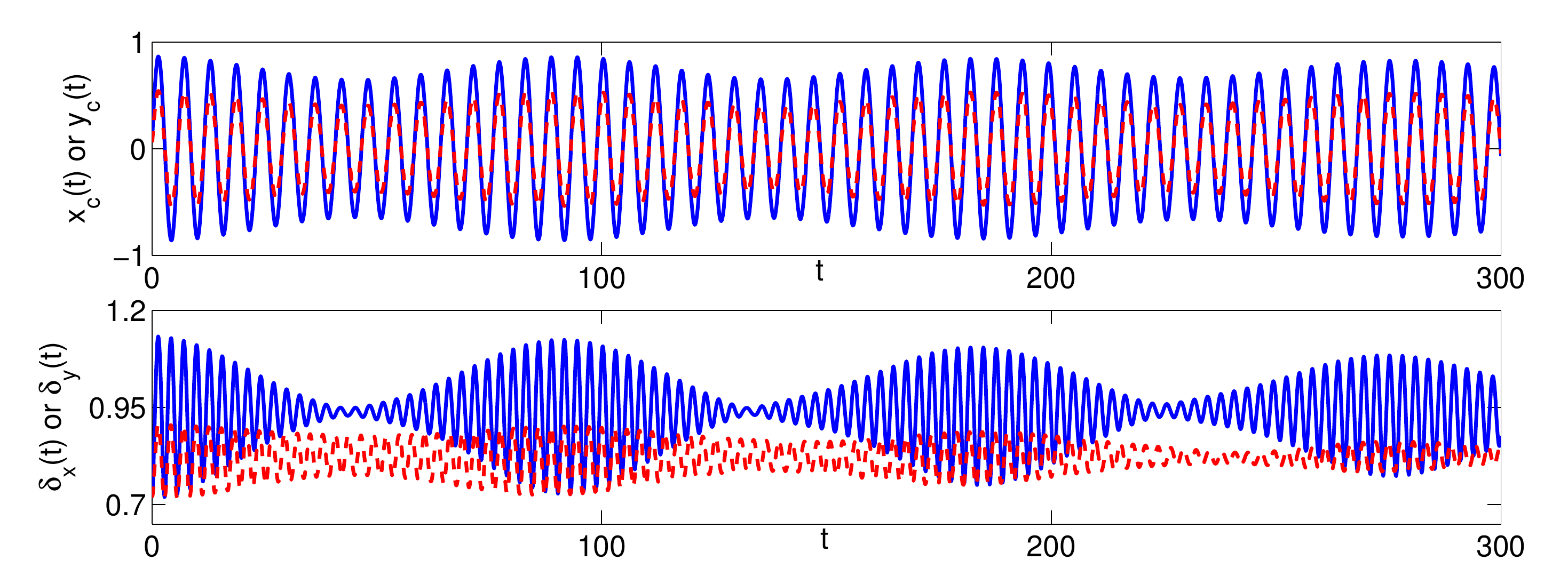,
height=4.6cm,width=11cm,angle=0}
}
\caption{Time evolution of the energies, centre of mass and condensate 
width in example \ref{dy:exmp:1} for case II for:  (a)
$s=0.5$ (subdispersion case), (b) $s=1.1$ (superdispersion case). }
\label{fig:ex6_con}
\end{figure}

Figure \ref{fig:ex6} and  \ref{fig:ex6_con} show the dynamics of mass, energy, centre of mass, 
condensate widths of the FNLSE with different fractional orders $s$. 
We can observe that
(i) The mass and total energy are well conserved.
(ii) The fractional order significantly affects the dynamics of the FNLSE. 
As we know, for the classical NLSE ($s =1$), the density profile retains its initial shape,  meanwhile swings periodically in the harmonic trap
(cf. Fig. \ref{fig:ex6} (a)). However, for the fractional case ($s\neq 1$), the density profile is quite different from the initial profile. 
For the subdispersion case, $s<1$, the decoherence emerges, i.e. the loss of solitary profile,
 and it becomes stronger when $|s-1|$ is larger. 
For superdispersion, i.e. $s>1$, there is much less decoherence observed. The density profile would exhibit damped oscillations around what
appears to be a rescaled ground state, which behaves similarly as the breather solutions of the classical NLSE. 
(iii) For both cases, the decoherence is weak and turbulence (the high frequencies)  does not emerge, letting
 alone the chaotic dynamics. 
The  turbulence and/or chaotic dynamics might emerge if the initially imprinted momentum is
large enough.

%

\medskip
 
\begin{exmp}\label{dy:exmp:2}
{\bf Dynamics of the FNLSE with position shifts in initial data.}
With fixed $s=0.75$ (subdispersion) and $v_0=0$ in (\ref{ini_dyn}), we study the following four cases: 
\begin{itemize}
\item {\bf Case I}.  Linear fractional Schr\"{o}dinger equation.  Let  $\beta=\lambda=0$, $\bx_0=(1,1)^{T}$.
\item {\bf Case II}. Linear fractional Schr\"{o}dinger equation.   Let   $\beta=\lambda=0$, $\bx_0=(3,3)^{T}$.
\item  {\bf Case III}. FNLSE with purely short-range interaction.  Let $\beta=50$, $\lambda=0$,  $\bx_0=(3,3)^{T}$.
\item  {\bf Case IV}. FNLSE with purely long-range interaction.  Let $\beta=0$,  $\lambda=10$, $\bx_0=(3,3)^{T}$.
\end{itemize}
\end{exmp}

\begin{figure}[h!]
\centerline{(a)
\psfig{figure=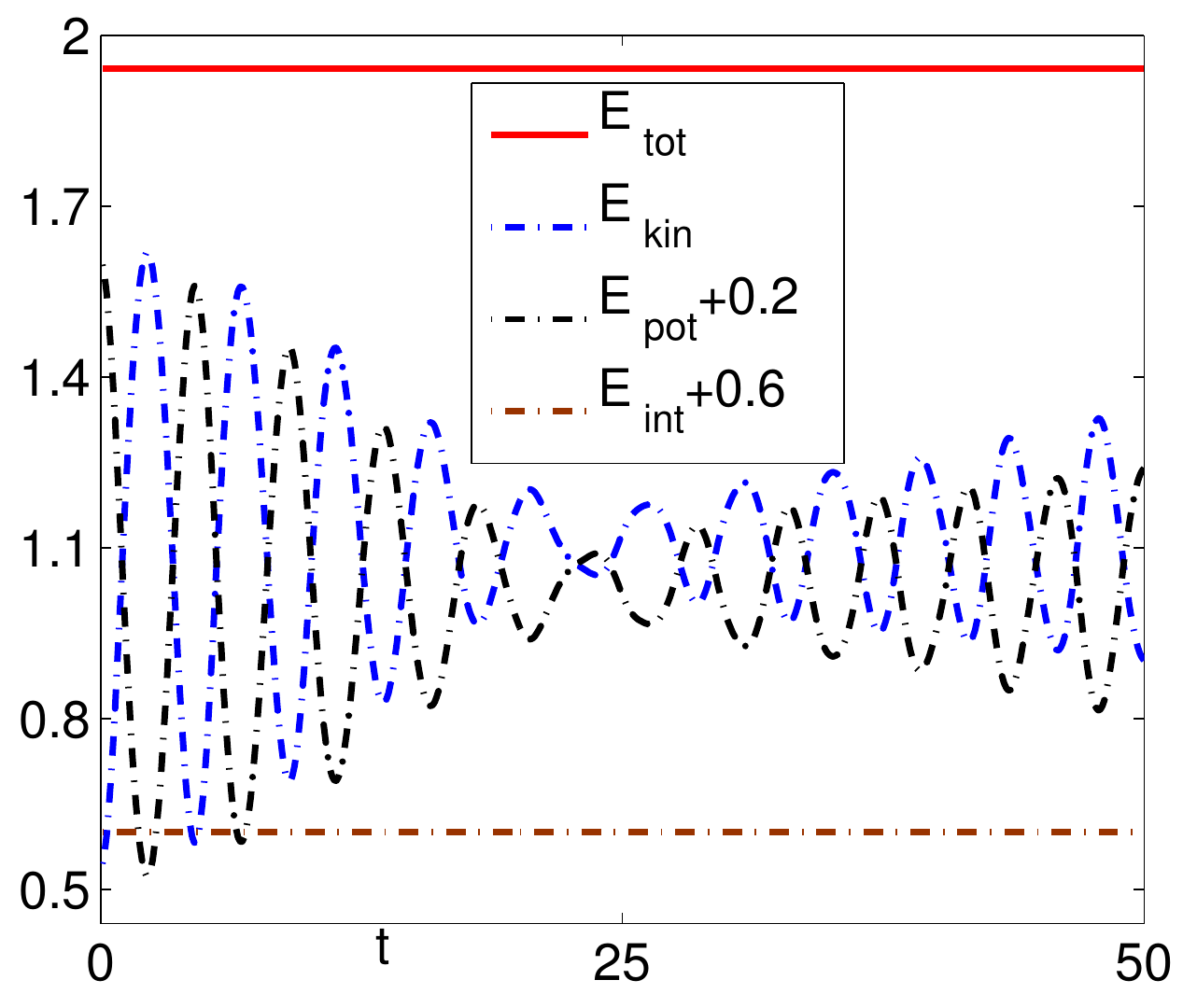,height=4.2cm,width=5.3cm,angle=0}\quad
\psfig{figure=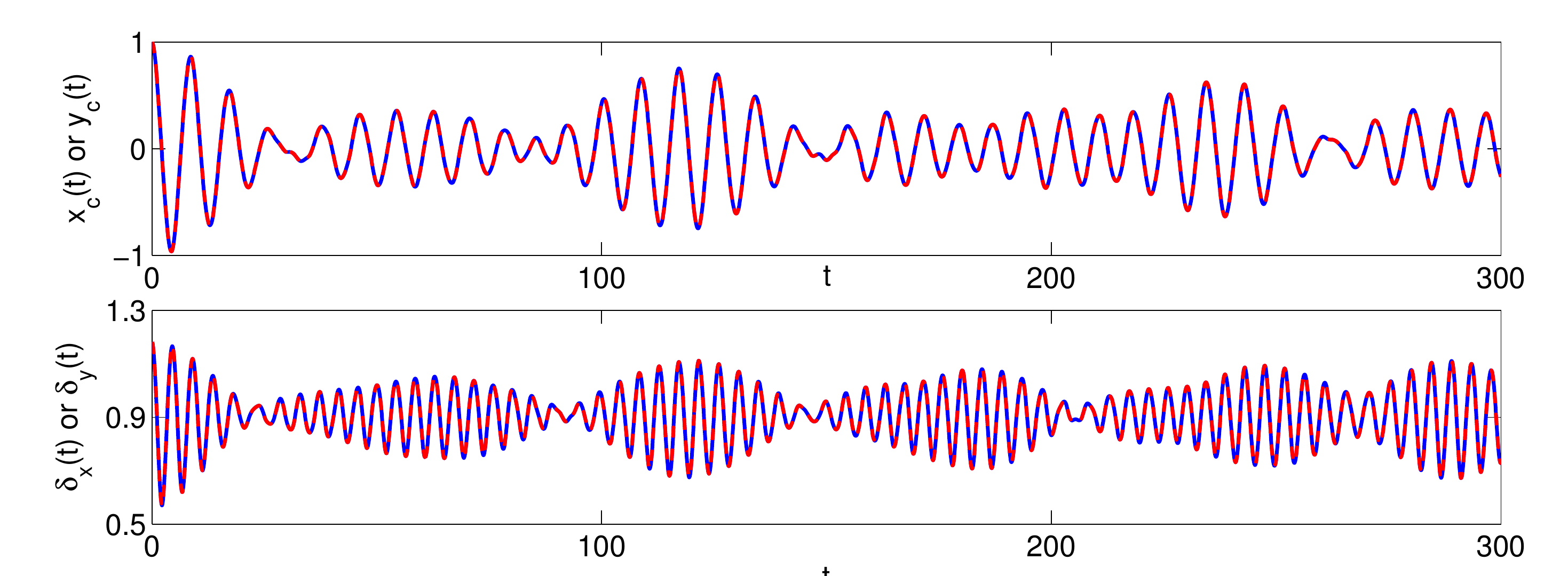,height=4.4cm,width=11cm,angle=0}
}
\centerline{(b)
\psfig{figure=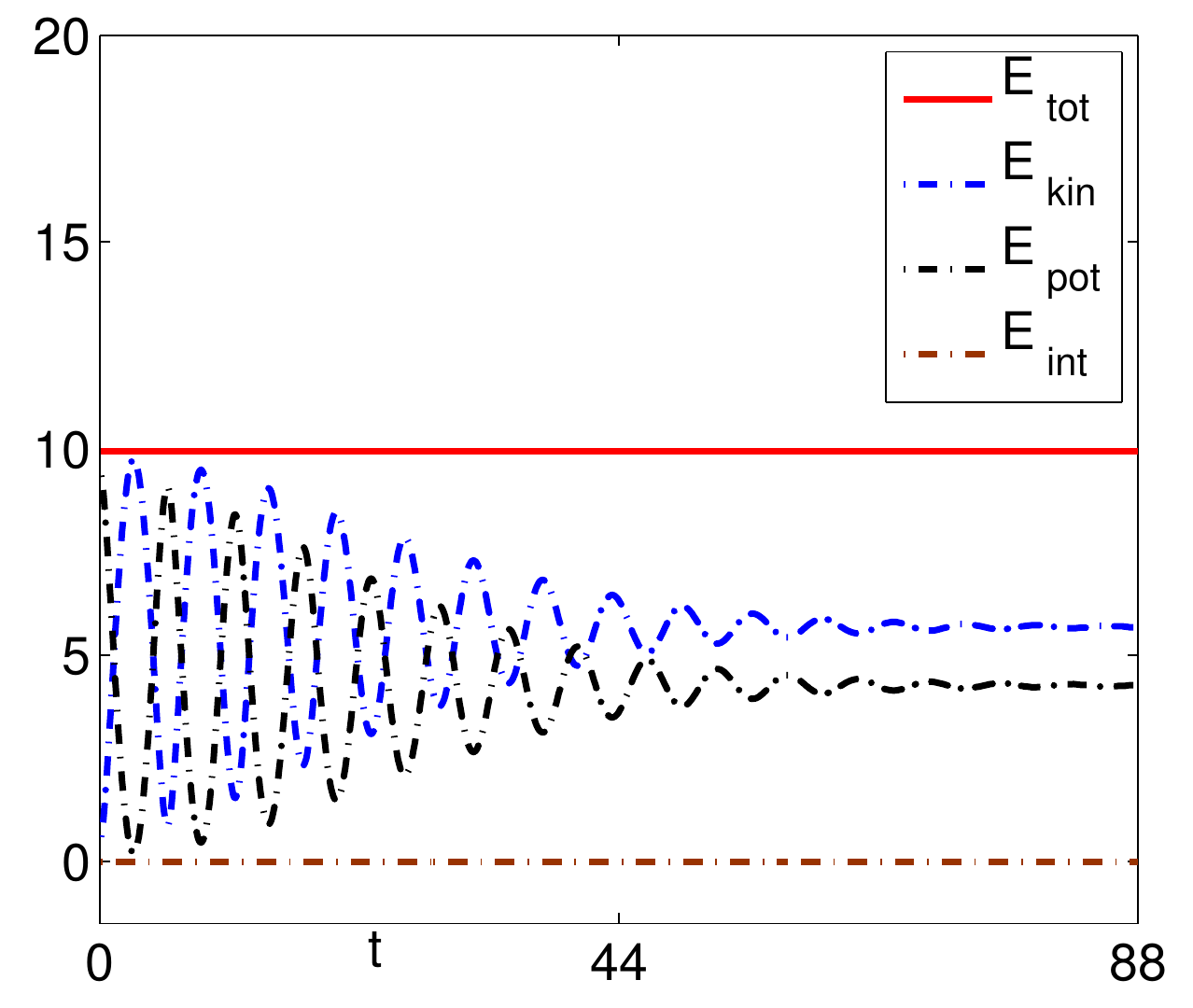,height=4.2cm,width=5.3cm,angle=0}\quad
\psfig{figure=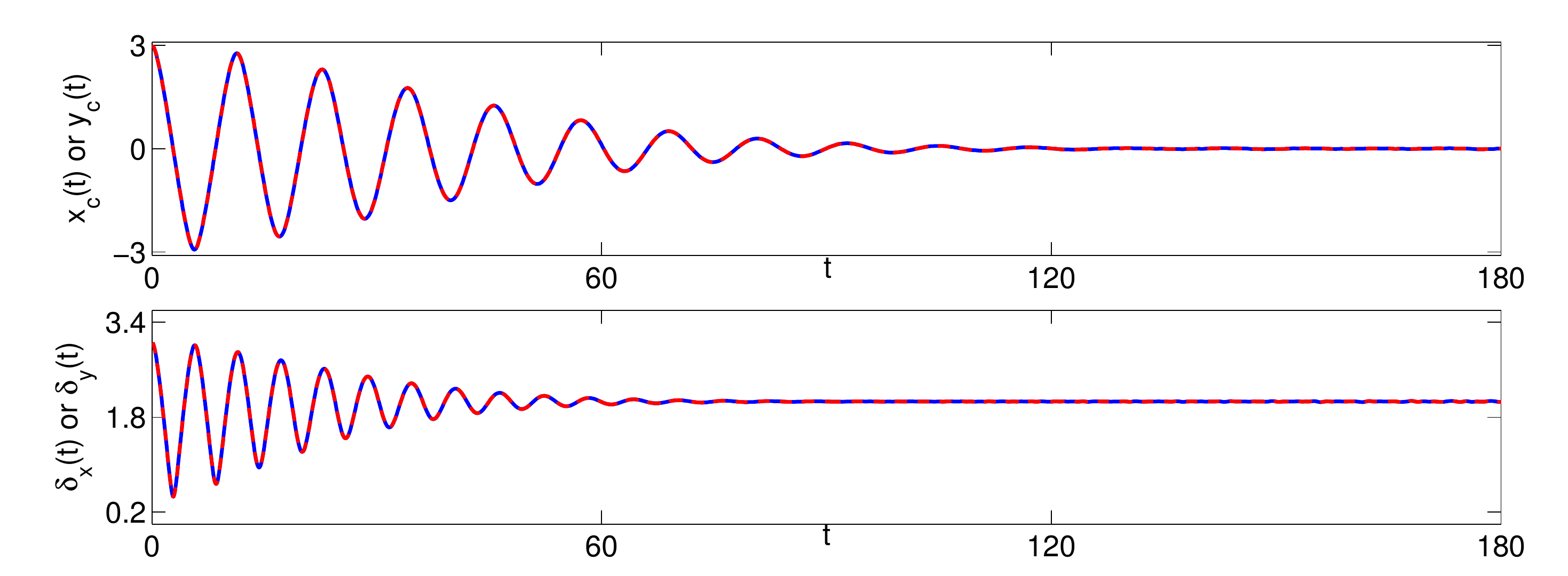,height=4.4cm,width=11cm,angle=0}}
\centerline{(c)
\psfig{figure=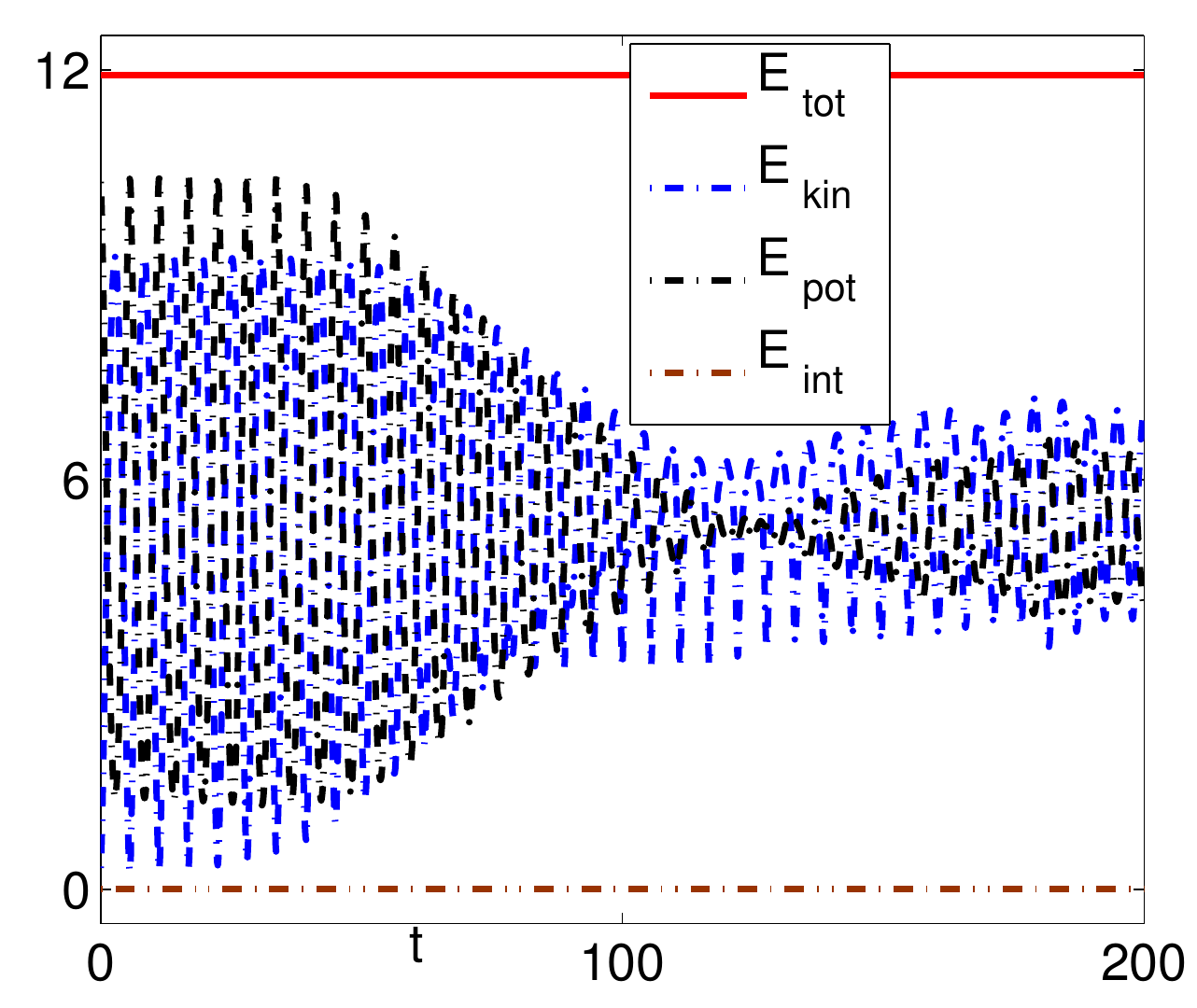,height=4.2cm,width=5.3cm,angle=0}\quad
\psfig{figure=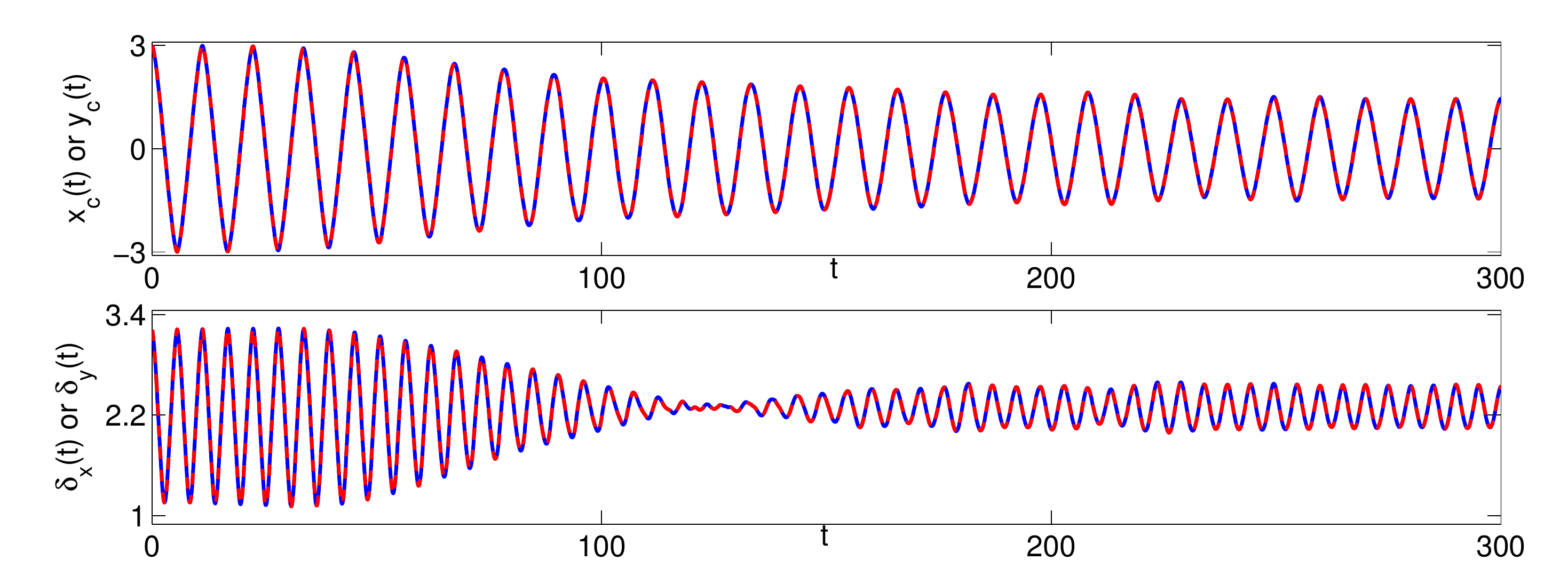,height=4.4cm,width=11cm,angle=0}}
\centerline{(d)
\psfig{figure=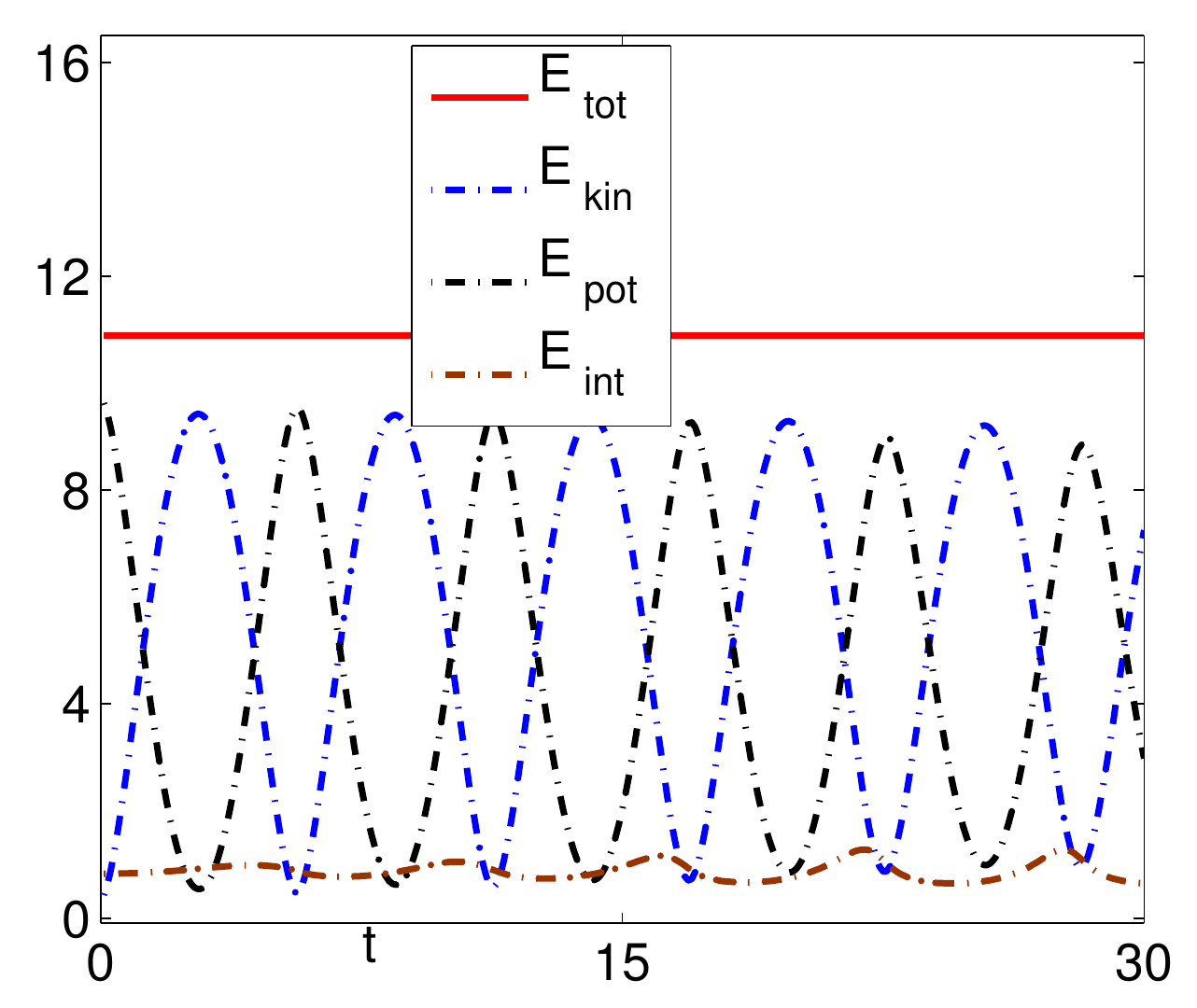,height=4.2cm,width=5.3cm,angle=0}\quad
\psfig{figure=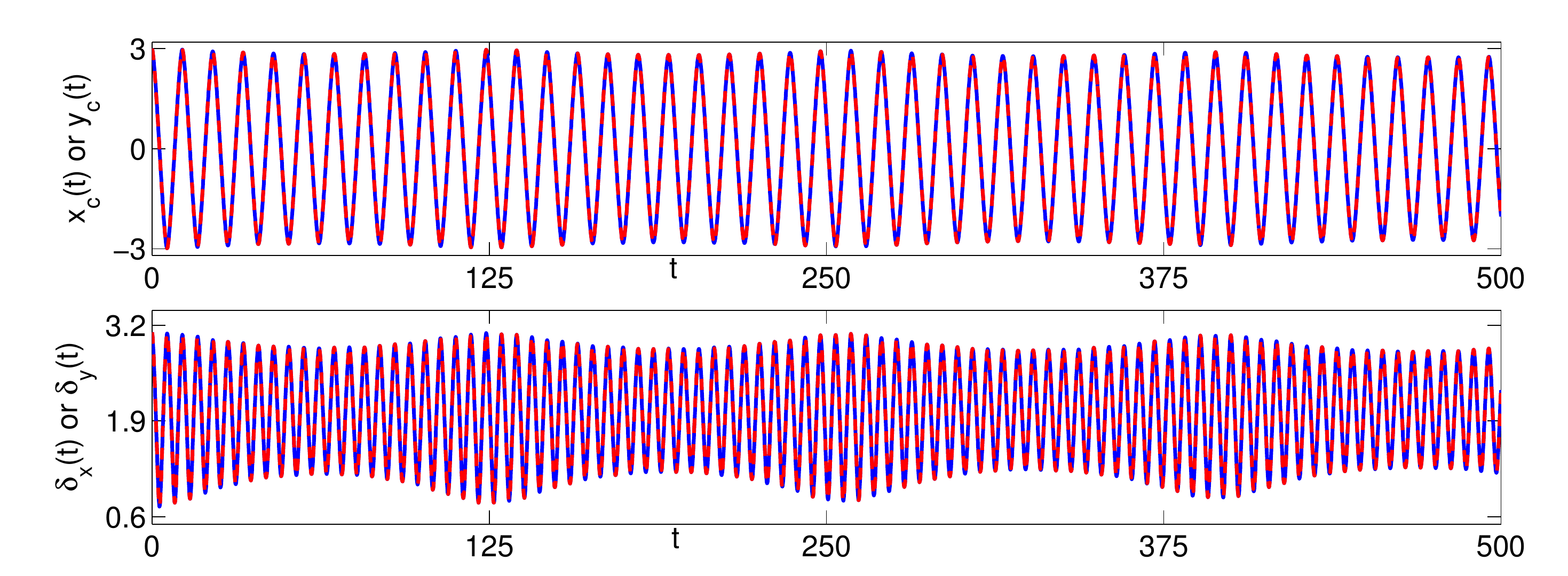,height=4.4cm,width=11cm,angle=0}
}
\caption{Time evolution of the energies, centre of mass and condensate 
width in example \ref{dy:exmp:2}  for cases I to IV (from top to bottom). Here, we consider a subdispersion case for $s=0.75$.
}
\label{fig:ex7}
\end{figure}

\begin{figure}[h!]
\centerline{(a)
\psfig{figure=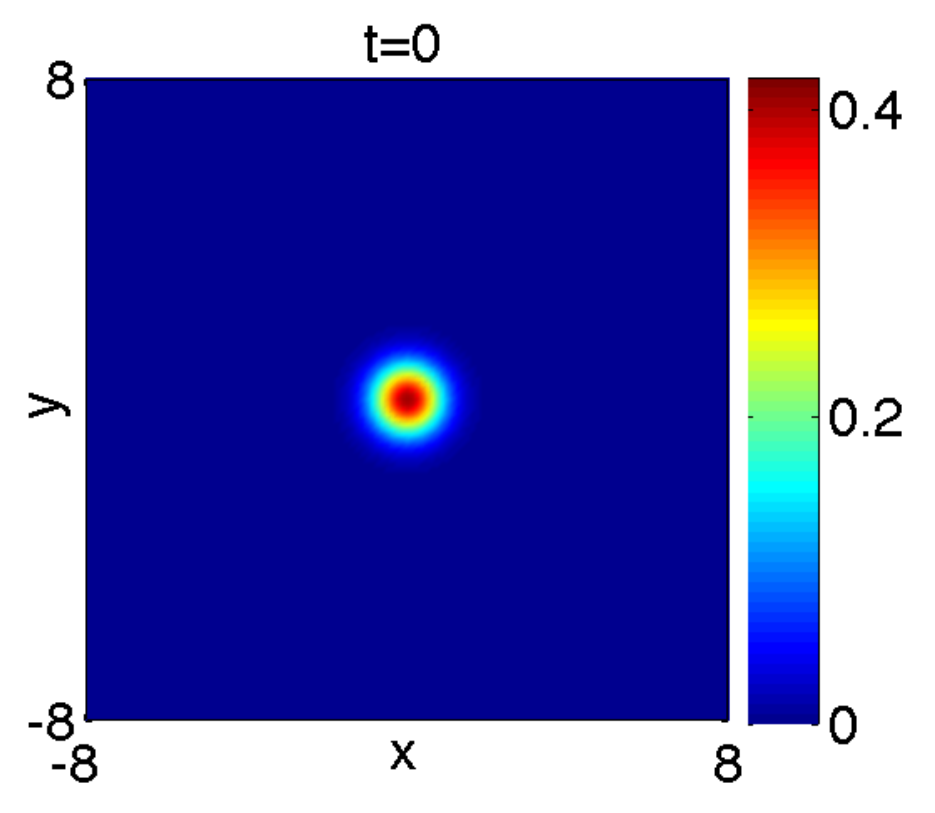,height=3.cm,width=3.3cm,angle=0}
\psfig{figure=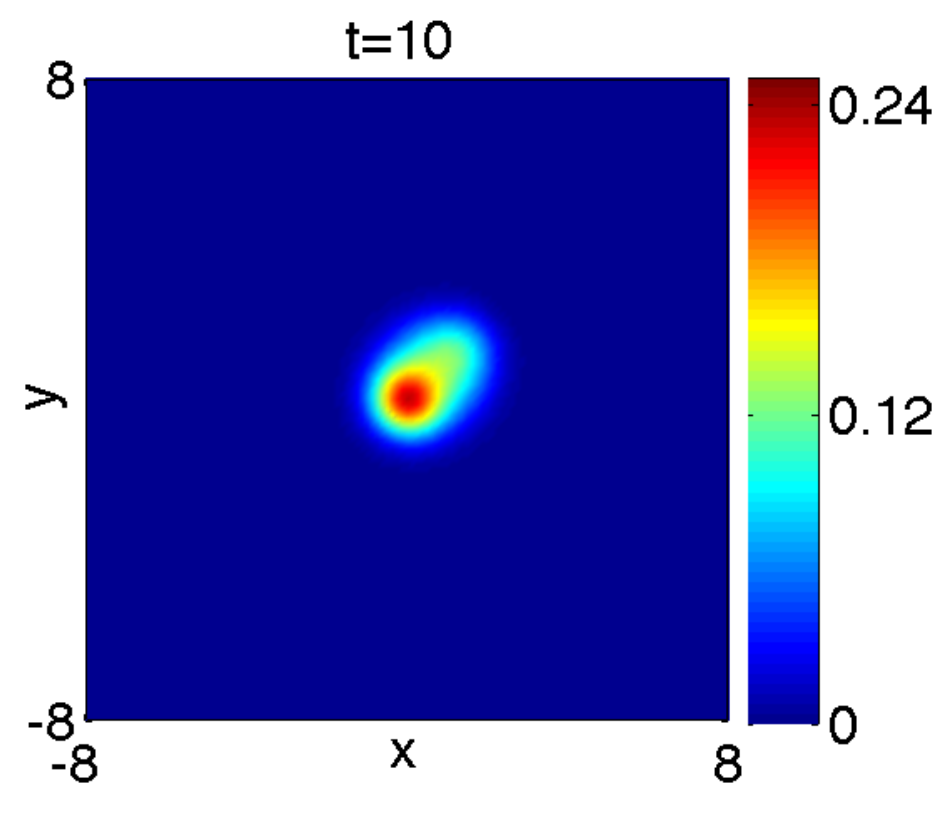,height=3cm,width=3.3cm,angle=0}
\psfig{figure=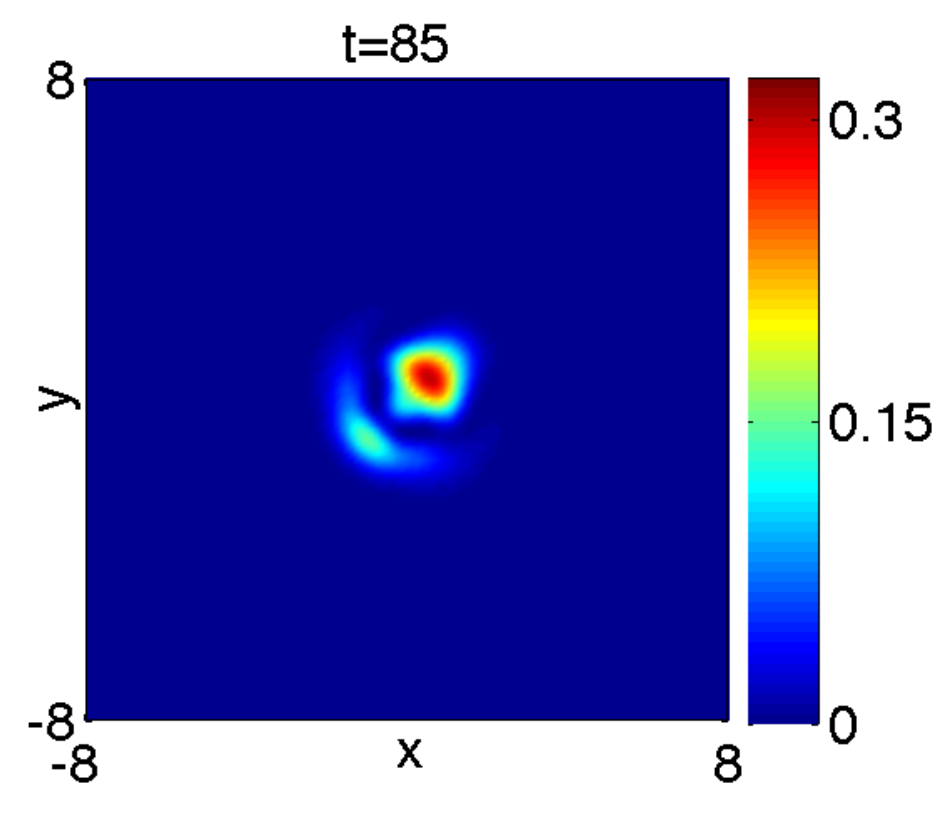,height=3.cm,width=3.3cm,angle=0}
\psfig{figure=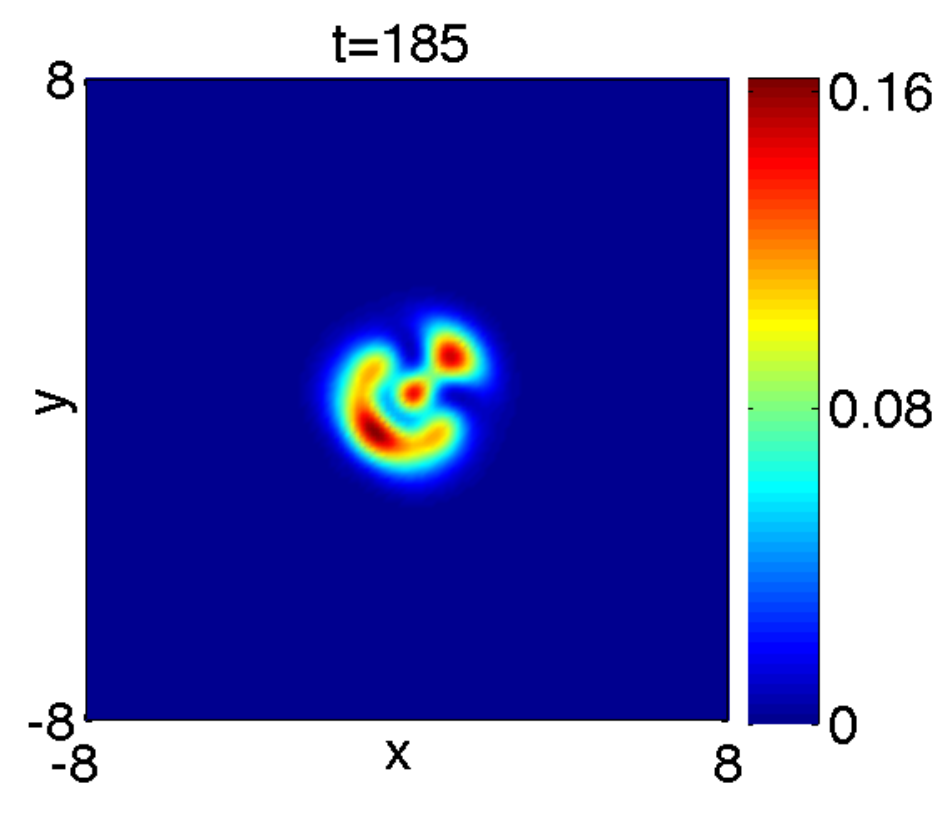,height=3.cm,width=3.3cm,angle=0}
\psfig{figure=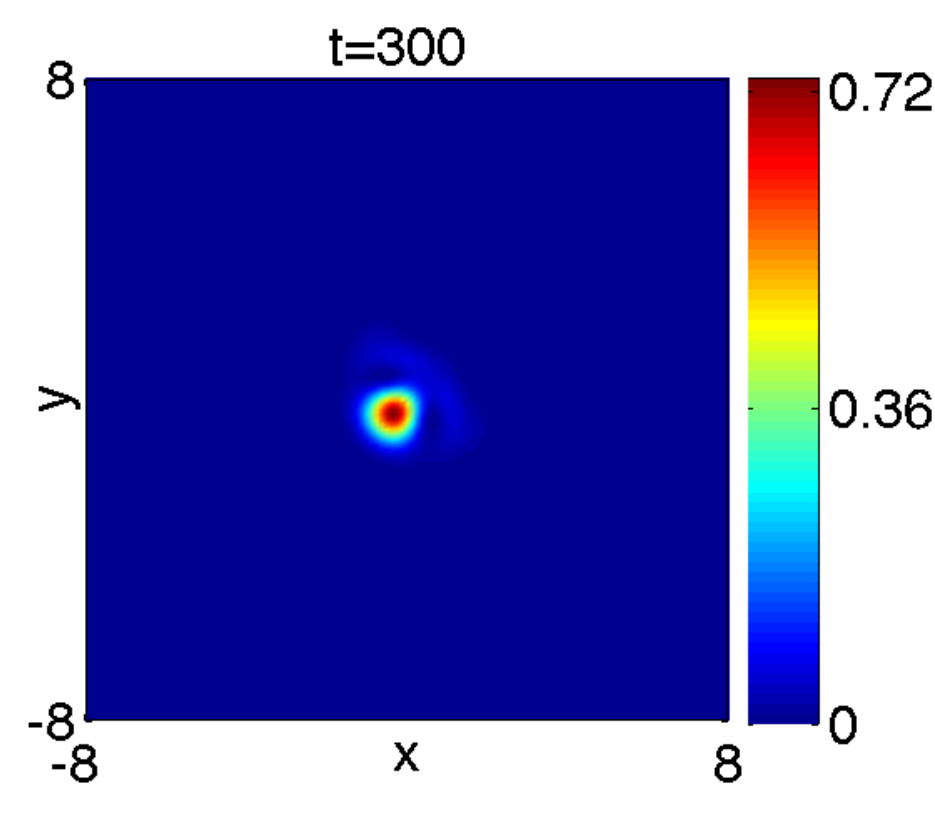,height=3.cm,width=3.3cm,angle=0}
}
\centerline{(b)
\psfig{figure=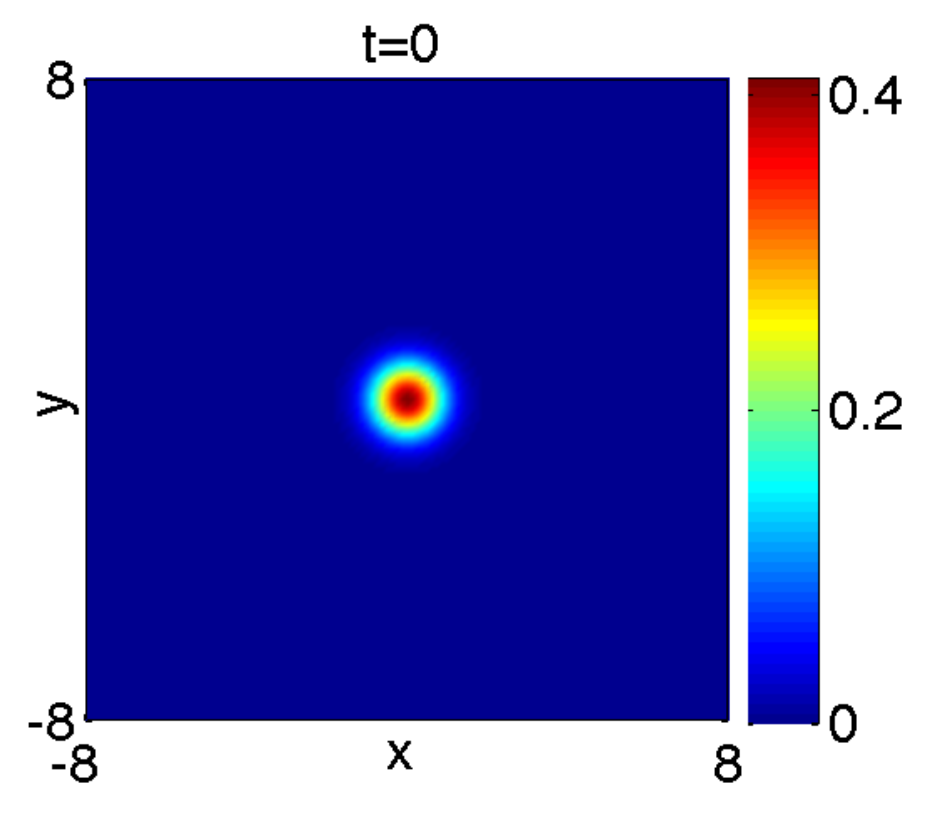,height=3.cm,width=3.3cm,angle=0}
\psfig{figure=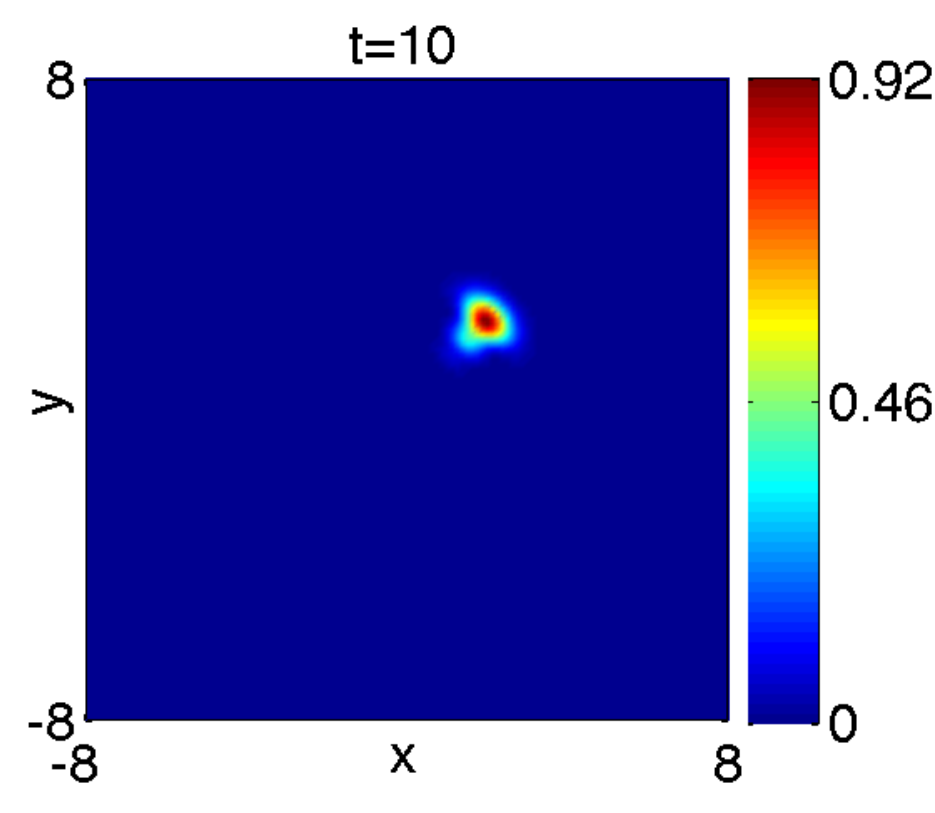,height=3cm,width=3.3cm,angle=0}
\psfig{figure=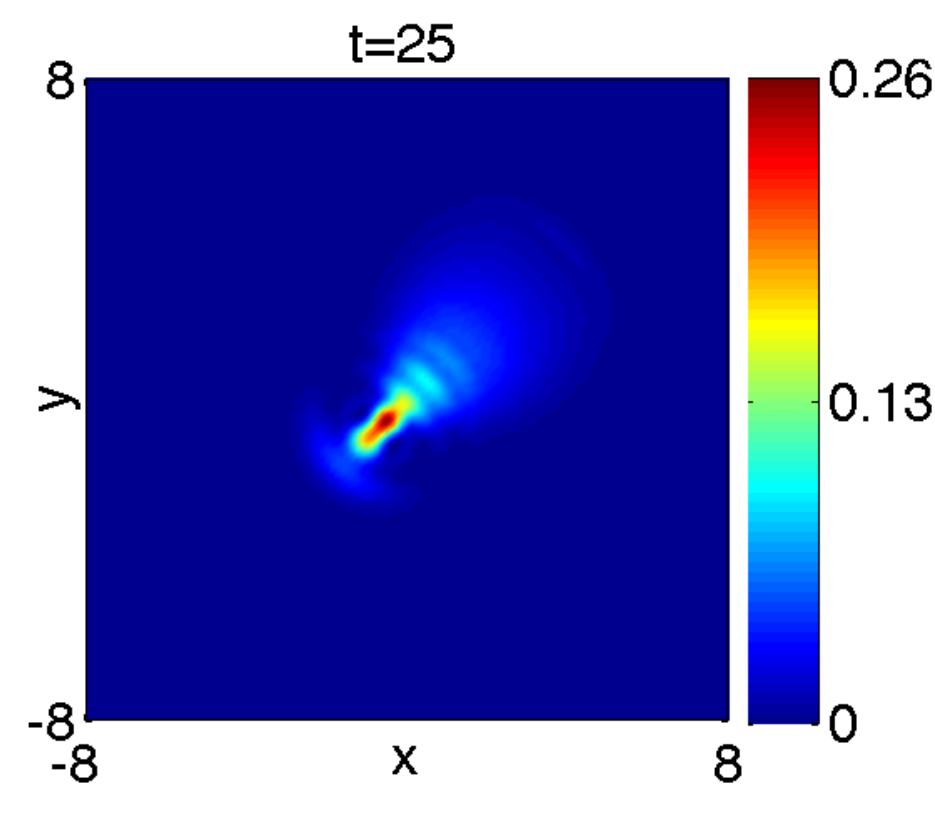,height=3.cm,width=3.3cm,angle=0}
\psfig{figure=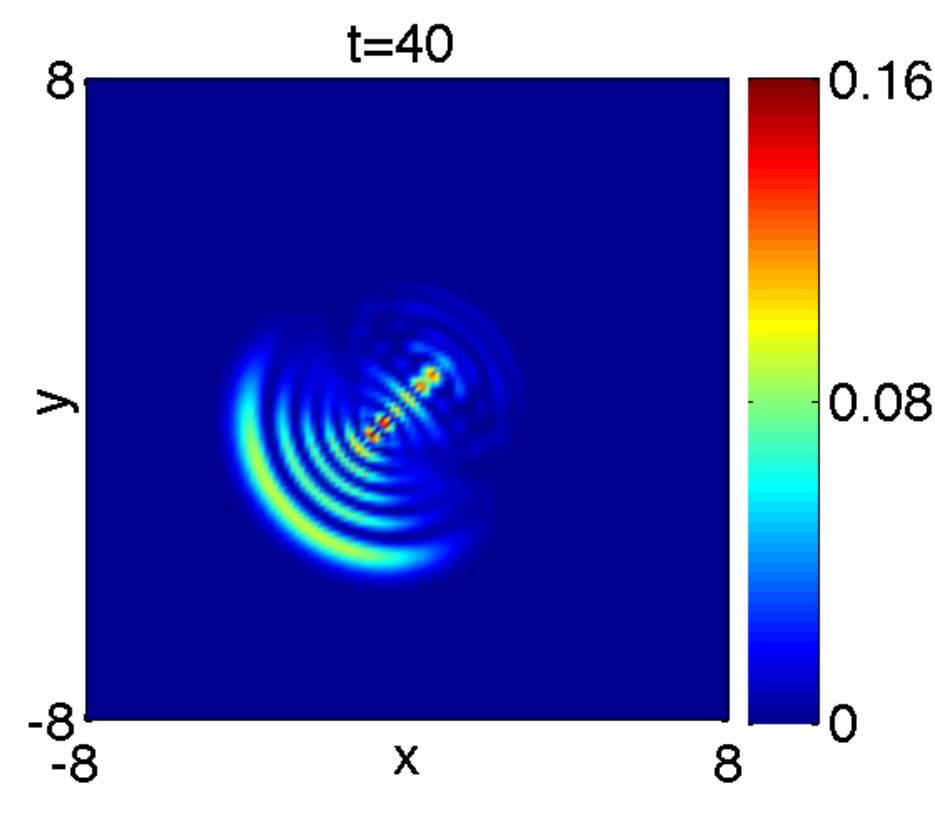,height=3.cm,width=3.3cm,angle=0}
\psfig{figure=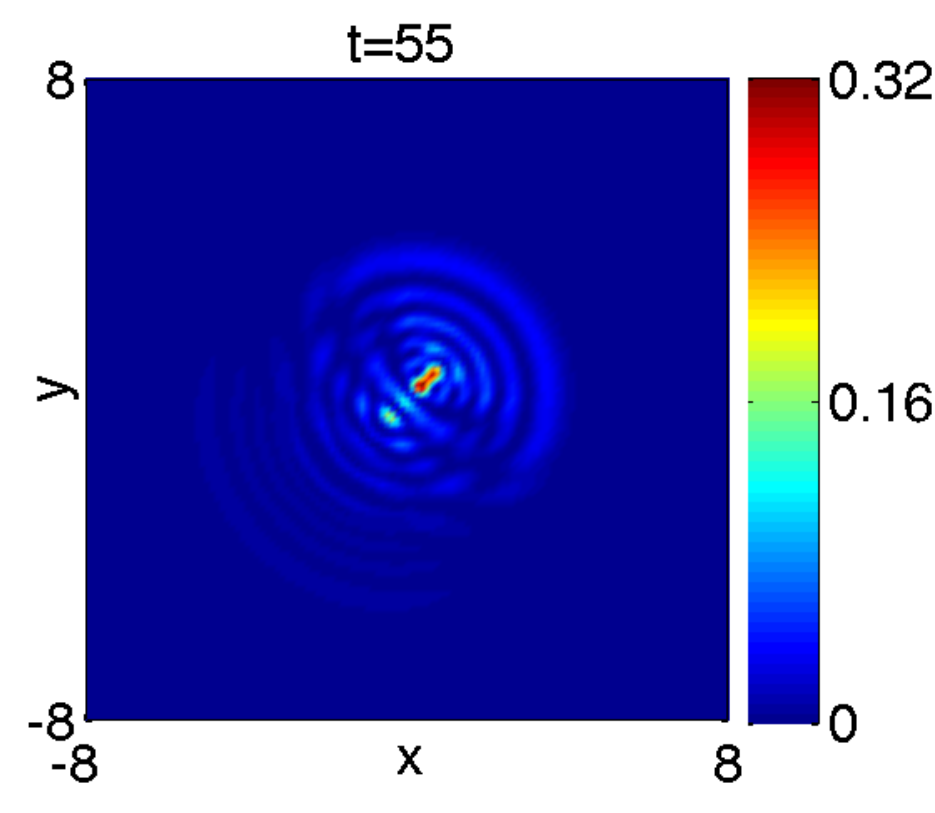,height=3.cm,width=3.3cm,angle=0}
}
\centerline{(c)
\psfig{figure=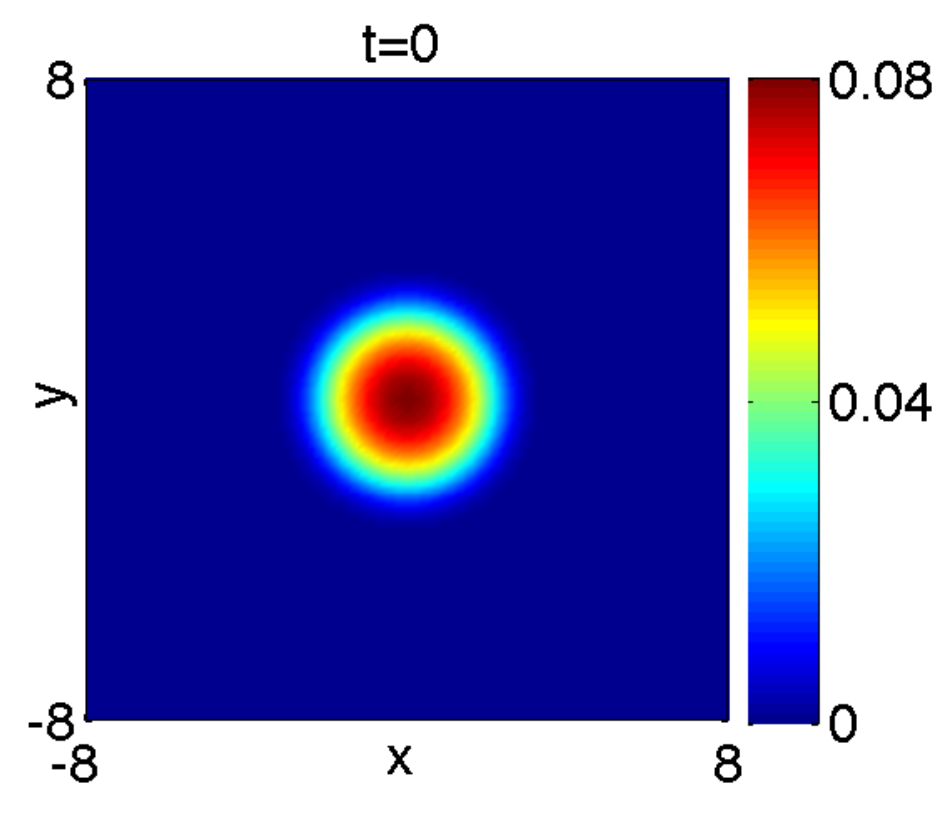,height=3.cm,width=3.3cm,angle=0}
\psfig{figure=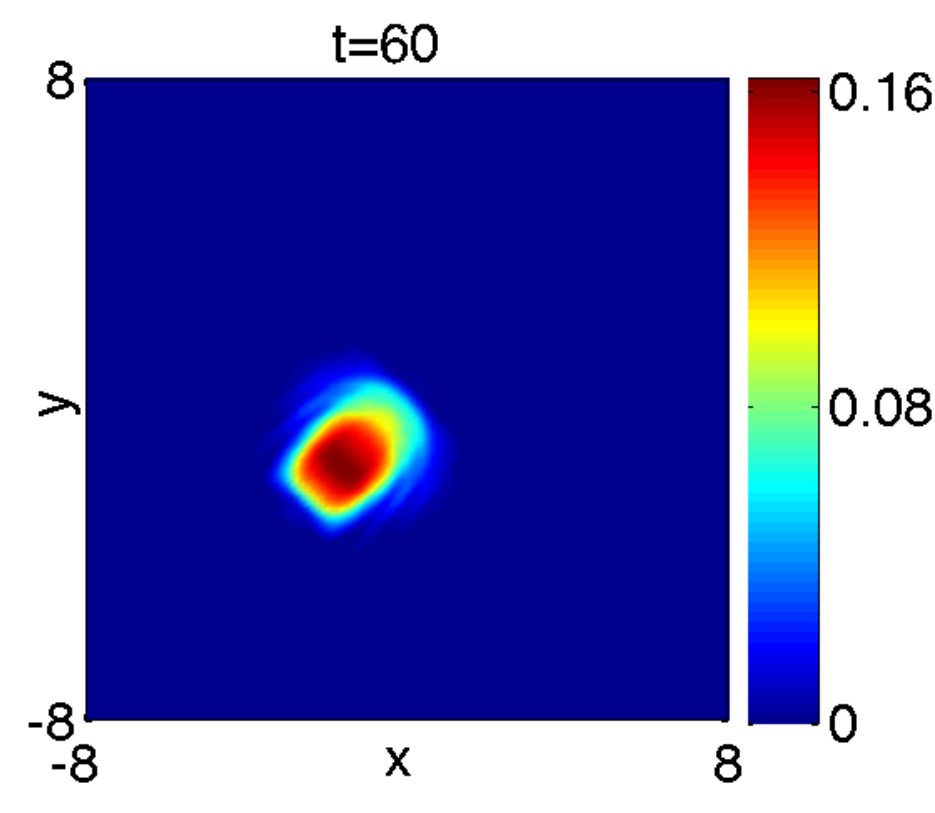,height=3cm,width=3.3cm,angle=0}
\psfig{figure=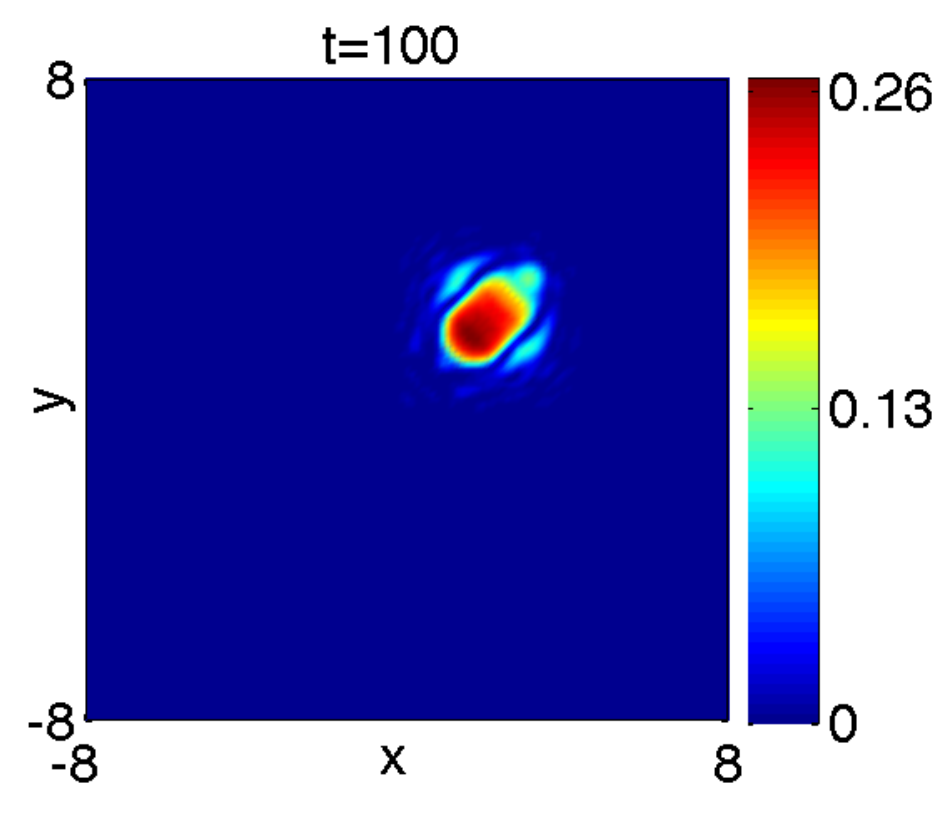,height=3.cm,width=3.3cm,angle=0}
\psfig{figure=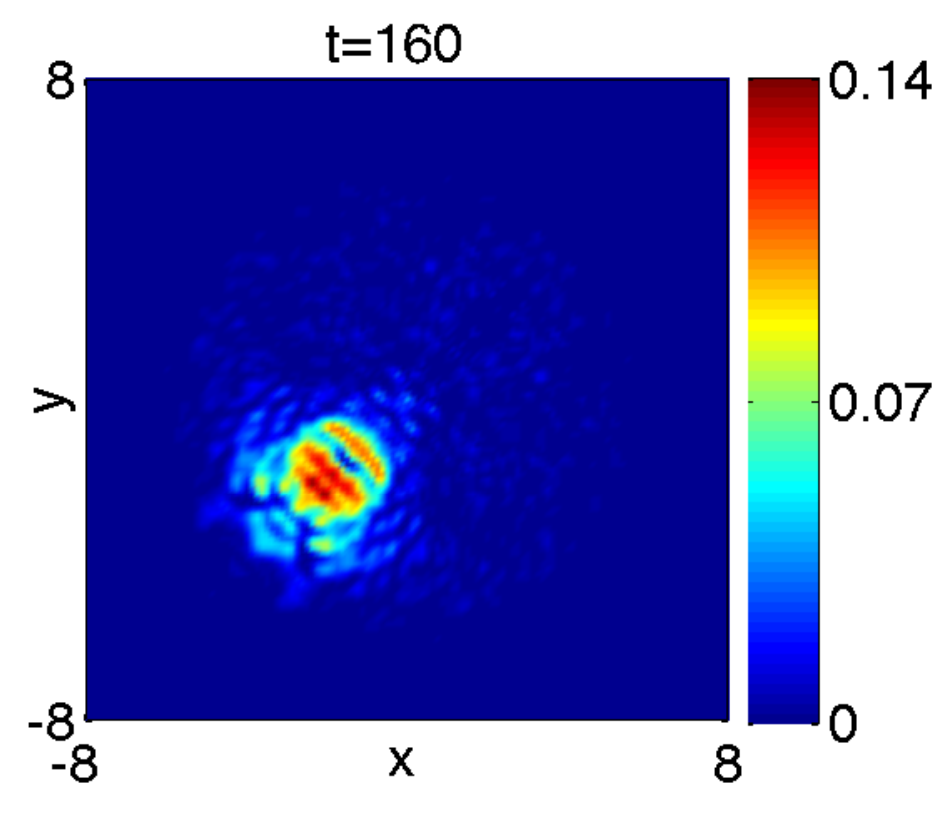,height=3.cm,width=3.3cm,angle=0}
\psfig{figure=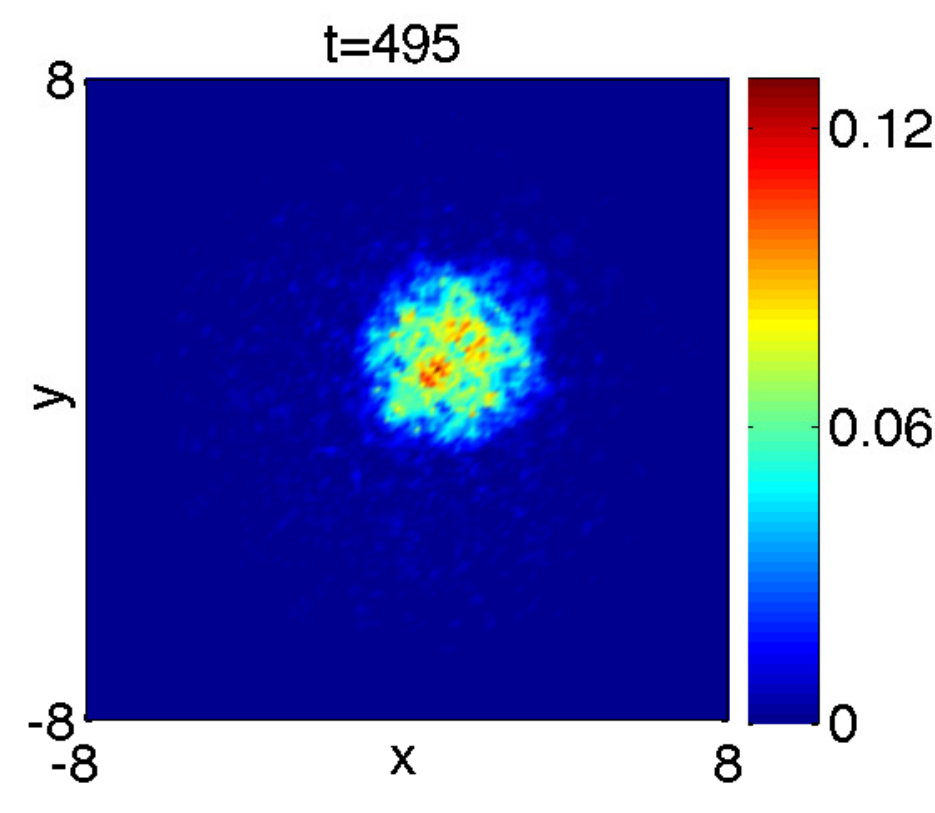,height=3.cm,width=3.3cm,angle=0}
}
\centerline{(d)
\psfig{figure=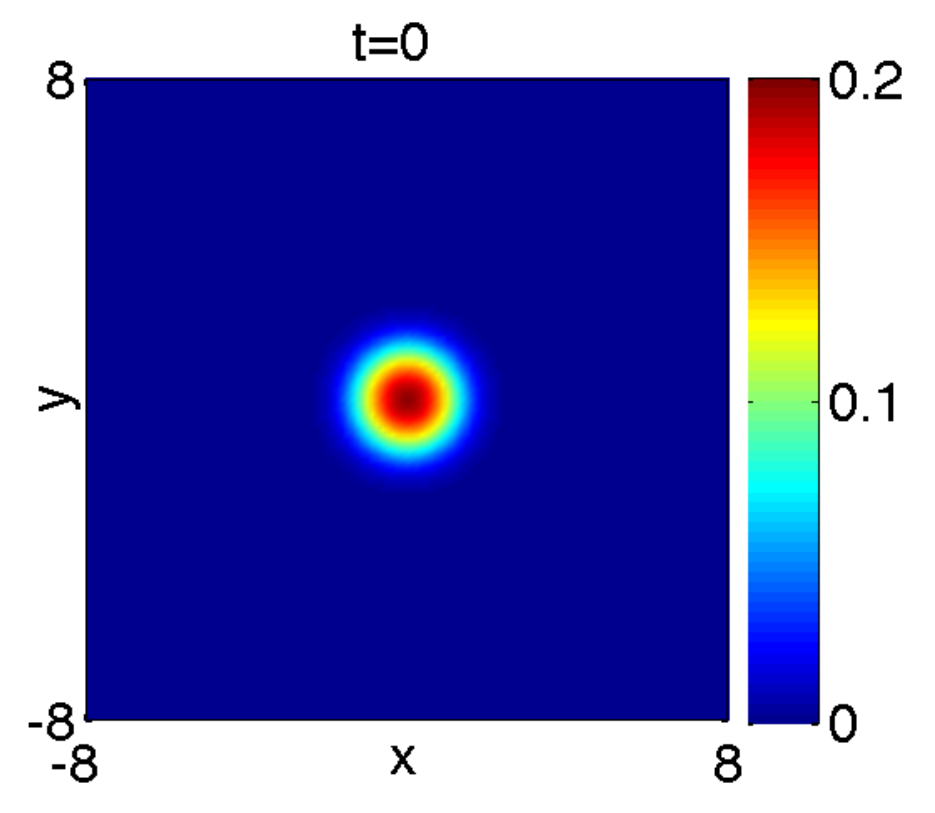,height=3.cm,width=3.3cm,angle=0}
\psfig{figure=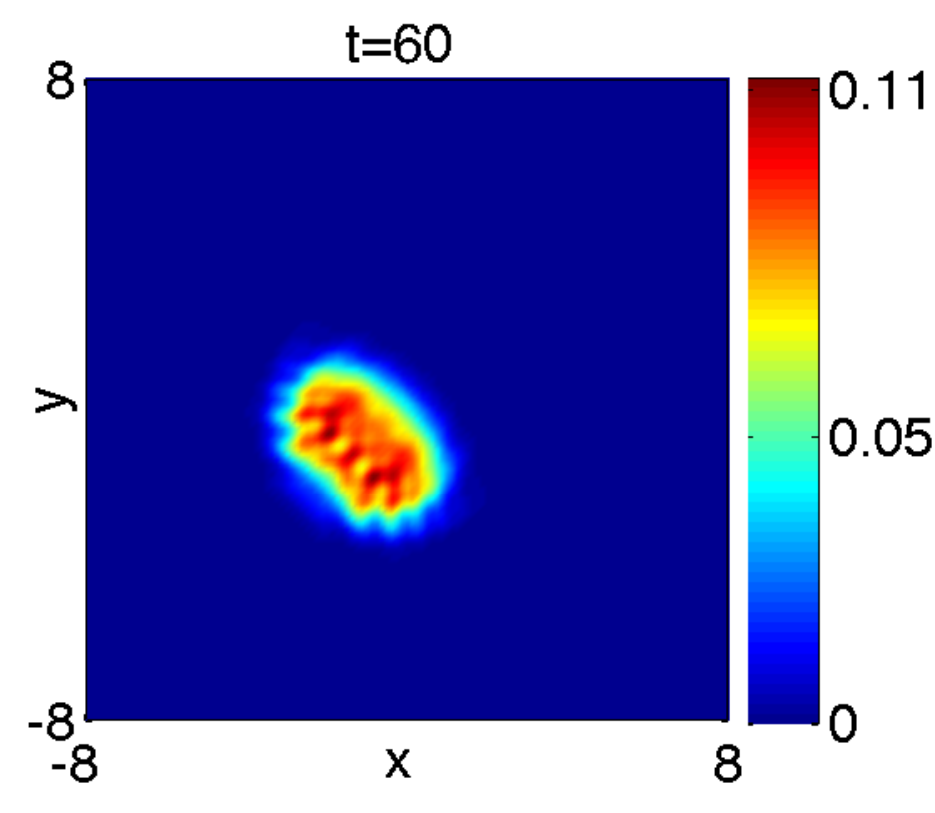,height=3cm,width=3.3cm,angle=0}
\psfig{figure=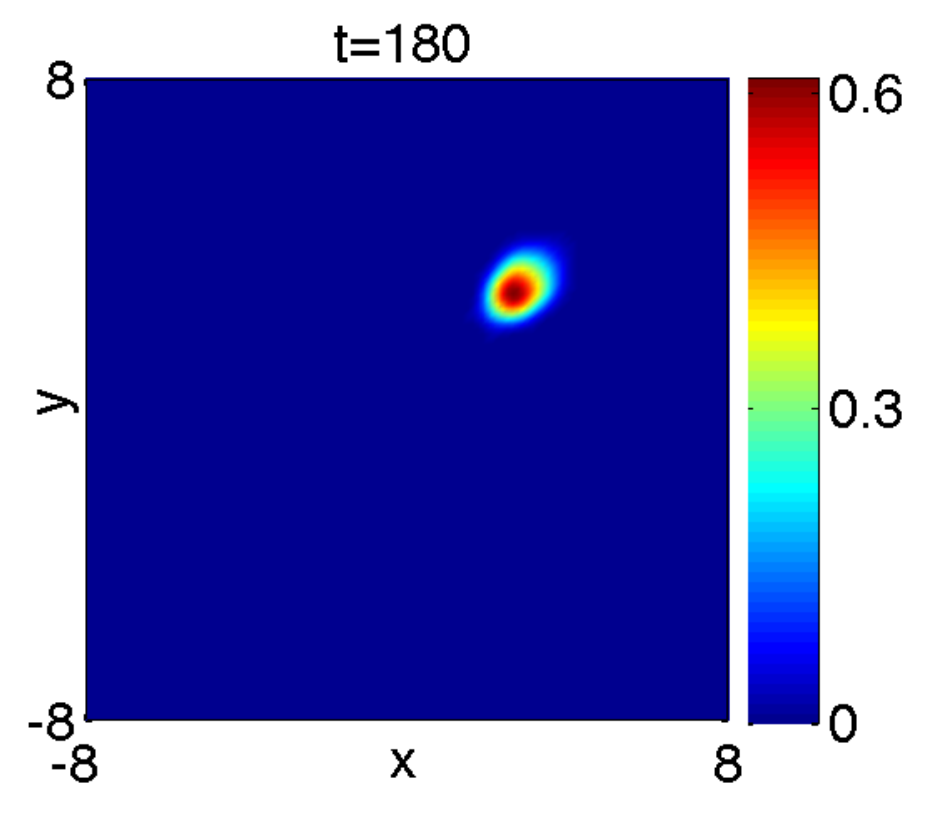,height=3.cm,width=3.3cm,angle=0}
\psfig{figure=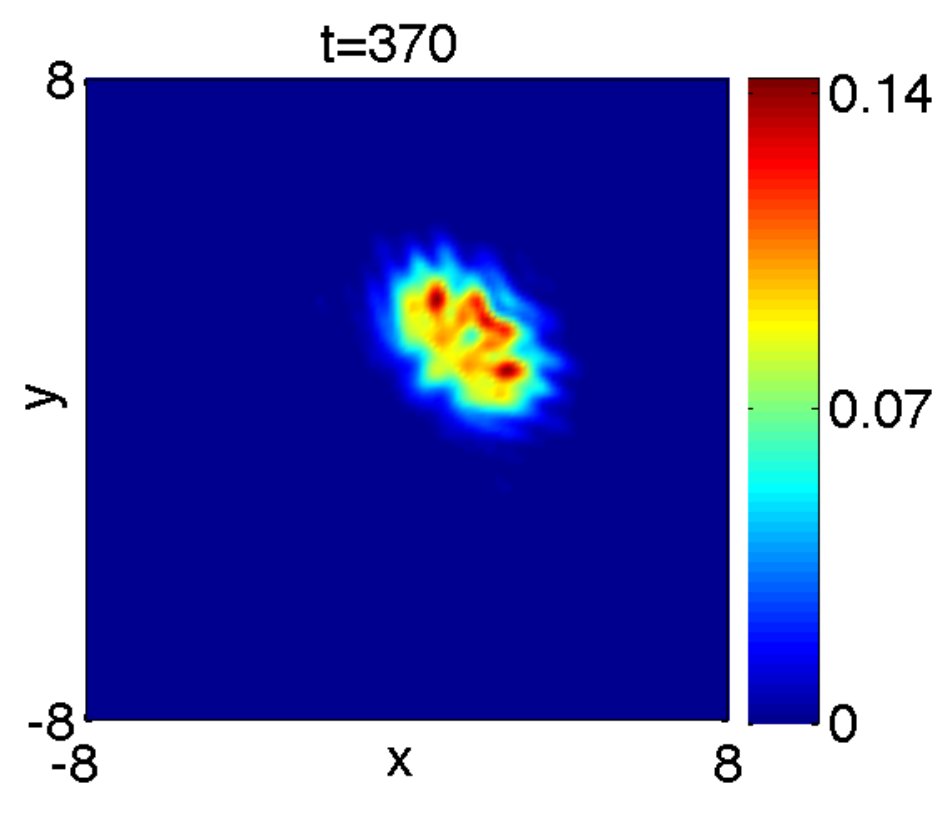,height=3.cm,width=3.3cm,angle=0}
\psfig{figure=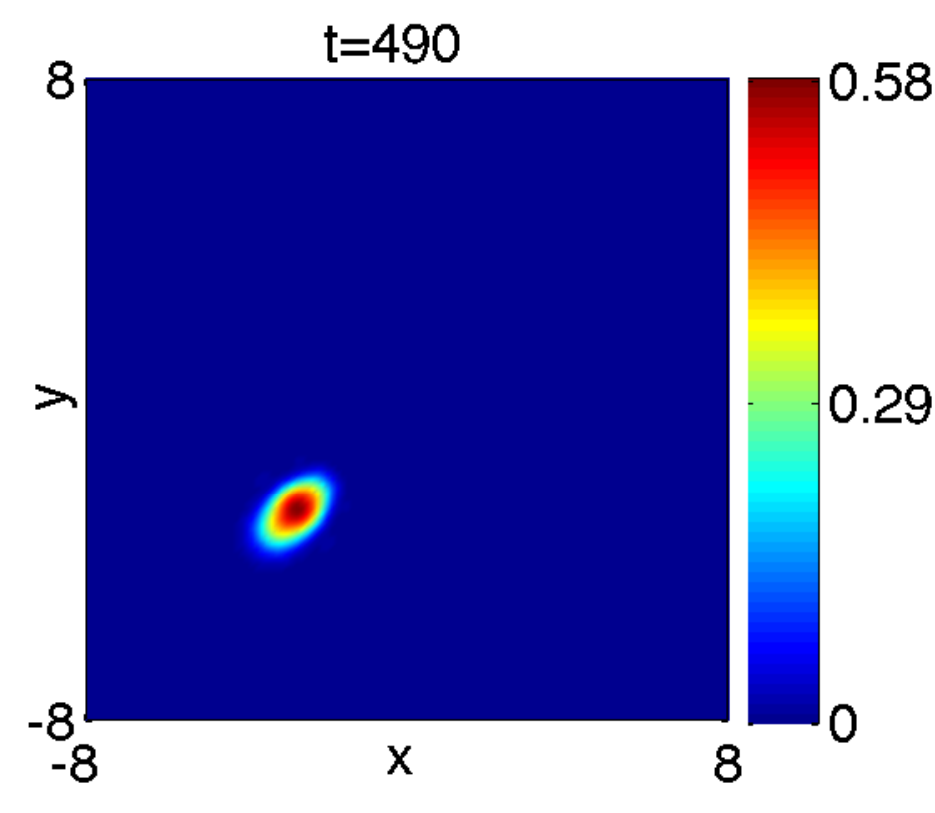,height=3.cm,width=3.3cm,angle=0}
} 
\caption{Contour plots of the density $|\psi(\bx,t)|^2$ at different times in example  \ref{dy:exmp:2} 
for cases I to  IV (from top to bottom). Here, we consider  $s=0.75$ (subdispersion).}
\label{fig:ex7_con}
\end{figure}

Figure \ref{fig:ex7} shows the  
dynamics of the mass, energy, centre of mass, condensate widths, 
while Figure \ref{fig:ex7_con} shows the contour  plot of the density $|\psi(\bx,t)|^2$ at different times.
Similarly to  Example \ref{dy:exmp:1},  we can see that
(i) For the FNLSE, the density profile no longer retains its initial shape as in the classical NLSE. 
The density profile also oscillates around the center of the  trap and decoherence emerges.
(ii) The dynamics of the wave function depends crucially on the initial shift $\bx_0$.
If the initial shift is small, the initial shape is changed slightly, i.e. the decoherence is small 
(cf. Fig. \ref{fig:ex7_con} (a)), while for large shifts, the decoherence appears very quickly. 
Turbulence and chaotic dynamics might also occur for a large $\bx_0$ in the linear FSE (cf. Fig. \ref{fig:ex7_con} (b)).
(iii) Both the short- and long-range  nonlinear interactions can 
reduce and/or  delay  the emergence of decoherence and suppress the wave function from
chaotic dynamics. Turbulence emerges in the FNLSE with pure local nonlinearity   (see Fig. \ref{fig:ex7_con} (c)), while 
the decoherence is weaker in the FNLSE with pure nonlocal nonlinearity. The density profile would actually oscillate like a breather
(cf. Fig. \ref{fig:ex7_con} (d)).
It would also be interesting to investigate the decoherence and turbulence properties in the superdispersion case
 $s>1$ and analyze how they are affected through a rotation effect. This will be analyzed in future research.
Our results are in accordance with those showed in \cite{KZ2014}.

%
%


%
%
%

%

\section{Conclusion}

In this paper, we proposed efficient and robust numerical methods for computing the ground states and dynamics of the 
FNLSE equation with an angular momentum and nonlocal interaction potentials.
Existence and non-existence of the ground states were presented and dynamical laws for the mass, energy, angular momentum and center of mass were obtained.

We then studied the ground states and dynamics of the FNLSE numerically. It was found that the 
fractional order $s$ affects both the ground states and dynamics in a significant way.  
The ground states become more peaked as  $s<1$ tends smaller, corresponding here
to subdispersion. For the superdispersion case, i.e. $s>1$, the creation of a giant vortex can be observed 
for a fast rotating system, which is totally different from the behavior of the classical GPE. 
Critical values of the rotating frequencies to create the first vortex solution are numerically found to depend on $s$.  
For the dynamics, decoherence as well as turbulence were observed in the FNLSE when an initial data is prepared
from a ground state with imprinted phase shift and/or position shift.  It is shown that the 
smaller the fractional 
exponent $s$ is, the easier the decoherence emerges. The larger the initial shift is, the easier the turbulence and
chaotic dynamics arise.  Furthermore, the presence of   repulsive  nonlinearities, both local and nonlocal,
can suppress the ``peaking''  effects of the ground states and the decoherence/turbulence observed in the dynamics. 

It is worthwhile to remark that the ground states of the FNLSE decay only 
algebraically as $|\bx|\rightarrow\infty$ when the external potential $V(\bx)$ is bounded \cite{FL2013} and $\beta<0$.
A very large computational domain is necessary for both the ground state computation and the dynamics  \cite{KSM2014}.
It would be interesting and crucial to derive a fractional version of the free boundary conditions such as 
the transient BC, absorbing BC and also the PML \cite{AABES2008} for the FNLSE.

Finally, let us emphasize that the  time and space  fractional NLSE, for $0< \gamma <1$,
\bea
\label{SFSEgamma}
&&i\p_t^{\gamma}\psi(\bx,t)=\left[\fl{1}{2}\left(-\nabla^2+m^2\right)^{s}
+V(\bx)+\beta|\psi(\bx,t)|^2+\lambda\Phi(\bx,t)
-\Omega L_z
\right]\psi(\bx,t),\\
\label{NonLocalPot2}
&& \Phi(\bx,t)=\mathcal{U}\ast|\psi(\bx,t)|^2,
\qquad \bx\in{\Bbb R}^d,\ t>0, \  d\ge2.
\eea
is also very interesting for some applications \cite{TimeFractionalErtik,Garrappa2015115,Laskin0,Laskin3, Naber2004,WX2007}. The next step of our work would consist in analyzing 
efficient and accurate numerical methods for solving FNLSEs both in space and time and understand their behavior and
properties.


\section*{Acknowledgements}

We acknowledge the support from the  ANR project BECASIM
 ANR-12-MONU-0007-02 (X. Antoine and Q. Tang), the ANR-FWF Project Lodiquas ANR-11-IS01-0003,  the
ANR project Moonrise ANR-14-CE23-0007-01 and the Natural Science Foundation of China grants  11261065, 91430103 and 11471050 (Y. Zhang).
We are also grateful  to Prof. Weizhu Bao and Dr. Yongyong Cai for valuable suggestions.


\appendix
\setcounter{equation}{0}  

\section{Proof of  lemma \ref{lawsDyn_AME}}
\label{Appendix_A}
\setcounter{equation}{0}
\renewcommand{\theequation}{A.\arabic{equation}}
\noindent
Let us introduce $\bk=\big(k_x, k_y, k_z \big)^T$ if $d=3$, and $\bk=\big(k_x, k_y\big)^T$ if $d=2$. 
 Let $J_z=y\p_x-x\p_y$ and $\widehat{J}_{z}=k_y\p_{k_x}-k_x\p_{k_y}$. Then, we have $\widehat{J_z \psi} =\widehat{J_z}
 \,\widehat{\psi} $.
Differentiating (\ref{AME}), noticing  (\ref{SFSE}) and using the Plancherel's formula, we have
\bea\label{apd-1}
\fl{d}{d t}\langle L_z \rangle (t)
\nn
&=& \langle L_z\psi_t, \psi \rangle  + \langle L_z\psi, \psi_t \rangle 
	=  -\langle i \psi_t, J_z \psi \rangle  - \langle J_z\psi, i\psi_t \rangle \\
&=& \frac{1}{(2\pi)^d} \left\{-\langle\, \fl{1}{2}\big( |\bk|^2+m^2\big)^s \widehat{\psi}+\widehat{\mathcal{V}\psi},  \widehat{J}_{z} \widehat{\psi} \,\rangle        
-\langle\, \widehat{J}_{z} \widehat{\psi}, \fl{1}{2}\big( |\bk|^2+m^2\big)^s \widehat{\psi}+\widehat{\mathcal{V}\psi} \,\rangle     \right\},                
\eea
with $\mathcal{V}\psi: = V\psi + \lambda \Phi \,\psi $. The rotation and local nonlinear terms cancel. We omit both for brevity.

By integrating the above equation by parts, we have
\bea
\fl{d}{d t}\langle L_z \rangle (t)
\nn	
&=&\frac{1}{(2\pi)^d} \left\{
-\langle\, \fl{1}{2}\big( |\bk|^2+m^2\big)^s \widehat{\psi}+\widehat{\mathcal{V}\psi},  \widehat{J}_{z} \widehat{\psi} \,\rangle
	+ \langle\, \widehat{\psi}, \widehat{J}_{z}\left(\fl{1}{2}\big( |\bk|^2+m^2\big)^s \widehat{\psi} + \widehat{\mathcal{V}\psi} \right)  \,\rangle   \right\}\\	
\nn	
&=&	\frac{1}{(2\pi)^d} \left\{ \langle\, \widehat{\psi}, \widehat{J}_{z}(\widehat{\mathcal{V}\psi})  \,\rangle 
	-\langle\, \widehat{\mathcal{V}\psi},  \widehat{J}_{z} \widehat{\psi} \,\rangle\right\} 
	=\frac{1}{(2\pi)^d} \left\{
       -\langle\, \widehat{J}_{z}\widehat{\psi}, \widehat{\mathcal{V}\psi}  \,\rangle 	
       	 -\langle\, \widehat{\mathcal{V}\psi},  \widehat{J}_{z} \widehat{\psi} \,\rangle  \right\}  \\	
&=& -\langle\, J_{z}\psi, \mathcal{V}\psi  \,\rangle  	 -\langle\, \mathcal{V}\psi,  J_{z} \psi \,\rangle
	 =\langle\,|\psi|^2 , \, J_{z}\mathcal{V}\,\rangle 
= 	\int_{\mathbb{R}^d} |\psi|^2 (y\p_x-x\p_y)\Big(V(\bx)+\lambda \Phi(\bx,t) \Big)  d\bx.  \qquad\quad
\eea
Therefore,  by adapting the polar/cylindrical coordinates transformation in 2D/3D and noticing 
$y\p_x-x\p_y=-\p_{\theta}$, one can obtain 
\be
I_1=:\int_{\mathbb{R}^d} |\psi|^2 (y\p_x-x\p_y)V(\bx)d\bx=0,
\ee
provide that $V(\bx)$ is radially/cylindrically symmetric in 2D/3D. Now that 
\be
\label{appA}
I_2=: \int_{\mathbb{R}^d} |\psi|^2 (y\p_x-x\p_y)\Phi(\bx,t) d\bx=\frac{1}{(2\pi)^d} \langle\widehat{|\psi|^2},\; \widehat{J}_{z}\widehat{\Phi} \rangle
=\int_{\mathbb{R}^d} \widehat{\mathcal{U}}(\bk)\widehat{|\psi|^2} (k_y\p_{k_x}-k_x\p_{k_y})\widehat{|\psi|^2} d\bk,
\ee
applying the polar/cylindrical coordinates transformation in 2D/3D in the Fourier space, it is easily to get $I_2=0$ if 
$ \widehat{\mathcal{U}}(\bk)$ in (\ref{Kernel}) is chosen as the Coulomb-type interaction or DDI with $\bn=(0, 0, 1)^T.$

\hfill $\square$

\section{Proof of  lemma \ref{lawsDyn_COM}}
\label{Appendix_B}
\setcounter{equation}{0}
\renewcommand{\theequation}{B.\arabic{equation}}

\noindent
{\bf Step 1:}
By differentiating (\ref{CoM}) and noticing  (\ref{SFSE}), we have
\bea
\dot{ \bx}_c(t)  
\nn	
&= &\frac{d}{dt} \langle \bx \psi,\psi \rangle = \frac{1}{i}\langle \bx \,i\psi_t,\psi \rangle 
	+ i \langle \bx \psi,i\psi_t \rangle =  i  \left[  \langle \bx \psi,i\psi_t \rangle -   \langle \bx\, i\psi_t, \psi \rangle  \right]\\
&=& \frac{i}{2} \left[ \langle \bx \psi,(-\Delta + m^2)^{s} \psi  \rangle -   
	\langle \bx\,(-\Delta + m^2)^{s} \psi, \psi \rangle  \right]
	-\Og \left[ \langle \bx J_z\psi, \psi \rangle +  \langle  \psi, \bx J_z\psi \rangle \right].
\eea
An integration by parts and an application of Plancherel's formula lead to 
\bea
\nn
\dot{ \bx}_c(t)  
&=&  \frac{i}{2} \frac{1}{(2\pi)^d} \left\{\langle i \nabla_\bk \widehat{\psi},
	(|\bk|^2+m^2)^{s} \widehat{\psi}  \rangle -   
	\langle i \nabla_\bk[(|\bk|^2+m^2)^{s}\widehat{\psi}(\bk)] ,  \widehat{\psi}\rangle  \right\}
	+\Og  \langle \psi J_z\bx, \psi \rangle \\ 
&=& \frac{s}{(2\pi)^d}  \, \langle (|\bk|^2+m^2)^{s-1}  \bk \widehat{\psi},\widehat{\psi}\rangle +\Og J \bx_c.
\label{PCoM2}
\eea
Note that (\ref{PCoM2}) is well-defined for  $\forall\,s>0$.  
If $s=1$,  (\ref{PCoM2}) yields
\be
\dot{ \bx}_c(t)- \Og J \bx_c= \frac{1}{(2\pi)^d}\langle \bk\,\widehat{\psi}, 
\widehat{\psi} \rangle=i \langle \psi, \nabla\psi \rangle=i \langle G\ast\psi, \nabla\psi \rangle,
\ee
with $G(\bx)=\delta(\bx)$.  
If $0<s<1$,  we have
\be
\Big(|\bk|^2+m^2 \Big)^{s-1}=\fl{1}{c_s}\int_{0}^{\infty}\lambda^{-s} e^{-\pi(|\bk|^2+m^2)\lambda}d \lambda,\qquad {\rm with}\qquad c_s=\Gamma\big(1-s\big)/\pi^{1-s}.
\ee
Hence, one gets
\bea
\nn
&&s\,\mathcal{F}^{-1}\Big( \big(|\bk|^2+m^2 \big)^{s-1}\widehat{\psi} \Big)
 =s\,\frac{1}{(2\pi)^d} \int_{\Bbb R^d} \Big(|\bk|^2+m^2 \Big)^{s-1}
 \widehat{\psi}\,e^{ i \bk\cdot\bx} d\bk\\
 \nn
&&=\fl{s}{c_s}  \int_{0}^{\infty}  \lambda^{-s}\, e^{-\pi \lambda m^2} 
	\left[ \frac{1}{(2\pi)^d} \int_{\Bbb R^d}  \widehat{\psi}
	\,e^{-\pi \lambda |\bk|^2}\,e^{i \bk\cdot\bx}\, d\bk \right] d\lambda
 = \fl{s}{c_s}  \int_{0}^{\infty}  \lambda^{-s}\, e^{-\pi \lambda m^2} 
	\left[ \psi \ast \mathcal{F}^{-1}\left(e^{-\pi \lambda |\bk|^2}\right)\right] d\lambda   \\
\nn
&&= \fl{s}{c_s (2\pi)^d}  \int_{0}^{\infty} \lambda^{-s}\, e^{-\pi \lambda m^2} 
	\int_{\Bbb R^d} \lambda^{-\fl{d}{2}}e^{-\fl{|\bx-\by|^2}{4\pi \lambda}} \psi(\by) d\by d\lambda \\	
\label{green_frac}	 
&&=\fl{s}{c_s (2\pi)^d}  \int_{\Bbb R^d} \psi(\by) \left[  \int_{0}^{\infty} \lambda^{-\fl{2s+d}{2}}
	 e^{-\pi \lambda m^2} \,e^{-\fl{|\bx-\by|^2}{4\pi\lambda}} d\lambda  \right] d\by 
=: \Big( G\ast \psi\Big) (\bx),	 
\eea
with 
\be\label{Gx}
G(\bx)
= \fl{s}{c_s (2\pi)^d} \int_{0}^{\infty} \lambda^{-\fl{2s+d}{2}}
	 e^{-\pi \lambda m^2} \,e^{-\fl{|\bx|^2}{4\pi\lambda}} d\lambda 
=\fl{2^{s-d/2}s}{\Gamma\big(1-s \big)\, \pi^{d/2}}
\left(\fl{m}{|\bx|}\right)^{\fl{d}{2}+s-1}K_{\fl{d}{2}+s-1}\Big(m |\bx|\Big),
\ee
where  $K_{v}(z)$ is the modified Bessel function of the second-kind and order $v$ defined by (\ref{Bes2}).
Finally, we obtain
\bea
\dot{ \bx}_c(t)- \Og J \bx_c
=\frac{s}{(2\pi)^d} \langle(|\bk|^2+m^2)^{s-1}  \bk \widehat{\psi},\widehat{\psi}\  \rangle
	=  \Big\langle s\mathcal{F}^{-1}\big(|\bk|^2+m^2)^{s-1}\widehat{\psi}\,\big), 	
	\mathcal{F}^{-1}\big( \bk\widehat{\psi}\, \big)  \Big\rangle  
= i\big\langle (G\ast\psi), \nabla\psi \big\rangle.\;\;
\eea 

\medskip

\noindent
{\bf Step 2:}
Let us consider the second-order derivative of (\ref{PCoM2}). By (\ref{green_frac}), we have 
\bea
\nn 
\ddot{ \bx}_c(t)- \Og J \dot{\bx}_c   
&=& 2s\, {\rm Re} \Big( C_d \langle\, \bk (|\bk|^2+m^2)^{s-1}   \widehat{\psi}_t,\widehat{\psi}\, \rangle\Big)  
=  2s\,  {\rm Im}\Big( C_d\langle    
	\, i \widehat{\psi}_t,\bk (|\bk|^2+m^2)^{s-1} \widehat{\psi}\,  \rangle \Big)   \\
\nn 
&=& 2\, {\rm Im}  \Big( C_d \langle s \mathcal{F}^{-1}\Big( (|\bk|^2+m^2)^{s-1} 
	\widehat{ \mathcal{V} \psi }	\Big) , \, \mathcal{F}^{-1}\big(\bk \widehat{\psi}\,\big)\,   \rangle\Big)
	+ s\,\Og\,\Big[C_d \langle\, \widehat{\psi},  \bk(|\bk|^2+m^2)^{s-1} \widehat{J}_{z}\widehat{\psi}  \,  \rangle  \\
\nn		
&&+C_d \langle\, \widehat{\psi},  \widehat{\psi}\, \widehat{J}_{z}\big( \bk(|\bk|^2+m^2)^{s-1}\big)\,  \rangle 
      -C_d \langle\, \bk(|\bk|^2+m^2)^{s-1} \widehat{\psi},  \widehat{J}_{z}\widehat{\psi}    \, \rangle \Big]\\
 \nn     
&=& 2\, {\rm Re}\Big(  \big \langle G\ast ( \mathcal{V} \psi),   \nabla\psi \big \rangle\Big)
	+ s\,\Og\,C_d \langle\, (|\bk|^2+m^2)^{s-1}\widehat{\psi},  \widehat{\psi}\, \widehat{J}_{z}\bk\,  \rangle \\
 \nn     
&=& 	2\, {\rm Re}\Big(  \big \langle G\ast ( \mathcal{V} \psi),   \nabla\psi \big \rangle\Big)
	+\Og J \big( \dot{\bx}_c-\Og J\bx_c \big).
\eea
Hence, we obtain
\be
\ddot{ \bx}_c(t)- 2\Og J \dot{\bx}_c   +\Og^2 J^2\bx_c=2\, {\rm Re}\Big(  \big \langle G\ast (\mathcal{V} \psi), \nabla \psi\rangle \Big),
\ee
ending hence the proof.
\hfill $\square$


\end{document}